\def\useAlgorithmic{0}

\documentclass[11pt]{article}

\usepackage{fullpage,graphicx,authblk,ifthen}
\usepackage{titlesec}
\titlespacing*{\subsection}{0pt}{2.25ex plus 1ex minus .2ex}{1ex plus .2ex}
\titlespacing*{\paragraph}{0pt}{2.25ex plus 1ex minus .2ex}{1em}
\usepackage[authoryear]{natbib}
\usepackage{amsmath,amsfonts,amsthm,amssymb}
\newtheorem{theorem}{Theorem}
\newtheorem{proposition}{Proposition}
\newtheorem{corollary}{Corollary}
\newcounter{def}
\newtheorem{definition}[def]{Definition}
\newcommand{\argmin}{\operatorname{argmin}}
\newcommand{\argmax}{\operatorname{argmax}}
\newtheorem{assumption}{Assumption}
\newcommand{\myproof}[1]{\begin{proof}#1\end{proof}}

\usepackage{scalerel}
\ifthenelse{\useAlgorithmic = 1}{
\usepackage{algorithm}
}{
\usepackage[ruled,vlined,linesnumbered]{algorithm2e}
}
\usepackage[noend]{algpseudocode}
\usepackage{hyperref}
\usepackage[table,xcdraw]{xcolor}
\usepackage{xspace}
\usepackage{enumitem}
\usepackage{caption,subcaption}
\hypersetup{
    colorlinks = true,
    citecolor={blue},
    linkcolor = {blue},
    menucolor = {black},
}
\newcommand{\mycitet}[2]{\citet{#2}}
\newcommand{\mycitep}[1]{\citep{#1}}
\newcommand{\proj}{\operatorname{proj}}
\newcommand{\conv}{\operatorname{conv}}

\newcommand{\cl}{\operatorname{cl}}
\newcommand{\inter}{\operatorname{int}}
\renewcommand{\Re}{\mathbb{R}}
\newcommand{\Z}{\mathbb{Z}}
\newcommand{\Q}{\mathbb{Q}}
\newcommand{\QQ}{\mathcal{Q}}
\newcommand{\B}{\mathbb{B}}
\newcommand{\F}{\mathcal{F}}
\newcommand{\G}{\mathcal{G}}

\newcommand{\C}{\mathcal{C}}

\renewcommand{\S}{\mathcal{S}}

\newcommand{\J}{L}
\newcommand{\RR}{\mathcal{R}}

\renewcommand{\P}{\mathcal{P}}

\newcommand{\R}{\mathcal{V}}

\newcommand{\noprint}[1]{}
\newcommand{\code}[1]{\texttt{#1}}
\renewcommand{\t}[1]{\texttt{#1}}
\newcommand{\MIBS}{\texttt{MibS}}
\newcommand{\BLIS}{\texttt{BLIS}}
\newcommand{\CHiPPS}{\texttt{CHiPPS}}
\newcommand{\SYMPHONY}{\texttt{SYMPHONY}}
\newcommand{\nonl}{\renewcommand{\nl}{\let\nl\oldnl}}
\newcommand{\midd}{\;\middle|\;}
\newcommand{\RNum}[1]{\uppercase\expandafter{\romannumeral #1\relax}}
\newcommand{\NPcomplexity}{\ensuremath{\mathsf{NP}}}
\newcommand{\Pcomplexity}{\ensuremath{\mathsf{P}}}

\newtheorem{feascon}{Feasibility Condition}

\makeatletter
\let\OldStatex\Statex
\renewcommand{\Statex}[1][3]{%
  \setlength\@tempdima{\algorithmicindent}%
  \OldStatex\hskip\dimexpr#1\@tempdima\relax}
\makeatother
\algnewcommand{\An}{\textbf{and}\xspace}
\algnewcommand{\Or}{\textbf{or}\xspace}

\graphicspath{{images/}}

\makeatletter
\def\input@path{{images}}
\makeatother

\newcommand{\MYTITLE}{Valid Inequalities for Mixed Integer Bilevel Linear Optimization
  Problems}

\newcommand{\MYABSTRACT}{Despite the success of
branch-and-cut methods for solving mixed integer bilevel linear optimization
problems (MIBLPs) in practice, there are still gaps in both the theory and
practice surrounding these methods. In the first part of this paper, we lay
out a basic theory of valid inequalities and cutting-plane methods for MIBLPs
that parallels the existing theory for mixed integer linear optimization
problems (MILPs). We provide a general scheme for classifying valid
inequalities and illustrate how the known classes of valid inequalities fit
into this categorization, as well as generalizing several existing classes.

In the second part of the paper, we assess the computational effectiveness of
these valid inequalities and discuss the myriad challenges that arise in
integrating methods of dynamically generating inequalities valid for MIBLPs
into a branch-and-cut algorithms originally designed for solving MILPs.
Although branch-and-cut methods for solving for MIBLPs are in principle
straightforward generalizations of those used for MILP, there are subtle but
important differences and there remain many unanswered questions regarding how
to suitably modify control mechanisms and other algorithmic details in order
to ensure performance in the MIBLP setting. We demonstrate that performance of
version \t{1.2} of the open-source solver \MIBS{} was substantially improved over
that of version \t{1.1} through a variety of improvements to the previous
implementation.}

\begin{document}

\title{\MYTITLE} 

  \author{Sahar Tahernejad\thanks{\texttt{satn@simcorp.com}}}
  \affil{SimCorp, Copenhagen, Denmark}
  \author{Ted K. Ralphs\thanks{\texttt{ted@lehigh.edu}}}
  \affil{Department of Industrial and System Engineering, Lehigh University,
    Bethlehem, PA}

\maketitle

\begin{abstract}
\MYABSTRACT
\end{abstract}

\section{Introduction}

Bilevel optimization (and multilevel optimization, more generally) provides a
framework for modeling and solution of optimization problems in which
decisions are made in multiple time stages by multiple (possibly competing)
decision-makers (DMs). Such optimization problems arise in a wide array of
applications, with ever more coming to light as computational tools for
solution of such problems become more widely available. For an overview of
bilevel optimization (the case in which the number of decision stages is two)
and its applications, we refer the reader to~\mycitet{Dempe}{DemFoundations02}.

The focus of this paper is on solution of \emph{mixed integer bilevel linear
optimization problems} (MIBLPs) in which some of the variables must take on
integer values. The first algorithm proposed for solution of MIBLPs was the
LP-based branch-and-bound algorithm of~\mycitet{Bard and Moore}{BarMooBranch90}.
Building on this basic framework, \mycitet{DeNegre and Ralphs}{DeNRal09}
showed that it could be successfully used as the basis for a branch-and-cut
algorithm much like those that revolutionized the solution of \emph{mixed integer
linear optimization problems} (MILPs). Just as in the MILP case, approximation
of the convex hull of feasible solutions by valid inequalities has proved to
be an effective means by which to improve bounds obtained by solving the
(typically very weak) linear optimization problem (LP) relaxation employed in
the original branch-and-bound algorithm.

The open-source solver \MIBS{}~\mycitep{MIBS1.2} was built on top of an
existing branch-and-cut framework intended for solving MILPs called
\BLIS{}~\mycitep{BLIS}. The original goal of \MIBS{} was to enable testing of
algorithmic ideas related to the various components of the branch-and-cut
approach in a ``white box'' framework. The initial version of the solver,
which targeted the pure integer case, was described by~\mycitet{DeNegre and
Ralphs}{DeNRal09}. Since then, \MIBS{} has continued to evolve and is now a
full-featured solver capable of handling general MIBLPs. The implementation of
version \t{1.1} was described in detail by~\mycitet{Tahernejad et.
al.}{TahRalDeN20}. This paper is a follow-on that discusses details of the
implementation of the newest version \t{1.2}, which features an improved control
mechanisms tuned specifically for the MIBLP setting and improved strategies
for generation of valid inequalities. 

The rest of the paper is divided into two parts. In the first part, a basic
theory of valid inequalities for MIBLPs is developed. A primary theme is that
much of the theory that has been developed for the MILP case and has been
exploited so successfully in the development of MILP solvers can be
generalized to the MIBLP case. However, while the techniques for generating
valid inequalities for MILPs have been well-systematized, with a rich and
well-developed theory underpinning them, the same cannot be said in the MIBLP
case. As far as we know, there is, for example, no generalization of Meyer's
Theorem~\mycitep{meyer74} for MIBLPs. In Section~\ref{sec:theory}, we attempt
to fill some of the existing gaps by laying out a theory of valid inequalities
paralleling the one that already exists in the MILP case. Then in
Section~\ref{sec:classes}, we propose a systematization that encompasses the
existing classes of valid inequalities and provides an overview of the known
techniques for generation of valid inequalities.

In the second part of the paper, we first discuss details of the practical
implementation of a branch-and-cut algorithm for MIBLPs in
Section~\ref{sec:implementation}, focusing on the integration of practical
methods of separating solutions to the LP relaxation from the convex hull of
feasible solutions by the dynamic generation of strong valid
inequalities~\mycitep{gomory:58,balas93,Cornuejols2008,Wolter2006}.
From a complexity-theoretic standpoint, the principle of equivalence of
optimization and separation tells us that cut generation is necessarily much
more challenging in the case of MIBLPs in the worst case. Our empirical
analysis shows that this is also true from the standpoint of practical
computation.

Finally, in Section~\ref{sec:computation}, we conclude that substantial
progress in improving effectiveness of practical algorithms has been made. We
describe computational experiments carried out with \MIBS{} \t{1.2} aimed at
assessing the overall effectiveness of this class of algorithms on solving the
instances in a set of standard benchmarks from the literature, the progress
that has been made over time, and the relative effectiveness of the various
classes of inequalities on different classes of problems. We also discuss the
importance of various parameters and how details of the cut generation
strategy and other components of the algorithm (particularly branching) impact
the effectiveness of the valid inequalities.

\section{Mixed Integer Bilevel Optimization}\label{sec:miblp} 

Bilevel optimization problems are difficult to solve, both in an empirical and
theoretical sense, even in the simple case of (continuous) bilevel linear
optimization problems (BLPs), in which the constraint and objective functions
are linear and the variables are continuous. MIBLPs are a specific class of
bilevel optimization problems in which the constraints and objective functions
must be linear, but some variables may also be required to take on integer
values. While BLPs are hard for the complexity class $\NPcomplexity$, MIBLPs
are hard for the class $\Sigma^\Pcomplexity_2$, one level higher in the
so-called polynomial time hierarchy~\mycitep{Stockmeyer76}.

\subsection{Formulation}

To formally define the class of optimization problems considered in this
study, let $x\in X$ denote the set of variables controlled by the
\emph{first-level} DM or \emph{leader} and $y\in Y$ denote the set of
variables controlled by the \emph{second-level} DM or \emph{follower}, where
$X=\Z^{r_1}_+ \times \Re^{n_1-r_1}_+$ and $Y=\Z^{r_2}_+\times
\Re^{n_2-r_2}_+$.  The general form of an MIBLP is
\begin{equation}\label{eqn:miblp}\tag{MIBLP}
  \min_{x \in X} cx + \Xi(x),
\end{equation}
where the function $\Xi$ is a \emph{risk function} that encodes the part of
the objective value of $x$ that depends on the response to $x$ in the second
level.

The function $\Xi$ may have different forms, depending on the precise
variant of the bilevel problem being solved.
In this paper, we consider deterministic MIBLPs and focus on the
so-called \emph{optimistic} case~\mycitep{LorMorWeak96}, in which
\begin{equation}\label{eqn:Xi} \tag{RF}
  \Xi(x) = \min \left\{ d^1y \midd y \in \P_1(x), y \in \argmin
  \{d^2y \midd y \in \P_2(x) \cap Y\} \right\},
\end{equation}
where
\begin{equation*}
\P_1(x) = \left\{y \in \Re_+^{n_2} \midd G^1 y \geq b^1 - A^1 x\right\}
\end{equation*}
is a parametric family of polyhedra containing points satisfying the linear
constraints of the first-level problem with respect to a given
$x \in \Re^{n_1}$ and
\begin{equation*}\label{eqn:lowerProblem}
\P_2(x) = \left\{y\in \Re_+^{n_2}\midd G^2y\geq b^2 - A^2 x \right\}
\end{equation*}
is a second parametric family of polyhedra containing points satisfying the
linear constraints of the second-level problem with respect to a given $x \in
\Re^{n_1}$. The input data is $A^1\in\Q^{m_1\times n_1}$, $G^1\in\Q^{m_1\times 
n_2}$, $b^1\in\Q^{m_1}$, $c \in \Q^{n_1}$, $d^1, d^2 \in \Q^{n_2}$,
$A^2\in\Q^{m_2\times n_1}$, $G^2\in\Q^{m_2\times n_2}$ and $b^2\in\Q^{m_2}$.
We note that the formulation allows participation of the second-level
variables in the first-level constraints. The matrix $G^1$ can be taken to be
a matrix of zeros (with appropriate dimensions) in the case in which the
first-level constraints are independent of the second-level variables.

In the optimistic case, the evaluation of $\Xi$ is a lexicographic
optimization problem, with the choice among alternative optima to the inner
optimization problem made according to the first-level objective function, as
specified formally in~\eqref{eqn:Xi}. In this case,~\eqref{eqn:miblp} can be
reformulated, in principle, as a (single-level) mathematical optimization
problem by considering the \emph{value function} of the second-level problem
instead of the risk function. This formulation is given by
\begin{equation}\label{eqn:miblp-vf}\tag{MIBLP-VF}
\min \left\{cx+d^1y\midd x\in X, y\in \P_1(x) \cap \P_2(x)\cap Y, 
d^2y\leq \phi(b^2 - A^2x)\right\},
\end{equation}
where $\phi$ represents the value function of the second-level problem and 
is defined as  
\begin{equation}\label{eqn:phi} \tag{VF}
  \phi(\beta) = \min \left\{d^2y\midd G^2y\geq \beta, y\in Y \right\} \quad
  \forall \beta \in \Re^{m_2}.
\end{equation}
The function $\phi$ yields the optimal value of the second-level problem
corresponding to a given $\beta\in\Re^{m_2}$. For a fixed $\beta\in\Re^{m_2}$,
the problem of evaluating $\phi$ is called either \emph{the follower's
problem} or the \emph{second-level problem}.

A wide variety of special cases of MIBLPs have been studied in the
literature. \emph{Interdiction problems} are one of the most important classes
among these special cases. In these problems, a subset of first-level
variables indexed by $\J = \{1, \dots, k_1\}$ ($k_1 \leq
\min\{r_1, r_2\}$), are the \emph{interdiction variables}, binary variables
that are each associated with a corresponding second-level variable whose
value must be zero if the associated interdiction variable's value is one. A
standard formulation for instances in this class of problem is
\begin{equation}\label{eqn:MIPINT}\tag{MIPINT}
\min \left\{d y\midd x \in \P^{INT}_1\cap X, x_\J \in \B^\J,
 y\in\argmax\{d y \midd y \in \P^{INT}_2(x) \cap Y\}\right\},
\end{equation}
where
\begin{equation*}
\P^{INT}_1 = \left\{x\in \Re_+^{n_1}\midd A x\geq b \right\},
\end{equation*}
\begin{equation*}
\P^{INT}_2(x) = \left\{y\in \Re_+^{n_2}\midd G y \geq g, y_\J \leq
\text{diag}(u)(e-x_\J) \right\}.
\end{equation*}
Here, $A \in\Q^{m_1\times n_1}$, $b \in\Q^{m_1}$, $d \in \Q^{n_2}$, $G \in\Q^{m_2\times n_2}$,
$g\in\Q^{m_2}$, $u\in\Re_+^\J$ represents the (variable) upper bound vector
for the second-level variables indexed by $\J$, and $e$ is the
$|\J|$-dimensional vector of ones. In Section~\ref{sec:classes}, we describe
specialized inequalities that are valid for this more general class. \MIBS{}
has specialized methods for solving~\eqref{eqn:MIPINT} and other special
classes of MIBLPs.

Finally, we introduce some additional sets that will be needed for discussing
the relaxations used in the remainder of the paper.
Dropping the integrality constraints and the optimality constraint of the
second-level problem from~\eqref{eqn:miblp-vf} results in the \emph{LP
relaxation}, with feasible region
\begin{equation*}
  \P = \left\{(x,y) \in \Re_+^{n_1 + n_2} \midd y \in \P_1(x) \cap \P_2(x)
  \right\}.
\end{equation*} 
This set includes all $(x,y)\in \Re_+^{n_1 + n_2}$ that satisfy the
linear constraints of the first- and second-level problems. The subset of $\P$
containing points that also satisfy the integrality constraints in both levels
is 
\begin{equation*}
\S = \P \cap (X \times Y).
\end{equation*}
and forms the feasible region of the \emph{MILP relaxation}.

\subsection{Bilevel Feasible Region and Certificates of Infeasibility}

\begin{figure}
\begin{center}
\begin{minipage}{.45\textwidth}
\scalebox{0.6}{\input{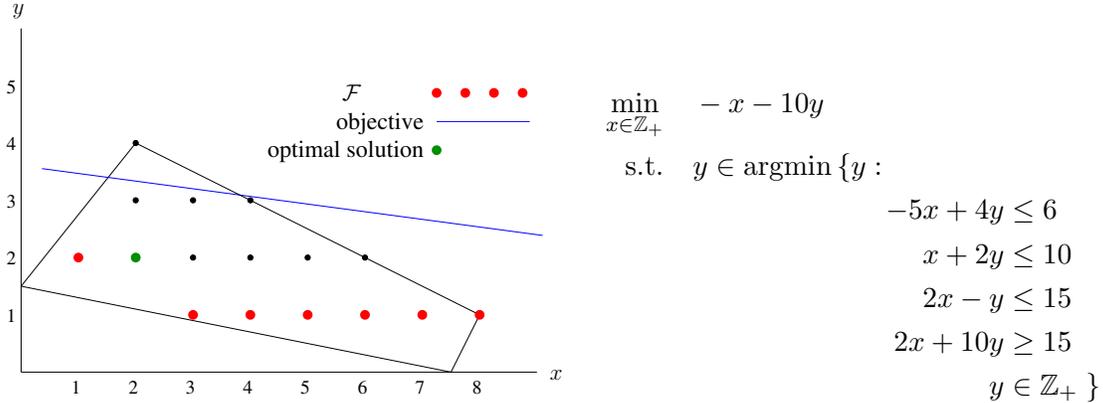}}
 \end{minipage}
\begin{minipage}{.45\textwidth}
\vskip .3in
\begin{alignat*}{3}
\min_{x\in \mathbb{Z}_+} & \quad -x - 10y\notag\\
\textrm{s.t.} &\quad y\in \operatorname{argmin}\left\{ y:\right.\notag\\
 & & -5 x + 4 y &\leq 6\notag\\
 & & x + 2y &\leq 10\label{eqn:moorece}\\
 & &  2x - y & \leq 15\notag\\
 & & 2x + 10y & \geq 15\notag\\
 & & y &\in \mathbb{Z}_+\left.\right\}\notag
\end{alignat*}
\end{minipage}
\caption{The feasible region and optimal solution of the the example from~\mycitep{MooBarMixed90}.\label{fig:mooreexample}}
\end{center}
\end{figure}

Because feasibility conditions are more complex and more difficult to verify
for a given point than in the MILP case, it is important to formally define
what is meant by the feasible region in this case. With respect to a given
$x\in\Re_+^{n_1}$, the \emph{rational reaction set} is defined as
\begin{equation*}
\RR(x) = \argmin \left\{d^2y \midd y \in \P_2(x)\cap Y\right\}
\end{equation*}
and contains all $y\in Y$ that are optimal to the second-level problem arising
from fixing the first-level variables to $x$. Note that $(x, y)$ may not be
bilevel feasible, even if $y \in \RR(x)$ if either $x \not\in X$ or
$y \not\in \P_1(x)$. As such, the
conditions for bilevel feasibility of $(x,y)\in\Re_+^{n_1 + n_2}$ can be
stated as
\begin{feascon} \label{fc:first-level}
  $x \in X$.
\end{feascon}
\begin{feascon} \label{fc:second-level}
  $y \in \P_1(x) \cap \RR(x)$.
\end{feascon}
Based on these conditions, the \emph{bilevel feasible region} is
\begin{equation*}
\F = \left\{(x,y) \in X \times Y \midd y \in \P_1(x)\cap\RR(x) \right\}.
\end{equation*}
The problem~\eqref{eqn:miblp} can thus be re-cast as the optimization problem
\begin{equation}\label{eqn:miblp-f}\tag{MIBLP-F}
  \min_{(x, y) \in \F} cx + d^1 y.
\end{equation}
The bilevel feasible region and optimal solution of the
well-known example from~\mycitep{MooBarMixed90}, shown in
Figure~\ref{fig:mooreexample}, illustrates these concepts.

In the context of MIBLPs, the infeasible solutions that arise as a result of
solving the LP relaxation can violate either Feasibility
Condition~\ref{fc:first-level} or~\ref{fc:second-level} (or both). Identifying
violations of Condition~\ref{fc:first-level} is easy, but in contrast with
MILPs, identifying violations of Condition~\ref{fc:second-level} may not be
easy. Condition~\ref{fc:second-level} can in fact be broken down into three
sub-conditions, as follows.
\begin{description}
\item[Feasibility Condition 2a] $y \in \P_1(x)$; 
\item[Feasibility Condition 2b] $y \in Y$; and
\item[Feasibility Condition 2c] $d^2 y \leq \phi(b^2 - A^2x)$. 
\end{description}
Like Condition~\ref{fc:first-level}, Conditions 2a and 2b are easy to check,
but verifying condition 2c requires solving an MILP. When condition 2c
is \emph{not} satisfied, the proof is a \emph{certificate} that generally takes one
of two forms.
\begin{definition}[Certificate of Infeasibility] \label{def:certificate-infeas}
With respect to $(\hat{x}, \hat{y}) \in \P \setminus \F$,
\begin{itemize}
\item An \emph{improving solution} is $y^* \in \P_2(\hat{x}) \cap Y$ such that
$d^2 \hat{y} > d^2 y^*$.
\item An \emph{improving direction} is $\Delta y \in Y$ such that
$d^2 \Delta y < 0$ and $y + \Delta y \in \P_2(x)$.
\end{itemize}
\end{definition}
The certificate proves the infeasibility of the given point and this
information is used directly in generating valid inequalities that separate
the point from $\F$, as described later in Section~\ref{sec:classes}. These
two certificates are theoretically equivalent for points in $\S \setminus \F$
and it is easy to convert one to the other (importantly, this is not the case
for points in $\P \setminus \S$). Despite the theoretical equivalence---each
can be generated by solving an MILP---the two approaches to proving
infeasibility are quite different in practice.

\subsection{Assumptions}

In the rest of the paper, we make the following assumptions.
\begin{assumption} \label{as:boundedness}
$\P$ is bounded.
\end{assumption}
This assumption ensures the boundedness of~\eqref{eqn:miblp}, but it is made
primarily for ease of presentation and can be straightforwardly relaxed.
\begin{assumption}
$\left\{r \in \Re_+^{n_2} \midd G^2 r \geq 0, d^2 r < 0 \right\}=\emptyset$.
\end{assumption}
This second assumption prevents the unboundedness of the second-level problem
(regardless of the first-level solution) and can be checked in a
pre-processing step.
\begin{assumption} \label{as:setJ}
All first-level variables with at least one non-zero coefficient in
the second-level problem (the \emph{linking variables}) are integer, i.e.,
\begin{equation*}
\J = \left\{i \in \{1, \dots, n_1\}\midd A^2_i \neq 0
\right\}\subseteq\left\{1,...,r_1\right\},
\end{equation*}
where $A^2_i$ represents the $i^{th}$ column of matrix $A^2$.
\end{assumption}
This third assumption guarantees that the optimal solution value of~\eqref{eqn:miblp} is
attainable whenever the optimal solution value is finite~\mycitep{VicSavJudDiscrete96}. The linking
variables are structurally important and recognition of this can
streamline the process of solving MIBLPs in several ways. 
Proposition~\ref{thm:linking-var-property} formally states the special role
played by the linking variables.
\begin{proposition}\label{thm:linking-var-property}
For the vectors $x^1$ and $x^2\in\Re_+^{n_1}$ with
$x_{\J}^1=x_{\J}^2\in\Z^{\J}$, we have
\begin{equation*}
  \phi(b^2 - A^2 x^1) = \phi(b^2 - A^2 x^2),
\end{equation*}
where $x_{\J}^1$ and $x_{\J}^2$ represent the subvector of $x^1$ and $x^2$,
respectively, corresponding to the linking variables.
\end{proposition}
This theorem states formally that in the optimistic setting, when the linking
variables are fixed, \eqref{eqn:miblp} (which is a non-linear optimization
problem in general, due to the non-linearity of the optimality constraint of
the second-level problem) becomes an MILP. This is because $\phi(b^2 - A^2x)$
is a constant whenever the linking variables are fixed.
Corollary~\ref{cor:opt-over-fixed-linking} states this more explicitly.
\begin{corollary}\label{cor:opt-over-fixed-linking}
For $\gamma \in \Z^\J$, we have
\begin{equation}\label{eqn:computeBestUB}\tag{UB}
\begin{aligned}
\min \left\{cx+d^1y \midd (x,y) \in\F, x_{\J} = \gamma \right\}  = 
\min \left\{cx+d^1y \midd (x,y) \in \S, d^2y\leq \phi(b^2 - A^{2}x), x_{\J} 
= \gamma \right\}. 
\end{aligned}
\end{equation}
\end{corollary} 
Hence, the best bilevel feasible solution $(x,y)\in\F$ with
$x_\J=\gamma\in\Z^\J$ can be obtained by solving the
problem~\eqref{eqn:computeBestUB}, which is an MILP. As already observed,
Corollary~\ref{cor:opt-over-fixed-linking} holds only in the optimistic
setting and in this setting, the corollary makes it clear that the non-linking
first-level variables can in fact all be assumed w.l.o.g to be part of the
lower-level problem. It seems unlikely that this observation has important
algorithmic implications, however.

\section{Theory of Valid Inequalities for MIBLPs \label{sec:theory}}

The theory of valid inequalities for polyhedra and other closed convex sets is
well-known and supported by deep foundations developed over several decades
(see, e.g.,~\mycitet{Gr\"otschel, Lov\'asz, and 
Schrijver}{grotschel.lovasz.schrijver:88}). Although much of the theory was
developed specifically with the case of MILPs in mind, it can be applied in
other contexts, such as MIBLPs. The main tools used in
applying this theory to MILPs are \emph{convexification} and the well-known
result of~\mycitet{Gr\"otschel, Lov\'asz, and 
Schrijver}{grotschel.lovasz.schrijver:88} that optimization over a closed
convex set is polynomially equivalent to the so-called \emph{separation
problem} associated with the set. Essentially, we optimize over
the convex hull of feasible solutions rather than the original feasible
region, thereby transforming the original non-convex problem into an
equivalent convex one.

\subsection{Convexification}

Convexification considers a conceptual reformulation of the problem obtained
by taking the convex hull of $\F$, as shown Figure~\ref{fig:mooreexample-polyhedral}.
\begin{figure}
\begin{center}
\begin{minipage}{.6\textwidth}
\begin{center}
\scalebox{0.6}{
\input{moore_both_hull_new3_pspdf}
}
\end{center}
\end{minipage}
\begin{minipage}{.3\textwidth}
\begin{align*}
\min \quad & c x + d^1 y\\
\textrm{s.t.} \quad & (x, y) \in \conv(\F)
\end{align*}
\end{minipage}
\caption{The polyhedral reformulation of the example from~\mycitep{MooBarMixed90}.\label{fig:mooreexample-polyhedral}}
\end{center}
\end{figure}
It may not be immediately obvious that convexification and separation can be
employed in the context of MIBLPs, so we first show that MIBLPs can, in
theory, be solved by a so-called \emph{cutting-plane method}. Showing this
formally involves showing that convexifying the feasible region does not
change the optimal solution value. This can be done in several steps. 

\begin{proposition}\label{prop:Fclosed}
Under Assumption~\ref{as:setJ}, $\F$ is closed.
\end{proposition}
\myproof{
Let $\S_{\gamma}$ be the feasible region of the 
problem~\eqref{eqn:computeBestUB} for $\gamma\in\Z^{\J}$. 
Under Assumption~\ref{as:setJ}, $\F$ is the union of (possibly infinite) disjoint 
sets $\S_{\gamma}$ for $\gamma\in\F_{x_\J} = \proj_{x_\J}(\F)$, i.e., 
\begin{equation}\label{eqn:f-closedness-proof1}
\F = \bigcup_{\gamma\in\F_{x_\J} }\S_\gamma.
\end{equation}
Furthermore, under Assumption~\ref{as:setJ}, for all $(\hat{x},\hat{y})\in\Re^{n_1+n_2}$, 
there is at least one neighborhood that intersects at most one of the sets 
$\S_{\gamma}$ with $\gamma\in\F_{x_\J}$. The radius $\bar{r}$ of this neighborhood 
can be defined as $0< \bar{r} < \min\left\{||x_\J-\hat{x}_\J||^2\midd x_\J\neq 
\hat{x}_\J, x_\J\in\F_{x_\J}\right\}$ (such $\bar{r}$ exists due to Assumption~\ref{as:setJ}). 
Therefore, the collection of sets 
$\S_{\gamma}$ for $\gamma\in\F_{x_\J}$ is a locally finite collection of 
$\Re^{n_1+n_2}$ under Assumption~\ref{as:setJ}. From this 
and~\eqref{eqn:f-closedness-proof1}, it follow that~\mycitep{munkres14}
\begin{equation}\label{eqn:f-closedness-proof2}
\cl(\F) = \cl\left(\bigcup_{\gamma\in\F_{x_\J} }\S_\gamma\right) = 
\bigcup_{\gamma\in\F_{x_\J} }\cl\left(\S_\gamma\right).
\end{equation}
We have $\cl\left(\S_\gamma\right) = \S_\gamma$ since this set is the feasible 
region of~\eqref{eqn:computeBestUB}, which is an MILP. Therefore, 
from~\eqref{eqn:f-closedness-proof2}, we have
\begin{equation*}
\cl(\F) = \bigcup_{\gamma\in\F_{x_\J} }\S_\gamma = \F.
\end{equation*}
Closedness of $\F$ follows. 
}
\begin{theorem}\label{thm:milp-feasregion-polyhedr}
Under Assumptions~\ref{as:boundedness} and ~\ref{as:setJ},
$\conv(\F)$ is a rational polyhedron. 
\end{theorem}
\myproof{
~\mycitet{Basu et al.}{basuetal18} showed that when input data are rational,
the closure of $\F$ can be written as the finite union of MILP-representable sets 
(with rational data). Since $\F$ is closed by Proposition~\ref{prop:Fclosed},
we have that
\begin{equation}\label{eqn:miblp-feasregion-repr}
\F = \cl\left(\F\right) = \bigcup_{i=1}^k L\left(\S_i\right),
\end{equation}
where $k$ is a scalar, $L$ denotes a linear transformation (specifically projection) and 
$\S_i$ is the feasible region of an MILP for $i=1,...,k$. 

From the fundamental theorem of integer programming~\mycitep{meyer74} and 
\eqref{eqn:miblp-feasregion-repr}, it follows that
\begin{equation}\label{eqn:polyhedr-miblp-feasregion-proof1}
\begin{aligned}
\conv\left(\F\right) = \conv\left(\bigcup_{i=1}^k L\left(\S_i\right)\right) = 
\conv\left(\bigcup_{i=1}^k\conv\left(L\left(\S_i\right) \right)\right)
\hfill\\\hfill
= \conv\left(\bigcup_{i=1}^kL\left(\conv\left(\S_i\right) \right)\right) = 
\conv\left(\bigcup_{i=1}^k \QQ_i\right),
\end{aligned}
\end{equation}
where $\QQ_i = L(\conv\left(\S_i\right))$ is a rational polyhedron for $i=1,...,k$. 

Furthermore, since $\QQ_i$ is bounded for $i=1,...k$ by
Assumption~\ref{as:boundedness}, the result follows from the fact that the
convex hull of the union of a finite number of bounded polyhedra is a
polyhedron~\mycitep{balas85,balas1998disjunctive}.
}

\begin{theorem}
Under Assumption~\ref{as:setJ}, we have that 
\begin{equation}\label{eqn:prob-opt-convhull}
\min_{(x,y)\in \F} cx+d^1y = \min_{(x,y)\in\conv(\F)} cx+d^1y.
\end{equation}
\end{theorem}
\myproof{
Since the objective function is linear, we have~\mycitep{ConCorZam14}
\begin{equation*}
\inf_{(x,y)\in \F} cx+d^1y = \inf_{(x,y)\in\conv(\F)} cx+d^1y
\end{equation*}
and the infimum of $cx+d^1y$ is attained over $\F$ if and only if it is attained 
over $\conv(\F)$. This follows the result because  the infimum of $cx+d^1y$ 
is attained over $\F$ under Assumption~\ref{as:setJ}.
}

The problem~\eqref{eqn:prob-opt-convhull} would be an LP in principle, if we
knew a complete description of $\conv(\F)$. In such case, a solution of the
MIBLP could be obtained by producing an extremal solution to this LP by, e.g.,
the simplex algorithm. Since we generally do not know a complete description
of $\conv(\F)$ and cannot construct one efficiently (this would be at least as
difficult as solving the original optimization problem), the cutting-plane
method is to employ the well-known technique of generating an approximation of
$\conv(\F)$ by solving the separation problem to generate valid inequalities,
as described next.

\subsection{Cutting-plane Method}

Although they are well-known, we first review several standard definitions and
results before describing the cutting-plane method in broad outline.
\begin{definition}
A \emph{valid inequality} for the set $\F$ is a triple
$(\alpha^x,\alpha^y, \beta)$, where $(\alpha^x,\alpha^y) \in \Q^{n_1+n_2}$ 
is the \emph{coefficient vector} and $\beta \in \Q$ is a \emph{right-hand
side}, such that
\begin{equation*}
\F\subseteq\left\{(x,y) \in \Re^{n_1\times n_2} \midd \alpha^x x + \alpha^y y \geq \beta \right\}.
\end{equation*}
\end{definition}
It is easy to see that an inequality is valid for $\F$ if and only if it is
valid for the convex hull of $\F$.
The so-called
\emph{separation problem} for $\conv(\F)$ is to either generate an 
inequality valid for $\conv(\F)$, but
violated by a given vector $(\hat{x},\hat{y})\notin\conv(\F)$ or to show that 
the vector is actually a member of $\conv(\F)$.
Formally, we define the problem as follows.
\begin{definition}
The \emph{separation problem} for $\conv(\F)$ with respect to a given
$(\hat{x},\hat{y}) \in \Re^{n_1+n_2}$ is to determine whether or not  
$(\hat{x},\hat{y}) \in \conv(\F)$ and if not, to produce an inequality 
$(\alpha^x, \alpha^y, \beta)\in \Q^{n_1+n_2+1}$
valid for $\conv(\F)$ and for which $\alpha^x \hat{x} + \alpha^y\hat{y} 
< \beta$.
\end{definition}
The process of solving an MIBLP by a standard cutting-plane method is
initiated by solving a convex relaxation of the original
problem~\eqref{eqn:miblp}. In contrast to MILPs, the problem obtained by
removing integrality constraints on all variables (at both levels) is not a
relaxation (its feasible region does not necessarily contain $\F$). On the
other hand, we have that $\F\subseteq\S\subseteq\P$, so the problem of
optimizing over either $\P$ or $\S$ is a valid relaxation. In \MIBS{}, $\P$ is
used as the feasible region of that starting relaxation and the relaxation
problem is
\begin{equation}\label{eqn:LRR}\tag{LR}
\min_{(x, y) \in \P} cx+d^1y.
\end{equation}
In the remainder of the paper, we take this as the initial relaxation. 

A valid inequality for $\F$ that is violated by a point in the feasible region
of the relaxation (typically an extreme point) is called a \emph{cutting
plane} or \emph{cut}. 
\begin{definition}
A \emph{cut} is an inequality valid for $\F$, but violated by some
$(\hat{x},\hat{y})\in\P\setminus\F$ (typically an optimal solution to a
relaxation defined with respect to $\P$).
\end{definition}
The term ``cut'' is often used interchangeably with the term ``valid
inequality,'' but they are not technically synonymous. Nevertheless, we will
use them interchangeably to some extent in the remainder of the paper in order
to be consistent with existing terminology. 

Let $(x^0, y^0)$ be an extremal optimal solution of~\eqref{eqn:LRR}. Then the
cutting-plane method consists of the following loop beginning with $t = 0$.

\begin{enumerate}

\item \label{cp-1} \textbf{\underline{Determine whether $(x^t,
y^t) \in \conv(\F)$ or not}.} Determining whether a given arbitrary point is
in $\conv(\F)$ is a $\Sigma^\Pcomplexity_2$-hard problem in general, but just
as in the MILP setting, we can exploit the fact that $(x^t,y^t)$ is an
extremal member of $\P$. Such a point must be a member of $\F$ in order to be
a member of $\conv(\F)$. Therefore, we need only determine whether it
satisfies the constraints that were relaxed. In this case, we must check
whether $(x^t,y^t)\in X\times Y$ and $y^t\in\RR(x^t)$. If so, then $(x^t,
y^t)\in \F$ (and $(x^t,y^t)$ is an optimal solution), so the solution process
is terminated. Otherwise $(x^t,y^t) \not\in \F$ and we move to
Step~\ref{cp-2}. As such, the problem of checking whether
$(x^t,y^t)\in\conv(\F)$ is in $\NPcomplexity$ when $(x^t,y^t)$ is an extreme
point of $\P$.

\item \label{cp-2} \textbf{\underline{Separate $(x^t,y^t) \not\in \F$ from
$\conv(\F)$}.} 
To do so, we construct an inequality $(\alpha^x, \alpha^y, \beta)$ that is
valid for $\conv(\F)$ but violated by $(x^t,y^t)$. This can be done by the
procedures discussed in Section~\ref{sec:classes}.

\item \label{cp-3} \textbf{\underline{Iterate}.} Add the constraint $\alpha^x
x + \alpha^y y \geq \beta$ to the relaxation and re-solve to obtain $(x^{t+1},
y^{t+1})$. Set $t = t+1$ and go to Step~\ref{cp-1}.

\end{enumerate}

As usual, we iterate until either the method terminates or until one of a
specified set of termination criteria are met, e.g., the number of iterations
exceeds some pre-defined limit. For the pure integer case ($n_1 = r_1$, $n_2 =
r_2$), a finite cutting-plane algorithm for~\ref{eqn:miblp} was described
by~\mycitet{DeNegre and Ralphs}{DeNRal09}.

Whether or not this method converges finitely depends on exactly how the valid
inequalities are generated and what properties they are guaranteed to have. In 
the case of MILPs, finite cutting-plane algorithms for the pure integer and general cases 
under mild assumptions were, respectively, given by~\mycitet{Gomory}{gomory:58} 
and~\mycitet{Del Pia and Weismantel}{delpiaWeismantel12} 
(see~\mycitep{gadeKucukyavus11} for more detailed discussion on the convergence 
of the cutting-plane algorithm). Since the feasible region $\P$ of
our assumed initial relaxation is a polyhedron, the bounding
problem~\eqref{eqn:LRR} is an LP and can be solved by standard algorithms.
Any extreme point of $\P$ is either a member of $\F$ or
is \emph{not} contained in $\conv(\F)$ and can be separated from it by an
inequality valid for $\conv(\F)$, as described above.
Hence, MIBLPs can, in principle, be solved by a cutting-plane method.

\subsection{Improving Valid Inequalities}

Traditionally, cutting-plane methods have been described theoretically as
generating only inequalities valid for the entire feasible set. In practice,
however, it is well-known that the addition of inequalities removing feasible
solutions is permitted, as long as this does not change the \emph{optimal
solution value} of the original problem. In the case of~\eqref{eqn:miblp-vf},
inequalities removing subsets of $\F$ are used routinely and we thus formally
define a notion of valid inequality that allows this.

When $(x^*,y^*) \in \F$ and we have that $cx^* + d^1y^* \leq \underset{(x,y) \in
(\G \cap \F)}{\min} cx + d^1y$ for some set $\G\subseteq\Re^{n_1+n_2}$, 
then we have that
\begin{equation*}
\min_{(x,y) \in \conv(\F)} cx + d^1y = \min\left\{cx^*+d^1y^*,
\min_{(x,y) \in \conv(\F)\setminus \G}
cx + d^1y\right\}. 
\end{equation*}
In this case, although $\G \cap\conv(\F)$ may contain a subset of the feasible 
region, it can be removed because it does not include any \emph{improving 
solutions} relative to $(x^*,y^*)$. 
Based on this discussion, the definition of valid inequality can be
generalized as follows.
\begin{definition}[Improving Valid
    Inequality]\label{def:valid-improving-inequality}
An \emph{improving valid inequality} for the set $\F$ with respect to an
incumbent $(x^*,y^*)\in\F$, is a triple
$(\alpha^x,\alpha^y,\beta)\in\Re^{n_1+n_2+1}$ such that 
\begin{equation*}
\left\{(x,y) \in \F \midd cx + d^1y < cx^* + d^1y^*  \right\} \subseteq
\left\{(x,y) \in \Re^{n_1\times n_2} \midd \alpha^x x + \alpha^y
y \geq \beta \right\}. 
\end{equation*}
\end{definition}
In the remainder of paper, the term ``valid inequality'' is taken to mean
``improving valid inequality'' unless otherwise stated.

\subsection{Methods for Construction \label{sec:gen-classes}}

In this section, we discuss general recipes for obtaining valid inequalities
for MIBLPs. Each recipe is generalized from a standard framework used to
derive valid inequalities in the MILP case. Note that these frameworks have
close connections to each other and most classes of inequalities can be
developed by application of more than one of these frameworks. Although the
general classes employed in solving generic MILPs were conceived specifically
for that purpose, the theoretical basis for many of the classes does not
actually depend on any specific properties of MILP and they can thus be easily
employed in other settings. The reader is referred to~\mycitet{Marchand et
al.}{Marchand2002},~\mycitet{Wolter}{Wolter2006},
and~\mycitet{Cornu\'ejols}{cornuejols08} for broad overviews of these
approaches and the relationships between the various classes. We describe each
of these classes here at a high level in the context of MIBLPs and then
discuss specific instantiations of these construction methods in
Section~\ref{sec:classes}.

\subsubsection{Disjunctive Procedure}

Disjunctive programming is both a modeling paradigm and a set of algorithmic
techniques introduced by~\mycitet{Balas}{balas:79} based on the concept of
what we refer to as a \emph{valid disjunction}.
\begin{definition}[Valid Disjunction]
A collection of disjoint sets $X_i\subseteq\Re^{n_1+n_2}$ for
$i=1,...,k$ represents a \emph{valid disjunction} for $\F$ if
\begin{equation*}
\F \subseteq \bigcup_{i=1}^{k}X_i.
\end{equation*}
\end{definition}
In this context, the collection $\{X_i\}_{1\leq i \leq k}$ is called
a \emph{disjunctive set}. The branch-and-cut algorithm depends crucially on
the identification of valid disjunctions that are violated by the solution to
some relaxation. The identified disjunctions are used both for branching and
cutting, two essential elements of the branch-and-cut algorithm.

Most algorithms for solving MILPs exploit in some way the fact that points in
$\P \setminus \S$ must violate one of the trivial disjunctions arising from
integrality requirements on the variables. Most algorithms for solving MIBLPs
exploit the fact that any point $(\hat{x}, \hat{y}) \in \S \setminus \F$ must
violate the more complex valid disjunction
\begin{equation}
\left(
\begin{aligned} 
& A^1 x \geq b^1 -  G^1 y^* \\
& A^2 x \geq b^2 - G^2 y^*\\
& d^2 y \leq d^2 y^*
\end{aligned}
\right)
\quad\quad \textrm{ OR } \quad\quad 
\left(
\begin{aligned} 
A^1x & \not\geq b^1 - G^1 y^* \\ \label{eq:miblp-disj}
& \textrm{OR } \\
A^2x &\not\geq b^2 - G^2 y^*,
\end{aligned}
\right)
\end{equation}
where $y^* \in \P_1(\hat{x}) \cap \RR(\hat{x})$ is an improving solution with
respect to an incumbent $\hat{y}$ (recall
Definition~\ref{def:certificate-infeas}). Note that such a $y^* \not= \hat{y}$
must exist. Although this disjunction can be written compactly using the
``$\not\geq$'' operator, there is no compact way to express it using linear
inequalities, which is an indication of what makes separation difficult in the
MIBLP case.

Just as with valid inequalities, the basic notion of valid disjunction can be
modified to allow for the possibility that the disjunction does not contain the
entire set $\F$, but possibly eliminates some solutions known to be suboptimal. 
This yields the concept of an \emph{improving valid disjunction}. 
\begin{definition}[Improving Valid Disjunction]
A collection of disjoint sets $X_i\subseteq\Re^{n_1+n_2}$ for
$i=1,...,k$ represents an \emph{improving valid disjunction}
for $\F$ with respect to an incumbent $(x^*,y^*)\in\F$ if
\begin{equation*}
\left\{(x,y) \in \F \midd cx + d^1y < cx^*+d^1y^* \right\} \subseteq
\bigcup_{i=1}^{k} X_i.
\end{equation*}
\end{definition}
In the remainder of paper, the term ``valid disjunction'' is taken to mean
``improving valid disjunction'' unless otherwise stated.  

\begin{definition}[Disjunctive Inequality]
A \emph{disjunctive (valid) inequality} for the set $\F$ with respect to
$\QQ \supseteq \F$ and a valid disjunction $\{X_i\}_{1 \leq i \leq k}$ is a
triple $(\alpha^x,\alpha^y,\beta)\in\Re^{n_1+n_2+1}$ such that
\begin{align*} 
\QQ\;\; \bigcap\;\; \left(\;\bigcup_{i=1}^{k} X_i\;\right)\subseteq
\left\{(x,y) \in \Re^{n_1\times n_2} \midd \alpha^x x
+ \alpha^yy \geq \beta \right\}.  
\end{align*}
\end{definition}
A subclass of the disjunctive inequalities arises from the disjunctions known
as \emph{split disjunctions} that have only two terms and are defined as
follows. 
\begin{definition}[Split Disjunction]~\label{def:split-disj}
Let $(\pi^x,\pi^y,\pi_0)\in\Z^{n_1+n_2+1}$ 
be such that $\pi^x_i = 0$ for $i \geq r_1 +1$ and $\pi^y_i = 0$ for $i \geq
r_2 +1$ (the coefficients of the continuous variables in both first and second levels
are zero). Then when
\begin{equation}\label{eqn:split-disjunc}
X_1 = \left\{(x,y) \in \Re^{n_1\times n_2} \midd \pi^x x + \pi^y y \leq \pi_0
- 1 \right\}\:\text{and}\: 
X_2 = \left\{(x,y) \in \Re^{n_1\times n_2} \midd \pi^x x + \pi^y
y \geq \pi_0  \right\}, 
\end{equation}
$\{X_1, X_2\}$ is a valid disjunction for $\F$ called a \emph{split
disjunction}.
\end{definition}
The validity of the above disjunction arises from the fact that the inner
product of any member of $\F$ with the coefficient vector is an integer and
thus all members of $\F$ must belong to either $X_1$ or $X_2$. 
The disjunctive inequalities derived from such a disjunction are known
as \emph{split inequalities}. 
\begin{definition}[Split Inequality] \label{def:split-ineq}
A \emph{split inequality} for the set $\F$ with respect to $\QQ \supseteq \F$
and a split disjunction $\{X_1, X_2\}$ is a disjunctive inequality with
respect to $\QQ$ and the disjunction $\{X_1, X_2\}$. 
\end{definition}

\subsubsection{Generalized Chv\'atal Procedure} The usual Chv\'atal
inequalities, as defined in the theory of MILPs (see~\mycitep{ConCorZam14}),
are a subclass of the split inequalities in which $X_1 \cap \QQ = \emptyset$
(as defined in Definition~\ref{def:split-ineq}). We introduce here a class we
refer to as \emph{generalized Chv\'atal inequalities} that differ slightly in
form, but are similar in spirit to the usual Chv\'atal inequalities.
\begin{definition}[Generalized Chv\'atal Inequality]
\label{def:generalized-chvatal}
Let a split disjunction $\{X_1, X_2\}$ for set $\F$ be given such that
$X_1 \cap \QQ = \{(\hat{x}, \hat{y})\}$ and $(\hat{x}, \hat{y}) \not\in \F$
for some given $\QQ \supseteq \F$. Then $(\pi^x,\pi^y, \pi_0)$ in~\eqref{eqn:split-disjunc} is
itself a valid inequality for the set $\F$, known as a \emph{generalized
Chv\'atal inequality} with respect to set $\QQ$. 
\end{definition}
The connection to the Chv\'atal inequalities in MILP should be clear. In the
MILP case, the usual Chv\'atal inequalities are split inequalities for which
we have a proof that $X_1$ does not contain any feasible solutions to the LP
relaxation of the MILP (and hence does not contain any solutions to the MILP
itself). We can thus easily conclude that the MILP feasible region is
contained in $X_2$ and derive the associated valid inequality. Here, we extend
this idea to allow that $X_1$ may contain a (single) solution to the
relaxation (whose feasible region we can think of as being $\QQ$), but that we
have an independent proof that the single solution to the relaxation is not
contained in the feasible region $\F$ of the MIBLP itself. We thus have a
similar proof that $\F$ is contained in the set $X_2$, yielding a valid
inequality as in the MILP case. Note that we could further extend this
definition to include cases where $X_1$ contains a set of feasible solutions to
the relaxation, all of which can be proven not to be in $\F$. However, as we
currently have no practical application of such a definition, we use the
simpler one here.

An obvious question that we address later is how to obtain a split disjunction
$(\pi^x,\pi^y, \pi_0)$ satisfying the requirements of the above definition in
a practical way. In MILP, a relevant split disjunction can be derived from the
tableau using the procedure of Gomory~\mycitep{gomory60}. One way of viewing
this procedure is that we first derive an inequality valid for the feasible
region of the LP relaxation by taking a combination of the inequalities in the
original formulation, then ensure that the coefficients of the left-hand side
satisfy the integrality requirements for a split disjunction either by
rounding or in some other fashion, and finally round the right-hand side to
obtain a valid inequality (see~\mycitet{Gomory}{gomory60}, for more details).

A very similar procedure is possible in the case of MIBLP by taking $\QQ$ in
Definition~\ref{def:generalized-chvatal} to be the LP relaxation $\P$
of~\eqref{eqn:miblp} and $(\hat{x}, \hat{y}) \in (X \times Y) \setminus \F$ to
be an extremal optimum to the LP relaxation. By taking a combination of the
constraints binding at $(\hat{x}, \hat{y})$, we can derive a split disjunction
satisfying Definition~\ref{def:generalized-chvatal}, as detailed later in
Theorems~\ref{thm:generalizedChvatalCut} and~\ref{thm:bilevelFeasibilityCut}.

\subsubsection{Intersection Inequalities}

Another well-known methodology for generating valid inequalities for MILPs
that can be generalized to the MIBLP setting is that of generating
so-called \emph{intersection inequalities}, which we refer to by their more
common name \emph{intersection cuts} (ICs). This method of constructing
inequalities was originally proposed by~\mycitet{Balas}{balas71} and has
already been extended to more general classes of optimization problems
by~\mycitet{Bienstock et al}{bienstocketal16} and was extended to MIBLPs
by~\mycitet{Fischetti et al}{FisLjuMonSinUse18,FisLjuMonSinNew17}.

The first step in constructing an IC
is to identify a so-called \emph{bilevel-free} set. The identification of such
sets underlies not only the generation of ICs, but also many other classes of
inequalities.
\begin{definition}
A \emph{bilevel-free set} (BFS) is a closed, convex set $\C \subseteq \Re^{n_1
+ n_2}$ such that $\inter(\C) \cap \F = \emptyset$.
\end{definition}
If $\C$ is a BFS, then inequalities valid for
$\conv(\overline{\inter(\C)}\cap \P)$\footnote{The notation
$\overline{\inter(\C)}$ denotes the complement of the interior of $\C$} are
also valid for $\F$. Such inequalities are called \emph{intersection cuts}
(ICs).

Given a point $(\hat{x}, \hat{y}) \not\in \F$ such that $(\hat{x}, \hat{y})$
is an extreme point of $\P$ and a BFS $\C$ containing $(\hat{x}, \hat{y})$ in
its interior, an IC violated by $(\hat{x}, \hat{y})$ can easily be generated
using the following general recipe.
\begin{enumerate}
\item Let $\R(\hat{x}, \hat{y})$ be a simplicial radial cone with vertex
$(\hat{x}, \hat{y})$ described by any set of $n_1 + n_2$ linearly
independent inequalities of $\P$ that are binding at $(\hat{x}, \hat{y})$.
\item Let the triple $(\alpha^x,\alpha^y, \beta)\in\Re^{n_1+n_2+1}$ be such
that the set $\{(x, y) \in \Re^{n_1 + n_2} \mid \alpha^x x + \alpha^y y
= \beta\}$ is the unique hyperplane containing the points of 
intersection of $\C$ with the extreme rays of $\R(\hat{x}, \hat{y})$.
\item Then $(\alpha^x,\alpha^y, \beta)$ is an inequality valid for $\F$ and
violated by $(\hat{x}, \hat{y})$. 
\end{enumerate}
Typically, the point to be separated is a basic feasible solution to an LP
relaxation and so the simplicial cone can be easily obtained by considering an
optimal basis. 

The challenge in generating ICs is primarily in the identification of the BFS.
Fortunately, a BFS is easily generated when one of the two certificates of
infeasibility from Definition~\ref{def:certificate-infeas} has already been
generated. Given $(\hat{x}, \hat{y}) \in \P$ and an associated improving
solution $y^* \in Y$, then
\begin{equation}\label{eqn:IS-BFS} \tag{IS-BFS}
\C(y^*) = \left\{(x,y)\in\Re^{n_1\times n_2}\midd 
d^2y \geq d^2 y^*, A^2 x \geq b^2 - G^2 y^* - 1 \right\}.
\end{equation}
is a BFS. Similarly, given an associated improving direction $\Delta y \in Y$,
    we have that
    \begin{equation}\label{eqn:ID-BFS} \tag{ID-BFS}
\C(\Delta y) = \left\{(x,y)\in\Re^{n_1\times n_2}\midd A^2 x+ G^2 y \geq 
b^2 - G^2 \Delta y - 1, y+\Delta y \geq -1\right\}
\end{equation}
is also a BFS.
Note that the size of the defined set $\C$ can impact the strength of the
generated IC (in general, bigger is better). Since the size is determined by the
particular certificate of infeasibility generated and the goal is to find the
largest set $\C$ possible, we must take into account how big the resulting BFS
will be when solving the problem of generating the certificate of
infeasibility. We discuss this more later. 

\subsubsection{Benders Cuts}

A final class of inequalities that are not
typically considered as a generic class in the theory of MILP, but are of
particular importance in solving MIBLPs are the Benders cuts. Benders
cuts arise from first applying a Benders reformulation, which involves
projecting out a subset of the variables. This reformulation introduces
optimality conditions that ensure the variables projected out must take optimal
values with respect to the values of the remaining variables. Valid
inequalities can then be derived by relaxing these optimality conditions. The
general theory is described in detail by~\mycitet{Bolusani and
Ralphs}{BolRal22}.

In the particular case of MIBLPs, the second-level optimality condition is
captured by the constraint
\begin{equation}
d^2 y \leq \psi(x) := \phi(b^2 - A^2 x),
\end{equation}
which can be seen as arising from a Benders-style projection operation.
Inequalities of the form
\begin{equation} \label{eqn:benders}
d^2 y \leq \bar{\psi}(x),
\end{equation}
where $\bar{\psi}: \Re^{n_1} \rightarrow \Re$ is such
that $\bar{\psi}(x) \geq \psi(x)$ for all $x \in \Re^{n_1}$ are therefore
generally called \emph{Benders cuts}. We introduce two classes of such
inequalities in the next section.

\section{Valid Inequalities for MIBLPs \label{sec:classes}}

As with MILPs, the main aim of generating a valid inequality is to remove from
the feasible region of the relaxation a solution which is not feasible for the
original problem, along with as much of the surrounding region as possible. In
the case of MILPs, all such solutions that arise in a typical cutting-plane
method violate integrality conditions and known classes primarily exploit
violations of simple disjunctions associated with single variables. In MIBLP,
the situation is more complex, as we have several different kinds of points we
are trying to separate and the methodologies for separating vary substantially
in both strength and computational expense.

In general, solutions to~\eqref{eqn:LRR} can be categorized as follows. 
\begin{enumerate}[label=C\theenumi.,ref=C\theenumi]
\item An extreme point $(\hat{x},\hat{y})$ of $\P$ such that 
$(\hat{x},\hat{y})\notin X\times Y$ and $d^2\hat{y} \leq \phi(b^2 - A^2\hat{x})$.
\label{case:goal-cut-1-a}
\item An extreme point $(\hat{x},\hat{y})$ of $\P$ such that 
$(\hat{x},\hat{y})\notin X\times Y$ and $d^2\hat{y} > \phi(b^2 - A^2\hat{x})$.
\label{case:goal-cut-1-b}
\item An extreme point $(\hat{x},\hat{y})$ of $\P$ such that $(\hat{x},\hat{y})\in 
X\times Y$ and $d^2\hat{y} > \phi(b^2 - A^2\hat{x})$.\label{case:goal-cut-2}
\end{enumerate}
Figure~\ref{fig:infeasible-points} illustrates the distinction between these
types of points.  
\begin{figure}
\begin{center}
\scalebox{0.8}{
\input{moore_both_hull_new5_pspdf}
}
\end{center}
\caption{Types of infeasible points in the example
from~\mycitep{MooBarMixed90}. \label{fig:infeasible-points}}
\end{figure}
Distinguishing between these cases is important, but not easy. Determining
when $(x, y) \not\in X \times Y$ is straightforward, but further
distinguishing between the cases~\ref{case:goal-cut-1-a}
and~\ref{case:goal-cut-1-b} is difficult. This can have important impacts in
practice, since distinguishing these cases is required to determine what classes
of inequality are needed to separate a given point. Case~\ref{case:goal-cut-2}
is unique to MIBLPs and involves removing a point that is integral. However,
in this case, we at least know that we are restricted to inequalities specific
to MIBLPs.

With respect to the goals of generating valid inequalities stated above, the
set of applicable valid inequalities for MIBLPs can be classified roughly into
three categories, as detailed in the next three paragraphs, in which $U$ is
taken to be a global upper bound on the optimal solution value. Note that the
presented classification is just to provide a rough idea about different
classes of valid inequalities and the classes are not entirely distinct.

\paragraph{Integrality Cuts.} These are inequalities that enforce integrality
and are violated by an extreme point of $\P$ for which $(x, y) \not \in X \times Y$
  (cases~\ref{case:goal-cut-1-a} and~\ref{case:goal-cut-1-b}). They are valid for
  \begin{equation*}
  \conv\left(\left\{(x, y) \in \S \midd cx + d^1 y < U \right\}\right).
  \end{equation*}
This set includes all inequalities valid for $\S$, the feasible region of the
MILP relaxation
\begin{equation*}\label{eqn:LR}
\min_{(x, y) \in \S} cx+d^1y.
\end{equation*}
These are the very same inequalities that are used in solving MILPs and
include all general families, as well as any specialized inequalities valid for
set $\S$ with particular structure. Since \MIBS{} utilizes the COIN-OR Cut
Generation Library (CGL)~\mycitep{Cgl}, all available separation routines in
this package can also be employed in \MIBS{}. Most of these inequalities are
themselves derived based on the principles described earlier in
Section~\ref{sec:gen-classes} and fall into those general categories described
there. The disjunctions utilized are typically split disjunctions, often those
involving a single variable.

\paragraph{Optimality Cuts.} These are inequalities that enforce optimality
conditions of the second-level problems and are violated by
extreme points of $\P$ of the form described in cases~\ref{case:goal-cut-1-b} and
\ref{case:goal-cut-2}. They are valid for
\begin{equation*}
\conv\left(\left\{(x, y) \in \F \midd cx + d^1 y < U \right\}\right).
\end{equation*}
Whereas the integrality cuts are meant to enforce the requirement that $(x,
y) \in \S$, the optimality cuts are meant to approximately enforce the
optimality condition
\begin{equation*}
d^2 y \leq \phi(x).
\end{equation*} 
 for the second-level problem. It is these optimality conditions that are
relaxed (in addition to the integrality conditions) in order to obtain a
tractable bounding problem. Optimality cuts are so named exactly because they
approximate the above non-linear inequality with linear inequalities (the name
is to evoke the optimality cuts arising in Benders decomposition, which are
closely related). Note that because $\phi$ is non-convex and non-concave in
general, this constraint cannot be approximated exactly by linear inequalities
(all discussed optimality cuts in this section are linear). When combined with
branching in a branch-and-cut algorithm, however, we can restrict the feasible
region to areas in which this function is convex.

\paragraph{Projected Optimality Cuts.} These are inequalities that are valid for
\begin{equation*}                      
  \conv\left(\left\{(x,y) \in \Re^{n_1\times n_2}
  \midd x\in\proj_{x}(\F), cx + \Xi(x) < U \right\}\right).
\end{equation*}
where $\Xi$ is as defined in~\eqref{eqn:Xi}.
As mentioned earlier, such inequalities can be generated to remove regions
containing only non-improving solutions, after
solving~\eqref{eqn:computeBestUB}.

In the remainder of the section, we describe the optimality cuts and projected
optimality cuts implemented in \MIBS{}.

\subsection{Generalized Chv\'atal Inequality \label{sec:generalized-chvatal}}

\subsubsection{General Chv\'atal Inequality}

\paragraph{Assumptions.}

None.

\begin{theorem}\label{thm:generalizedChvatalCut}
Let $(\hat{x},\hat{y}) \in \S \setminus \F$ be a basic feasible solution
to~\eqref{eqn:LRR} and let $H \subseteq \mathbb{N}^N$ consist of the indices
of a set of $n_1+n_2$ linearly independent inequalities (including
non-negativity) in the description of $\P$ that are binding at
$(\hat{x},\hat{y})$ (describe the polyhedral cone $\R(\hat{x}, \hat{y})$),
where $N = m_1 + m_2 + n_1 + n_2$. Let $u \in \Q^{N}$ be such that $u_i =
0$ for $i \not\in H$, $u_i > 0$ for $i \in H$, and such that
$(\alpha^x, \alpha^y, \beta)$ is a split disjunction where
\begin{align*}
(\alpha^x, \alpha^y) = u \left[\begin{array}{cc} A^1 & G^1 \\
A^2 & G^2 \\ I & 0 \\ 0 & I \end{array} \right],
&& \beta = u \left[\begin{array}{c} b^1 \\ b^2 \\ 0 \\ 0 \end{array}
\right] + 1.
\end{align*}
Then
\begin{equation*}
\alpha^x x + \alpha^y y \geq \beta \quad \forall (x, y) \in \F.
\end{equation*}
Furthermore, we have
\begin{equation*}
\alpha^x \hat{x} + \alpha^y \hat{y} = \beta - 1,   
\end{equation*}
so the inequality is violated by $(\hat{x}, \hat{y})$.
\end{theorem}
\myproof{
The inequality $(\alpha^x, \alpha^y, \beta - 1)$ is valid for $\P$, since it
is derived as a positive combination of inequalities in the description of
$\P$. Furthermore, $(\alpha^x, \alpha^y, \beta)$ is a split disjunction for
which
$$
\left\{(x,y) \in \Re^{n_1 + n_2} \midd \alpha^x x + \alpha^y
y \leq \beta - 1 \right\} \cap \P = \{(\hat{x}, \hat{y})\}.
$$
If not, then the face $\{(x, y) \in \P \mid \alpha^x x
+ \alpha^y y = \beta -1\}$ contains a non-trivial face of $\P$, which leads to
a contradiction of the fact that the set of inequalities indexed by $H$ are
linearly independent. Hence,
$$
\left\{(x,y) \in \Re^{n_1 + n_2} \midd \alpha^x x + \alpha^y
y \leq \beta - 1 \right\} \cap \F = \emptyset
$$
and $(\alpha^x, \alpha^y, \beta)$ is valid for $\F$.
}

\paragraph{Discussion.}
This inequality is a form of the generalized Chv\'atal inequality introduced
in Definition~\ref{def:generalized-chvatal}. In the stated form, it is not
clear how to choose $H$ and how to generate the weight vector $u$. Choosing
$H$ is not difficult, as it can be derived from an optimal basis in the
typical case in which $(\hat{x}, \hat{y})$ is an optimal basic feasible
solution to the LP relaxation. Choosing $u$ is more difficult. The integer
no-good inequality introduced next uses $u_i = 1 \; \forall i \in H$.
A more general way to generate such inequalities systematically
is through the use of an auxiliary optimization problem similar to the
cut-generating LP used in the case of lift-and-project for
MILPs~\mycitep{balas93}, as follows.
\begin{align*}
\min_{u \in \Q^{m_1 + m_2 + n_1 + n_2}} &\alpha^x \hat{x} + \alpha^y \hat{y}
- \beta\\ 
\alpha^x &= \left[ (A^1)^\top \;\; (A^2)^\top \;\; I \;\; 0 \right] u \\
\alpha^y &= \left[ (G^1)^\top \;\; (G^2)^\top \;\; 0 \;\; I \right] u \\
\beta &= \left[(b^1)^\top \;\;\; (b^2)^\top \;\;\; 0 \;\;\; 0 \right] u \\
u_i & \geq \epsilon (> 0) \quad \forall i \in H \\
u_i & = 0 \quad \forall i \not\in H \\
\alpha^x & \in \Z^{r_1} \times \{0\}^{n_1- r_1} \\
\alpha^y & \in \Z^{r_2} \times \{0\}^{n_2- r_2} 
\end{align*}
As written, the above problem is unbounded, so some appropriate normalization
is also needed. The general version of this cut has not been implemented
in \MIBS{}, since it seems clear that it would not be as effective as the
intersection cuts, which have the ability to eliminate infeasible integer
points and can also be used to separate non-integer points. 

\subsubsection{Integer No-good Cut \label{sec:integer-no-good}}

\paragraph{Assumptions.}
\begin{itemize}[topsep=0pt]
\item $r_1 = n_1$ and $r_2 = n_2$.
\item Vectors $b^1$ and $b^2$ and all matrices $A^1, A^2, G^1$, and $G^2$ are 
integer. 
\end{itemize}
\begin{theorem}[\mycitet{}{DeNRal09}]\label{thm:bilevelFeasibilityCut}
Let $(\hat{x},\hat{y}) \in \S \setminus \F$ be a basic feasible solution
to~\eqref{eqn:LRR} and $H \subseteq \mathbb{N}^N$ consist of the indices of a
set of $n_1 + n_2$ linearly independent inequalities (including
non-negativity) in the description of $\P$ that are binding at
$(\hat{x},\hat{y})$, as in Theorem~\ref{thm:generalizedChvatalCut}. Then,
under the stated assumptions, we have
\begin{equation*}
\alpha^x x + \alpha^y y \geq \beta \quad \forall (x, y) \in \F,
\end{equation*}
where
\begin{align*}
(\alpha^x, \alpha^y) = \mathbf{1}^N \left[\begin{array}{cc} A^1 & G^1 \\
A^2 & G^2 \\ I & 0 \\ 0 & I \end{array} \right],
&& \beta = \mathbf{1}^N \left[\begin{array}{c} b^1 \\ b^2 \\ 0 \\ 0 \end{array}
\right] + 1,
\end{align*}
where $\mathbf{1}^N$ is a row vector of dimension $N = m_1 + m_2 + n_1 + n_2$
withe all entries equal to one. Furthermore, we have
\begin{equation*}
\alpha^x \hat{x} + \alpha^y \hat{y} = \beta -1,   
\end{equation*}
so the inequality is violated by $(\hat{x}, \hat{y})$.
\end{theorem}
\myproof{
This is a special case of Theorem 5 in which $u = \mathbf{1}^N$
}

\paragraph{Discussion.}
Since we require $\alpha^x$, $\alpha^y$ and $\beta$ to integer-valued, this
inequality is a disjunctive inequality for $\F$ with respect to the split 
disjunction  
\begin{align*} 
& X_1 = \left\{(x, y) \in \Re^{n_1 + n_2} \midd \alpha^x x + \alpha^y y
  \leq \beta - 1 \right\} 
& X_2 = \left\{(x, y) \in \Re^{n_1 + n_2} \midd \alpha^x x + \alpha^y y
  \geq \beta  \right\}. 
\end{align*}
Since $\P \cap X_1 = (\hat{x}, \hat{y})$ and $(\hat{x}, \hat{y}) \not\in \F$,
it is also a generalized Chv\'atal inequality for $\F$ with respect to $\P$
and the above split disjunction. The utilized split disjunction is defined by
considering a combination of constraints of $\P$ binding at
$(\hat{x}, \hat{y})$, as described earlier.
The combination is guaranteed to yield a split disjunction because of our
assumption that the constraint matrix is integral. Essentially, this is a
simple case in which a vector $u$ satisfying the conditions of Theorem 5 can
be easily derived. 

This cut can be applied to remove an optimal solution
$(\hat{x},\hat{y}) \in \S \setminus \F$ to~\eqref{eqn:LRR} under the above
assumptions, but it does not eliminate any other bilevel infeasible members of
$\S$. It is easy to generate such inequalities because all required
information for its generation can be obtained from the optimal tableau
obtained when solving~\eqref{eqn:LRR}.

\begin{figure}[tbh]
\centering
\includegraphics[height=2.5in]{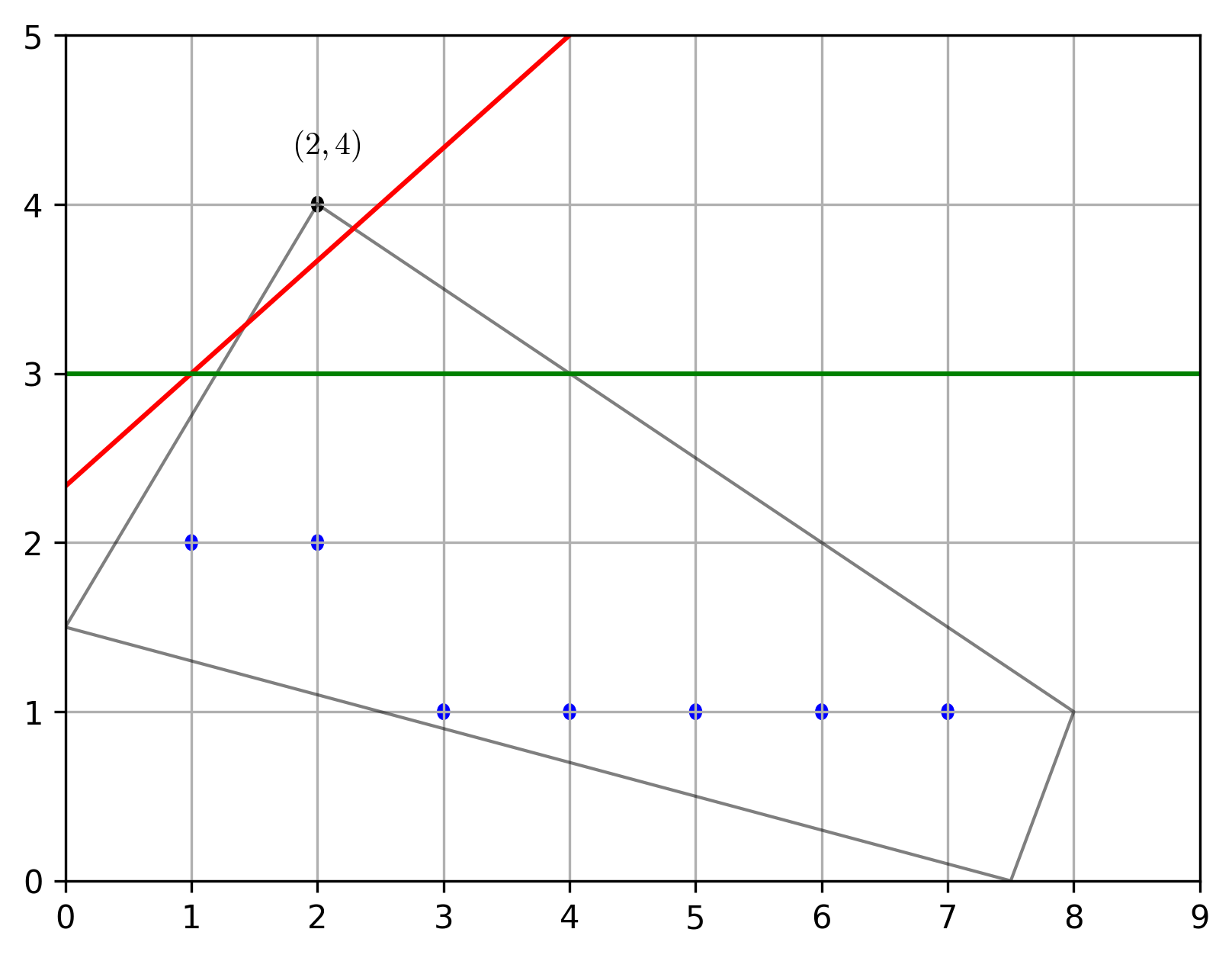}
\caption{The generated integer no-good cut for the example shown in Figure~\ref{fig:mooreexample}}
\label{fig:int-no-good-cut}
\end{figure}

\paragraph{Example.} The red line in Figure~\ref{fig:int-no-good-cut} shows the
generated integer no-good cut for removing the bilevel infeasible extreme
point $(2,4)$ of $\P$ in the example shown in Figure~\ref{fig:mooreexample}.
The inequality is obtained in two steps, as described above. We first sum the
first two inequalities in the formulation (which are the ones binding at $(2,
4)$) to obtain the inequality $-2x + 3y \leq 8$, valid for $\P$. We then
subtract one from the right-hand side to obtain the final inequality $-2x +
3y \leq 7$. As one can observe, this inequality does not remove any integer
points from $\P$ except $(2,4)$. It is not difficult to see that a deeper cut
can be obtained by using general weights to combine the inequalities.
Combining the inequalities with weight vector $u = [14, 70]^\top$ yields the
inequality $y \leq 4$, valid for $\P$. Rounding, we obtain the valid
inequality $y \leq 3$, which clearly dominates the original integer no good.  

\subsection{Benders Cuts}

\MIBS{} includes two classes of Benders cuts, one that applies to problems
with binary linking variables and a stronger class that applies only to
interdiction problems. 

\subsubsection{Benders Binary Cut} 

\paragraph{Assumptions.}

\begin{itemize}[topsep=0pt]
\item $(x, y) \in \F \Rightarrow x_L \in \B^L$ 
\end{itemize}
\begin{theorem}\label{thm:generalBendersCut}
Let $(\hat{x}, \hat{y}) \in \P$ such that $\hat{x}_\J \in \B^\J$ and
$d^2\hat{y} > d^2 y^*$ for some $y^* \in \P_1(\hat{x}) \cap \P_2(\hat{x}) \cap
Y$. 
Further let
\begin{itemize}
\item $\J^- = \left\{i \in \J \midd A^2_i \leq 0\right\}$, and
\item $\J^+ = \left\{i \in \J \midd A^2_i \geq 0\right\}$.
\end{itemize}
Then, under the assumption, we have
\begin{equation} \label{eqn:generalBendersCut}
d^2 y - M \left(\sum_{i \in \J \setminus \J^- : \;\hat{x}_i = 1} (1 - x_i) +
\sum_{i \in \J \setminus \J^+: \;\hat{x}_i = 0} x_i\right) \leq d^2y^*
\quad \forall (x,y) \in \F,
\end{equation}
where $M \geq \max\{d^2y \mid (x,y) \in \F\} - d^2y^*$. Furthermore, this
inequality is violated by all $(x, y) \in \P$ for which
$x_{\J}=\hat{x}_{\J}$ and $d^2 y > d^2y^*$. More specifically, it is violated
by all $(x, y) \in \S \setminus \F$ such that $x_{\J}=\hat{x}_{\J}$.
\end{theorem}
\myproof{
Note that $\J^+$ are the indices of linking variables for which an increase in
value enlarges the second-level feasible region. Similarly, $\J^-$ are the
indices of linking variables for which a decrease in value
enlarges the second-level feasible region.

Let $(\tilde{x},\tilde{y})\in\F$, $(\hat{x}, \hat{y}) \in \P$, and $y^* \in
\P_1(\hat{x}) \cap \P_2(\hat{x}) \cap Y$ be given as in the statement
of the theorem. We show that the inequality derived from $(\hat{x}, y^*)$ is
satisfied by $(\tilde{x},\tilde{y})$. We consider two cases.
\begin{itemize}
\item[(i)] Assume $\tilde{x}_i = 0$ for $i \in \J \setminus \J^+$ such that
$\hat{x}_i = 0$ and $\tilde{x}_i = 1$ for $i \in \J \setminus \J^-$ such that
$\hat{x}_i = 1$, so that
\begin{equation} \label{eqn:generalBenders-proof1}
\sum_{i \in \J \setminus \J^-: \;\hat{x}_i = 1} (1 - \tilde{x}_i) +
\sum_{i \in \J \setminus \J^+: \;\hat{x}_i = 0} \tilde{x}_i = 0
\end{equation}
Then we first claim that
\begin{equation*}
G^2 y^* \geq b^2 - A \hat{x} \geq b^2 - A \tilde{x},
\end{equation*}
where the first inequality in the chain follows from
$y^* \in \P_2(\hat{x})$ and the second inequality follows from our
assumption and the definition of $\J^-$ and $\J^+$. Then we have
$y^* \in \P_2(\tilde{x}) \cap Y$. It follows that we must have
\begin{align*}
d^2 \tilde{y} - M \left(\sum_{i \in \J \setminus \J^-: \;\hat{x}_i = 1} (1
- \tilde{x}_i) + \sum_{i \in \J \setminus \J^+: \;\hat{x}_i =
0} \tilde{x}_i \right) = d^2 \tilde{y} \leq d^2 y^*.
\end{align*}
Hence, $(\tilde{x},\tilde{y})$ satisfies the
inequality~\eqref{eqn:generalBendersCut} in this case.

\item[(ii)] Assume that either for some $i \in \J \setminus \J^+$, we have  
$\tilde{x}_i = 1$, $\hat{x}_i = 0$, or for some $i \in \J \setminus \J^-$, we
have $\tilde{x}_i = 0$, $\hat{x}_i = 1$. Then
\begin{align*}
d^2 \tilde{y} \leq d^2 y^* + M \leq d^2 y^* + M
\left(\sum_{i \in \J \setminus \J^- : \;\hat{x}_i = 1} (1 - \tilde{x}_i) +
\sum_{i \in \J \setminus \J^+: \;\hat{x}_i = 0} \tilde{x}_i \right).
\end{align*}
Hence, $(\tilde{x},\tilde{y})$ satisfies the
inequality~\eqref{eqn:generalBendersCut} in this case.
\end{itemize}
Since $(\tilde{x},\tilde{y})$ was arbitrary, we have that the
inequality~\eqref{eqn:generalBendersCut} is valid for $\F$.

Finally, let $(\bar{x},\bar{y})\in \P$ be given with $\bar{x}_\J=\hat{x}_\J$ and assume that
\begin{equation}
d^2\bar{y} > d^2y^*
\end{equation}
so that $(\bar{x},\bar{y}) \not \in \F$. Then by the same logic
as~\eqref{eqn:generalBenders-proof1}, we have that
\begin{align*}
d^2 \bar{y} - M \left(\sum_{i \in \J \setminus \J^-: \;\bar{x}_i = 1} (1
- \bar{x}_i) + \sum_{i \in \J \setminus \J^+: \;\bar{x}_i =
0} \bar{x}_i \right) = d^2 \bar{y} > d^2 y^*
\end{align*}
and the inequality~\eqref{eqn:generalBendersCut} is violated by $(\bar{x}, \bar{y})$. Since
$(\bar{x}, \bar{y})$ was chosen arbitrarily, we have that the inequality 
is violated by all $(x,y)\in\P$ with $x_\J=\hat{x}_\J$ and $d^2y > d^2y^*$.
}
\paragraph{Discussion.} This inequality is easily seen to be a disjunctive
inequality for $\F$ with respect to the valid disjunction
\begin{align*}
&& X_1 = \Biggl\{(x, y) \in \Re^{n_1 + n_2} \Biggl | &
\sum_{i \in \J \setminus \J^- : \;\hat{x}_i = 1} (1 - x_i) +
\sum_{i \in \J \setminus \J^+: \;\hat{x}_i = 0} x_i = 0,\; 
d^2y\leq d^2\hat{y} \Biggl\}\\
&& X_2 = \Biggl\{ (x, y) \in \Re^{n_1 + n_2} \Biggl | & 
\sum_{i \in \J \setminus \J^- : \;\hat{x}_i = 1} (1 - x_i) +
\sum_{i \in \J \setminus \J^+: \;\hat{x}_i = 0} x_i \geq 1
\Biggl\}.
\end{align*} 
Note that this inequality can be used to eliminate infeasible solutions
$(\hat{x}, \hat{y}) \in \P \setminus \F$ to the relaxation~\eqref{eqn:LRR}
with $\hat{x}_\J \in \B^\J$ even when $\hat{y} \not\in Y$, provided that we
can produce a point $y^* \in \P_1(\hat{x}) \cap \P_2(\hat{x}) \cap Y$ such
that $d^2 \hat{y} > d^2 y^*$. Such $y^*$ can be obtained, e.g., by
solving the second-level problem associated with $\hat{x}$ (the existence of
$y^* \in \P_1(\hat{x}) \cap \P_2(\hat{x}) \cap Y$ such that $d^2 \hat{y} > d^2
y^*$ is only guaranteed, however, when $\hat{y} \in Y$).

It is also important to note that this inequality is closely related to the
generalized no-good cut introduced later, which may be a preferable
alternative, due to the possible difficulty of computing an appropriate value
of $M$ (we solve an integer program to find the big-M and it is likely that
computing a valid big-M is formally NP-hard). The (possible) advantage of this
inequality is that in cases where some columns of the matrix $A^2$ have
uniform sign for first-level variables, this inequality may be stronger.

\subsubsection{Benders Interdiction Cut} 

\paragraph{Assumptions.} 
We assume the given problem is of the form~\eqref{eqn:MIPINT}.

In what follows, $G_i$ is the $i^\textrm{th}$ column of $G$ and $\P^{INT}
= \{(x, y) \in \Re^{n_1 + n_2} \mid x \in \P_1^{INT},
y \in \P_2^{INT}(x)\}$ the feasible region of the relaxation
of~\eqref{eqn:MIPINT} equivalent to~\eqref{eqn:LRR}, $\RR(x)$ is the rational
reaction set associated with $x \in \P_1^{INT}$ and $\F$ is the set of bilevel
feasible solutions, as usual.

\begin{theorem}\label{thm:bendersCut}
Let $(\hat{x}, \hat{y}) \in \P^{INT}$ such that $\hat{x}_\J \in \B^\J$ and
$d\hat{y} > dy^*$ for some $y^* \in \P_2^{INT}(\hat{x}) \cap Y$. Further, let 
\begin{itemize}
\item $\J^- = \left\{i \in \J \midd G_i \leq 0 \right\}$, and
\item $\J^+ = \left\{i \in \J \midd G_i \geq 0 \right\}$.
\end{itemize}
Then, under the assumptions, we have
\begin{equation} \label{eqn:bendersCut}
dy + \sum_{i \in \J^-} (d_i y^*_i) x_i
- \sum_{i \in \J^+: \;y^*_i = 0} d_i y_i
- M \left(\sum_{i \in \J \setminus \J^-: \;y^*_i > 0} x_i\right)
\leq d y^* \quad
\forall (x,y) \in \F,
\end{equation}
where $M \geq \max\{dy \mid (x,y) \in \F\} - d\hat{y}$.
Furthermore, this inequality is violated by all $(x, y) \in \P^{INT}$ such that
$x_{\J}=\hat{x}_{\J}$ and $d y > dy^*$ and more specifically, it is
violated by all $(x, y) \in \S \setminus \F$ such that $x_{\J}=\hat{x}_{\J}$.
\end{theorem}
\myproof{
The sets $\J^+$ and $\J^-$ play a role similar here to the role they play in
the Benders binary cut. However, in interdiction problems, it is the bounds on
the second-level variables that are parametric in the first-level variables
rather than the right-hand side of the second-level constraints. This is why
these sets are defined directly in terms of the matrix $G$ rather than the
matrix $A$. As such, $\J^+$ (resp., $\J^-$) represent indices of second-level
variables that can be increased (resp., decreased) in value without affecting
feasibility of $y^*$.

Let $(\tilde{x},\tilde{y})\in\F$, $(\hat{x}, \hat{y}) \in \P^{INT}$, and
$y^* \in \P_2^{INT}(\hat{x}) \cap Y$ be given as in the statement
of the theorem. We show that the inequality derived from $y^*$ is satisfied by
$(\tilde{x},\tilde{y})$. We consider two cases.
\begin{itemize}
\item[(i)] Assume
$\tilde{x}_i = 0$ for $i \in \J \setminus \J^-$ such that $y^*_i > 0$. This
means that only variables whose value can be decreased without affecting
feasibility are fixed to value zero by the choice of $\tilde{x}$. Then we
construct $y'\in Y$ 
\begin{equation*}
y'_i = 
\begin{cases}
0 & \textrm{if } i \in \J^- \textrm{ and } \tilde{x}_i = 1,\\
\tilde{y}_i & \textrm{if } i\in \J^+ \text{ and } y^*_i = 0, \\
y^*_i & \text{otherwise }.
\end{cases}
\end{equation*}
Then we first claim that $G y' \geq g$, since we have
\begin{itemize}
\item $G y^* \geq g$,
\item $y'_i \geq y^*_i$ for all $i \in \J^+$, 
\item $y'_i \leq y^*_i$ for all $i \in \J^-$, and
\item $y'_i = y^*_i$ otherwise.
\end{itemize}
We also have $\tilde{x}_i = 1 \Rightarrow y'_i = 0$.
Hence,  
$y'\in\P_2^{INT}(\tilde{x})\cap Y$ and since $\tilde{y}\in\RR(\tilde{x})$, we
have 
\begin{eqnarray*}
d\tilde{y} \leq d y' & = & d y^* - \sum_{i \in \J^-} (d_i y^*_i) \tilde{x}_i  +
\sum_{i \in \J^+: \;y^*_i = 0} d_i \tilde{y}_i \\
& = & d y^* - \sum_{i \in \J^-} (d_i y^*_i) \tilde{x}_i +
\sum_{i \in \J^+: \;y^*_i = 0} d_i \tilde{y}_i +
M \left(\sum_{i \in \J \setminus \J^-: \;y^*_i > 0} \tilde{x}_i\right),
\end{eqnarray*}
since $\sum_{i \in \J \setminus \J^-: \;y^*_i > 0} \tilde{x}_i = 0$. It follows that
$(\tilde{x},\tilde{y})$ satisfies the inequality~\eqref{eqn:bendersCut}.

\item[(ii)] Now assume $\tilde{x}_i = 1$ for some $i \in \J \setminus \J^-$
such that $y^*_i > 0$. Then $\sum_{i \in \J \setminus \J^-: \;y^*_i >
0} \tilde{x}_i \geq 1$ and we have
\begin{equation*}
d\tilde{y} \leq d y' + M \leq d y^* - \sum_{i \in \J^-} (d_i
y^*_i) \tilde{x}_i + \sum_{i \in \J^+: \;y^*_i = 0} d_i \tilde{y}_i +
M \left(\sum_{i \in \J \setminus \J^-: \;y^*_i > 0} \tilde{x}_i\right)
\end{equation*}
It follows that
$(\tilde{x},\tilde{y})$ satisfies the inequality~\eqref{eqn:bendersCut}.
\end{itemize}
Since $(\tilde{x},\tilde{y})$ was arbitrary, it follows that the inequality is
valid for $\F$. 

Moreover, since $y^* \in \P_2^{INT}(\hat{x})$, we have $y^*_i = 0$
for all $\{i\in L\mid \hat{x}_i = 1\}$. 
It follows that
\begin{equation}\label{eqn:benders-proof3}
dy^* = d y^* - \sum_{i \in \J^-} (d_i y^*_i) \hat{x}_i  +
\sum_{i \in \J^+: \;y^*_i = 0} d_i y^*_i +
M \left(\sum_{i \in \J \setminus \J^-: \;y^*_i > 0} \hat{x}_i\right)
\end{equation}
Finally, let $(\bar{x},\bar{y})\in\P^{INT}$ with $\bar{x}_\J=\hat{x}_\J$ and
\begin{equation}\label{eqn:benders-proof4}
d \bar{y} > dy^*
\end{equation}
be given. Then by~\eqref{eqn:benders-proof3} and~\eqref{eqn:benders-proof4},
the inequality~\eqref{eqn:bendersCut} is violated by $(\bar{x}, \bar{y})$. Since
$(\bar{x}, \bar{y})$ was arbitrary, the inequality is violated by all
$(x,y)\in\P^{INT}$ with $x_\J=\hat{x}_\J$ such that $d y > dy^*$.
}
\paragraph{Discussion.} \mycitet{Caprara et al.}{capraraetal16} proposed a
special case of this inequality for knapsack interdiction problems, though it
was not found to be useful in their algorithm. A generalized version, valid
for problems of the form~\eqref{eqn:MIPINT}, was released
in \MIBS{} \t{1.0}~\mycitep{MIBS1}, as discussed by~\mycitet{Ralphs}{Ral15}.
Similar general versions later appeared
in~\mycitet{Tahernejad}{Tahernejad2019} and~\mycitet{Fischetti et
al.}{FisLjuMonSinInterdiction19}. What is presented here is a further
generalization that does not make any assumptions on, e.g., monotonicity. This
inequality is a disjunctive inequality for $\F$ with respect to the valid
disjunction
\begin{align*}
&& X_1 = \Biggl\{(x, y) \in \Re^{n_1 + n_2} \Biggl | &
\sum_{i \in \J: \hat{x}_i = 0}x_i+\sum_{i \in \J: \hat{x}_i = 1}(1-x_i) = 0, 
dy\leq d\hat{y} \Biggl\}\\
&& X_2 = \Biggl\{ (x, y) \in \Re^{n_1 + n_2} \Biggl | & 
\sum_{i \in \J: \hat{x}_i = 0} x_i +\sum_{i \in \J: \hat{x}_i = 1} (1-x_i) \geq 1,
dy + \sum_{i \in \J^-} (d_i \hat{y}_i) x_i - \sum_{i \in \J^+: \;\hat{y}_i
= 0} d_i y_i \leq d\hat{y} \Biggl\}.
\end{align*} 

As with the Benders binary cut, the Benders interdiction cut can be exploited
to eliminate infeasible solutions
$(\hat{x}, \hat{y}) \in \P^{INT} \setminus \F$ to the relaxation~\eqref{eqn:LRR}
even when $\hat{y} \not\in Y$, provided that we can produce a point
$y^* \in \P_2^{INT}(\hat{x}) \cap Y$ such that $d \hat{y} > dy^*$ (see earlier
discussion in the previous section on the Benders binary cut).

\subsection{Intersection Cuts \label{sec:intersection-cuts}}

\MIBS{} includes several known classes of ICs that are closely
related, differing only in the form of the BFS used to generate
each of them. The first two classes are covered in this section and use
closely related BFSs. The first of these is derived from the
existence of an improving feasible \emph{solution} to the second-level problem, while
the second is derived based on the existence of improving
feasible \emph{direction} to the second-level problem. The third class is
based on an entirely different type of BFS and is discussed in
Section~\ref{sec:hyperIC}.

\subsubsection{Improving Solution Intersection Cuts \label{sec:isics}}

\paragraph{Assumptions.}
\begin{itemize}[topsep=0pt]
\item $A^2x + G^2y - b^2\in\Z^{m_2}$ for all $(x,y)\in\S$.
\item $d^2 \in \Z^{n_2}$ (Type \RNum{2} only).
\end{itemize}
\begin{theorem}[\mycitet{}{FisLjuMonSinUse18},
    \mycitep{FisLjuMonSinNew17}]\label{thm:sep1IC} Let
$(\hat{x},\hat{y}) \in \P \setminus \F$ be an optimal solution
of~\eqref{eqn:LRR} such that $d^2\hat{y} > d^2 y^*$ for some
$y^* \in \P_2(\hat{x}) \cap Y$. Then, under the stated assumptions, we have
\begin{equation}\label{eqn:sep1IC}
\alpha^x x + \alpha^y y \geq \beta \quad \forall (x, y) \in \F,
\end{equation}
where the inequality~\eqref{eqn:sep1IC} is the IC
associated with the sets $\R(\hat{x},\hat{y})$ and
\begin{equation}\label{eqn:sep1ICsetS}
\mathcal{C} = \left\{(x,y)\in\Re^{n_1\times n_2}\midd 
d^2y \geq d^2 y^*, A^2 x \geq b^2 - G^2y^* - 1 \right\}, 
\end{equation}
Furthermore, the inequality~\eqref{eqn:sep1IC}
is violated by $(\hat{x},\hat{y})$. 
\end{theorem}

\paragraph{Discussion.} Two distinct procedures can be used can to derive
violated valid inequalities from the above formula, by computing $y^*$ in two
different ways.
\emph{\textbf{Type \RNum{1}}} ISICs are obtained simply by solving the
second-level problem to determine $y^* \in \RR(\hat{x})$ (though any
$y^* \in \P_2(\hat{x}) \cap Y$ such that $d^2 y^* < d^2 \hat{y}$ will suffice).
In \emph{\textbf{Type \RNum{2}}} ISICs, $y^*$ is taken to be an optimal solution of the
following MILP, which is constructed with an objective function designed to
make the resulting BFS as large as possible. 
\begin{align} 
y^*\in&\text{argmin}\sum_{i=1}^{m_2}w_i\nonumber\\
&d^2y\leq \lceil d^2 \hat{y} - 1 \rceil \nonumber \\ 
&G^2y+(A^2\hat{x} - L)w\geq b^2 - L\label{eqn:define-y-3}\\
&y\in Y\nonumber\\
&w\in\{0,1\}^{m_2},\nonumber
\end{align}
where $L_i = \sum_{j=1}^{n_1}\min\{A^2_{ij}l_{x_j},A^2_{ij}u_{x_j}\}$ for $i =
1, ...,m_2$ and $A^2_{ij}$ represents the element of $i^{\text{th}}$ row and
$j^{\text{th}}$ column of $A^2$. Due to constraint~\eqref{eqn:define-y-3},
when $w^*_i = 0$ (where $w^*$ represents the optimal value of $w$), we can
conclude that $g^2_i y^* \geq b^2_i - a^2_i x$ for all bilevel feasible
solutions. Hence, this inequality can be removed from the definition of set
$\C$ in~\eqref{eqn:sep1ICsetS} in this case. By utilizing such property, the
defined set $\C$ for an ISIC of type \RNum{2} is hopefully larger than the one
for type \RNum{1} and may provide a stronger cut. However, it should also be
noted that finding $y^*$ for ISICs of type \RNum{2} requires solving an MILP,
while for ISICs of \RNum{1}, no additional effort is needed if the feasibility
check is already being done.

It is important to point out that these two types of ISICs make no requirement
on integrality of $(\hat{x}, \hat{y})$. We require that
$y^* \in \P_2(\hat{x}) \cap Y$, but $y^*$ need not be a member of
$\RR(\hat{x})$ (i.e., optimal to the second-level problem arising from fixing
$\hat{x}$). Type \RNum{1} ISICs are cheaply and conveniently generated when
we've already generated a $y^*$ satisfying the desired conditions either as a
by-product of a primal heuristic or a feasibility check. When such $y^*$ is
not already available, it must be generated as part of the cut generation.
Whether and when it is worthwhile to do this is an important empirical
question. Type \RNum{2} ISICs always require the solution of an auxiliary
subproblem, making them potentially more expensive.

With both types of inequalities, the separation procedure can fail in two
different ways: either $\P$ may be entirely contained in the BFS
or the subproblem for generating $y^*$
may be infeasible. The former case is detected when the rays of
$\R(\hat{x}, \hat{y})$ do not intersect the boundary of the BFS. This situation
is obviously unlikely if $(\hat{x}, \hat{y})$ is a solution the initial LP
relaxation, but when separation is done for a subproblem deeper in the
branch-and-bound tree, it can (and does) happen. The latter case arises if
either $d^2 \hat{y} \leq \phi(b^2 - A^2 \hat{x})$ or $\RR(\hat{x})
= \emptyset$. When $(\hat{x}, \hat{y}) \in \S \setminus \F$, we must have
$d^2 \hat{y} > \phi(b^2 - A^2 \hat{x})$ and $\RR(\hat{x}) \not= \emptyset$.
When $(\hat{x}, \hat{y}) \not\in \S$, we may be in
case~\ref{case:goal-cut-1-a} (the optimality conditions for the second-level
problem that were relaxed in~\eqref{eqn:LRR} are already satisfied by
$(\hat{x}, \hat{y})$) and we cannot use MIBLP inequalities to separate
$(\hat{x}, \hat{y})$. As a simple example, consider what happens in the
example from Figure~\ref{fig:mooreexample} if the first-level objective is
changed to minimizing the value of $3x - y$. Then the optimal solution to the
bilevel problem is $(1,2)$ and the solution to the relaxation~\eqref{eqn:LRR}
is $(0, \frac{3}{2})$. Since $\phi(0) = +\infty$, it is obvious that no MIBLP
cut can be generated (obviously, we have the option of separating
$(\hat{x}, \hat{y})$ with inequalities valid for $\S$ (MILP cuts)). If
$\hat{y} \in Y$ but $\hat{x} \not\in X$, then separation fails only if
$\hat{y} \in \RR(\hat{x})$ (there is no improving solution). More analysis of
the prevalence of failures in generating ICs is given later in
Section~\ref{sec:cg-failures}.

\begin{figure}[tbh]
\centering
\includegraphics[height=2.5in]{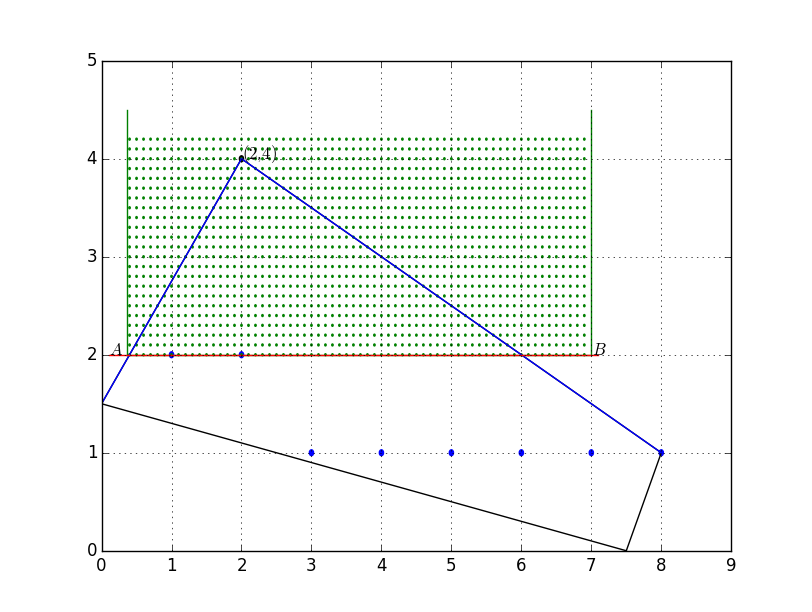}
\caption{The generated ISIC for the example shown in
Figure~\ref{fig:mooreexample}} 
\label{fig:IC}
\end{figure}

\paragraph{Example.} 
Figure~\ref{fig:IC} shows the generated ISIC (the red line) for
removing the extreme point (2,4) of $\P$ in the example shown in
Figure~\ref{fig:mooreexample}. The blue cone and the green dotted region show
$\R(2,4)$ and set $\C$, respectively. Note that in this case, the generated
inequalities and sets $\C$ are the same for both described types of ISIC.
\vskip -0.1 in 

\subsubsection{Improving Direction Intersection Cuts \label{sec:idics}}

\paragraph{Assumptions.}
\begin{itemize}[topsep=0pt]
\item $A^2x + G^2y - b^2\in\Z^{m_2}$ for all $(x,y)\in\S$.
\item $d^2 \in \Z^{n_2}$.
\end{itemize}

\begin{theorem}[\mycitet{}{FisLjuMonSinNew17}]\label{thm:watermelonIC}
Let $(\hat{x},\hat{y}) \in \P \setminus \F$ be an optimal
solution to~\eqref{eqn:LRR} and $\Delta\hat{y}\in \Z^{r_2}\times\Re^{n_2-r_2}$ 
satisfy these conditions:
\begin{itemize}
\item $d^2\Delta\hat{y} < 0$.
\item $\hat{y}+\Delta\hat{y} \in \P_2(\hat{x})$.
\end{itemize} 
Then, under the stated assumptions, we have
\begin{equation}\label{eqn:watermelonIC}
\alpha^x x + \alpha^y y \geq \beta \quad \forall (x, y) \in \F,
\end{equation}
where the inequality~\eqref{eqn:watermelonIC} is the IC
generated from the cone $\R(\hat{x},\hat{y})$ and
\begin{equation}\label{eqn:watermelonICsetS}
\C = \left\{(x,y)\in\Re^{n_1\times n_2}\midd A^2 x+ G^2(y+\Delta\hat{y})\geq 
b^2 - 1, y+\Delta\hat{y}\geq -1\right\}.
\end{equation}
Furthermore, the inequality~\eqref{eqn:watermelonIC} is violated by
$(\hat{x},\hat{y})$. 
\end{theorem}

\paragraph{Discussion.}
Generating an IDIC requires finding an improving direction
$\Delta\hat{y}\in \Z^{r_2}\times\Re^{n_2-r_2}$ for the second-level problem,
as in Definition~\ref{def:certificate-infeas} and
Theorem~\ref{thm:watermelonIC}. Given an improving direction, an IC can be
derived from the associated BFS defined in~\eqref{eqn:IS-BFS}. As with ISICs
of type \RNum{2}, the aim is to find the largest possible convex set and this
is attempted by solving the following subproblem.
\begin{equation} \label{eqn:define-deltay-watermelonIC}
\begin{aligned}
\Delta\hat{y}\in&\text{argmax}\sum_{i=1}^{m_2}w_i + \sum_{i=1}^{n_2}v_i\\
&d^2\Delta y\leq -1\\
&G^2\Delta y\geq b^2-A^2\hat{x} - G^2\hat{y}\\
&\Delta y\geq -\hat{y}\\
&G^2\Delta y\geq w\\
&\Delta y\geq v\\
&\Delta y\in \Z^{r_2}\times\Re^{n_2-r_2}\\
&w\leq 0\\
&v\leq 0.
\end{aligned}
\end{equation}
The optimal values of the vectors $w$ and $v$ can be employed to drop some of
the inequalities in the definition of set~\eqref{eqn:watermelonICsetS},
enlarging our convex set further. $w^*_i = 0$ ($v^*_i=0$, resp.) means that
$g_i^2(y+\Delta\hat{y})\geq b_i^2-a_i^2x$ ($y_i+\Delta\hat{y}_i\geq 0$, resp.)
for all bilevel feasible solutions, so this inequality can be removed from the
definition of set~\eqref{eqn:watermelonICsetS}.

As with ISICs of types \RNum{1} and \RNum{2}, the IDIC makes no assumption
about the point to be separated and it can thus be used to separate points
$(\hat{x}, \hat{y})$ for which $\hat{x} \not\in X$ and/or $\hat{y} \not \in
Y$. Also as with ISICs, separation can fail in the same two ways: $\P$ may be
entirely contained in the BFS or~\eqref{eqn:define-deltay-watermelonIC} may be
infeasible. In the latter case, the infeasibility can be a result of any one
of \emph{three} conditions:
\begin{enumerate}
\item $\RR(\hat{x}) = \emptyset$,
\item $d^2 \hat{y} \leq \phi(b^2 - A^2 \hat{x})$ (we are in
case~\ref{case:goal-cut-1-a}), or 
\item $d^2 \hat{y} > \phi(b^2 - A^2 \hat{x})$ but there is no improving
direction satisfying the integrality conditions
 of~\eqref{eqn:define-deltay-watermelonIC}.
\end{enumerate}
When $\hat{y} \in
Y \setminus \RR(\hat{x})$, then~\eqref{eqn:define-deltay-watermelonIC} must be
feasible.

To illustrate these concepts, it's useful to consider
again what happens in the example of Figure~\ref{fig:mooreexample} if the
first-level objective is changed to minimizing the value of $3x - y$. As
before, the optimal solution to the MIBLP is $(1,2)$ and the solution to the
relaxation~\eqref{eqn:LRR} is $(0, \frac{3}{2})$. In this case, it is obvious
that no cut can be generated simply because it is clear that no improving
direction exists, which means~\eqref{eqn:define-deltay-watermelonIC} is
infeasible. The generation of improving directions when
$(\hat{x}, \hat{y}) \not\in \S$ will play an important role in the
computational experiments later.

\begin{figure}[tbh]
\centering
\includegraphics[height=2.5in]{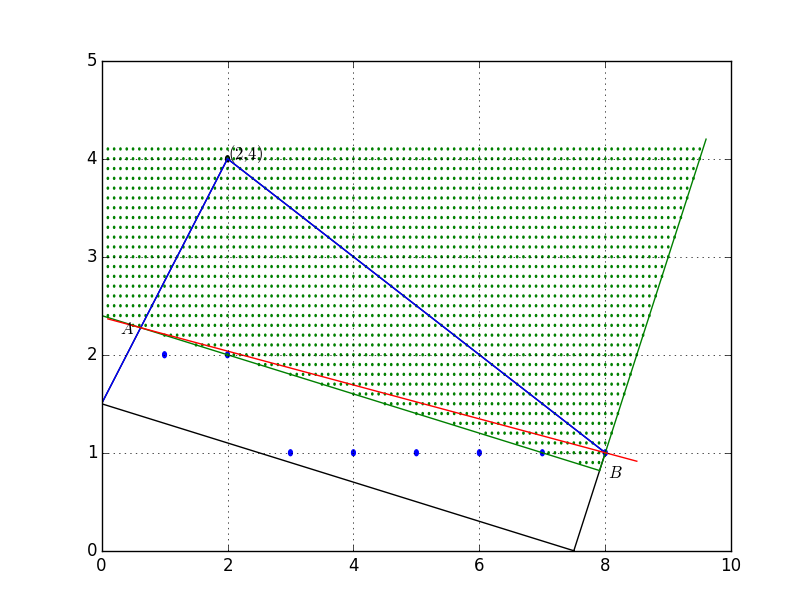}
\caption{The generated IDIC for the example shown in
Figure~\ref{fig:mooreexample}} 
\label{fig:watermelonIC}
\end{figure}

The close connection between ISICs and IDICs is evident in the fact that if
$y^* \in \P_2(\hat{x}) \cap Y$ is an improving solution, then $y^* - \hat{y}$
is an improving direction. Importantly, however, it can be the case that $y^*
- \hat{y} \not\in \Z^{r_2}\times\Re^{n_2-r_2}$ and this is the reason why some
points can be separated by ISICs but not by IDICs. In fact, neither class of
inequalities strictly dominates the other, as one can observe in
Figures~\ref{fig:IC} and~\ref{fig:watermelonIC}.

\paragraph{Example.} 
Figure~\ref{fig:watermelonIC} shows the
generated IDIC (the red line) for removing the 
extreme point (2,4) of $\P$ in the example shown in
Figure~\ref{fig:mooreexample}. The blue cone and
the green dotted region show $\R(2,4)$ and set $\C$, respectively. Note that
the points in the triangular region that is neither in the BFS nor in
$\conv(\F)$ are the points that can be separated by an ISIC in this case, but
not an IDIC. 
\vskip -0.1 in 

\subsubsection{Hypercube Intersection Cut}\label{sec:hyperIC}
\paragraph{Assumptions.}
\begin{itemize}[topsep=0pt]
\item None.
\end{itemize}
\begin{theorem}[\mycitet{}{FisLjuMonSinNew17}]\label{thm:hypercubeIC}
Let $(\hat{x},\hat{y}) \in \P$ be an optimal
solution of~\eqref{eqn:LRR} with $\hat{x}_{\J}\in\Z^{\J}$. Then, under the 
stated assumption, we have
\begin{equation}\label{eqn:hypercubeIC}
\alpha^x x + \alpha^y y \geq \beta \quad \forall (x, y) \in \F 
\textrm{ such that } x_{\J}\not = \hat{x}_{\J}, 
\end{equation}
where the inequality~\eqref{eqn:hypercubeIC} is the 
IC generated from the cone $\R(\hat{x},\hat{y})$ and convex set
\begin{equation*}\label{eqn:sep1HyperICsetS}
\mathcal{C} = \left\{(x,y)\in\Re^{n_1\times n_2}\midd \hat{x}_i-1\leq 
x_i\leq \hat{x}_i+1\quad \forall i\in \J \right\}.
\end{equation*}
Furthermore, the inequality~\eqref{eqn:hypercubeIC} is violated by
$(\hat{x},\hat{y})$. 
\end{theorem}

\paragraph{Discussion.} Since this inequality may eliminate bilevel feasible
solutions for which $\gamma = \hat{x}_{\J}$, the
problem~\eqref{eqn:computeBestUB} must first be solved with $\gamma =
\hat{x}_{\J}$ prior to generating such a cut. In this way, we ensure that all
removed bilevel feasible solutions will be non-improving. As with previous
classes, we can attempt to generate such a cut even when $\hat{y} \not\in Y$.
Even when~\eqref{eqn:computeBestUB} is infeasible, which can happen, the inequality
is still valid.

\begin{figure}[tbh]
\centering
\includegraphics[height=2.5in]{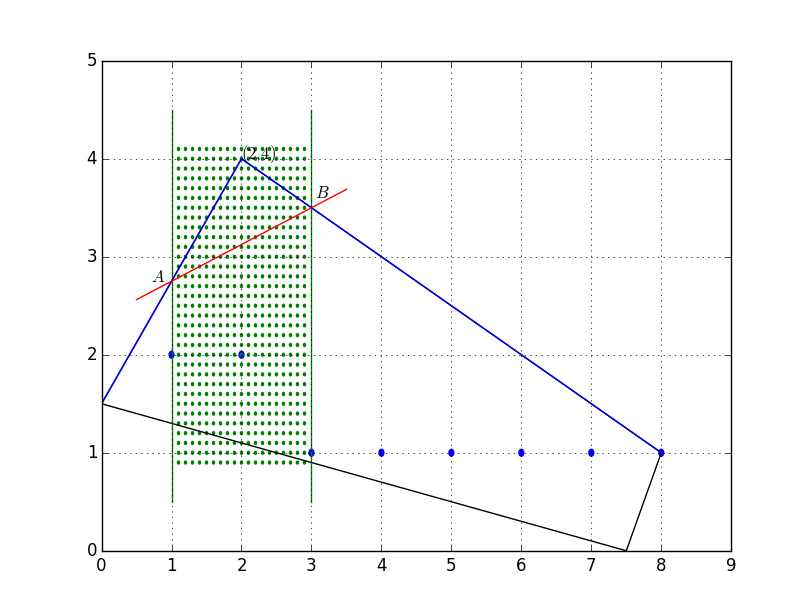}
\caption{The generated hypercube IC for the example shown in Figure~\ref{fig:mooreexample}}
\label{fig:hypercubeIC}
\end{figure}

\paragraph{Example.} Figure~\ref{fig:hypercubeIC} shows the
generated hypercube IC (the red line) for removing the 
extreme point (2,4) of $\P$ in the example shown in 
Figure~\ref{fig:mooreexample}. The blue cone and
the green dotted region show $\R(2,4)$ and set $\C$, respectively.
\vskip -0.1 in 

\subsection{Generalized No-good Cut}

\paragraph{Assumptions.}

\begin{itemize}[topsep=0pt]
\item $x_L \in \B^{\J}$.
\end{itemize}
\begin{theorem}\label{thm:binaryNoGoodCut}
Let $\gamma\in\B^{\J}$. Then, under the desired assumption, we have
\begin{equation}\label{eqn:genNoGoodCut}
\sum_{i \in \J: \gamma_i = 0} x_i + \sum_{i \in \J: \gamma_i = 1} (1-x_i) \geq
1 \quad \forall (x, y) \in \F \textrm{ such that } x_{\J}\not = \gamma.
\end{equation}
Furthermore, the inequality~\eqref{eqn:genNoGoodCut} is violated by all 
$(x,y)\in\P$ with $x_{\J}=\gamma$.
\end{theorem}
\myproof{
When $x_L = \gamma$, we have $\sum_{i \in \J: \gamma_i = 0} x_i  = 0$ and
$\sum_{i \in \J: \gamma_i = 1} (1-x_i)  = 0$, which follows that
\begin{equation}\label{eqn:noGoodCutProofEq}
x_L = \gamma \Rightarrow \sum_{i \in \J: \gamma_i = 0} x_i + 
\sum_{i \in \J: \gamma_i = 1} (1-x_i) =0.
\end{equation}
Furthermore, when $x_L \not= \gamma$, 
at least one of the following happens 
\begin{itemize}
\item $\exists i\in\J$ such that $\gamma_i =0$ and $x_i=1\Rightarrow 
\sum_{i \in \J: \gamma_i = 0} x_i  \geq 1$.

\item $\exists i\in\J$ such that $\gamma_i =1$ and $x_i=0\Rightarrow
\sum_{i \in \J: \gamma_i = 1} (1-x_i)\geq 1$.
\end{itemize}
Hence, we have
\begin{equation}\label{eqn:noGoodCutProofNotEq}
x_L \not= \gamma \Rightarrow \sum_{i \in \J: \gamma_i = 0} x_i + 
\sum_{i \in \J: \gamma_i = 1} (1-x_i) \geq 1.
\end{equation}
 The result follows from~\eqref{eqn:noGoodCutProofEq} 
 and~\eqref{eqn:noGoodCutProofNotEq}.
}

\paragraph{Discussion.} The generalized no-good cut is the only projected
 optimality cut currently in \MIBS{} and is a generalization of the no-good
cut suggested by~\mycitet{DeNegre}{DeNegre2011}. The cut proposed
in~\mycitet{DeNegre}{DeNegre2011} removes all solutions with the
same first-level values, while the new version strengthens the former one
by separating the solutions with the same linking part. Moreover, since
the no-good cut is valid only for the problems in which all first-level
variables are binary, the assumption of the generalized one is less strict.

As with the hypercube IC, the corresponding 
problem~\eqref{eqn:computeBestUB} should be solved prior to 
generating a generalized no-good cut. In contrast with the 
hypercube ICs, it is guaranteed that the generalized 
no-good cut eliminates all solutions with the target linking values. 
Unlike the generalized no-good cut, however, the hypercube IC
is applicable for general MIBLPs.

The generalized no-good cut can be seen as a disjunctive cut arising from the
disjunction
\[
x_L = \gamma  \quad\quad \textrm{ OR } \quad\quad \sum_{i \in L: \gamma_i = 0} x_i + 
\sum_{i \in L: \gamma_i = 1} (1-x_i) \geq 1.
\]
If we have already solved the problem~\eqref{eqn:computeBestUB} with the
linking variables fixed to $\gamma$, we can conclude that the disjunctive term
on the left above can contain no \emph{improving} solutions. Therefore, we can
impose the term on the right as a valid inequality. Technically speaking, this
is then actually a generalized Chv\'atal inequality with the nonzero weights
being on the bounds of the linking variables.

\subsection{Stronger Valid Inequalities \label{sec:strong-cut}}

As in the MILP case, there is a limit to the strength of inequalities that can
be obtained through ``simple'' procedures, such as the ones described in the
last section. In the case of MILPs, cut generation procedures are usually
limited to those with running times that are polynomial in the encoding size
of the given instance (generally at most solving an LP), though it is clear
from the equivalence of optimization and separation~\mycitep{GroLovSch88} that
inequalities generated by such procedures will necessarily be relatively weak
in general. Of course, recursive application can obviously lead to much
stronger inequalities (those of higher \emph{rank}).

Since MIBLPs are hard for the complexity class $\Sigma^\Pcomplexity_2$, little
can be expected from separation procedures with polynomial running times. Most
of the procedures described in this paper require instead the solution of
problems that are $\NPcomplexity-$hard (solving an MILP subproblem). Once
again, because of the equivalence of optimization and separation, even
procedures involving solving subproblems that are MILPs cannot, in general, be
expected to produce strong inequalities, except through recursive procedures
that yield inequalities of higher rank. Even inequalities
generated by solving MILP subproblems can be weak because the
relaxation~\eqref{eqn:LRR} is weak to begin with. The weakness of the
relaxation~\eqref{eqn:LRR} comes mainly from relaxing the optimality
conditions for the second-level problem that require
\begin{equation}
d^2 y \leq \phi(b^2 - A^2 x) \;\forall (x, y) \in \F.
\end{equation}
Thus, all valid inequalities that are cuts with respect to the
relaxation~\eqref{eqn:LRR} can be thought of as somehow approximating these
optimality conditions by imposing a collection of weaker (linear)
inequalities.

To obtain inequalities stronger than the ones discussed herein, it is
necessary to consider subproblems more difficult (in principle) than simple
MILPs. The most straightforward approach is to focus on Benders inequalities
of the form~\eqref{eqn:benders} in which the value function $\phi$ is replaced
by an upper approximation (a so-called \emph{primal function} (see, e.g.,
\mycitep{BolConRalTah20}). \mycitet{Tahernejad}{Tahernejad2019} investigated the
generation of one such class of inequalities in which $d^2y$ is bounded above
by a constant. This requires the solution of a (simpler) bilevel optimization
problem, which can be solved approximately to obtain a parametric class of
inequalities that can be strengthened after branching. Another class of
Benders inequalities utilized in a pure Benders algorithm is described
in~\mycitet{Bolusani and Ralphs}{BolRal22}. These inequalities require the
recursive approximation of the full value function. In general, the question
of how to generate such stronger inequalities efficiently is very much an open
one.

\subsection{Comparing Classes Analytically \label{sec:comparing}}

Having now discussed in detail the specific classes of inequalities
implemented in \MIBS{}, we summarize by grouping and comparing the various
classes discussed. A rough comparison can be done by considering the size of
the full set of bilevel infeasible solutions that are guaranteed to be
violated by a given inequality. Obviously, this is a rough way of comparing,
but we will see that this rough guide does in fact align well with the
empirical results reported in Section~\ref{sec:computation}. Roughly speaking,
\begin{itemize}
\item Generalized Chv\'atal cuts are guaranteed to remove only a single point
$(x, y) \in \S \setminus \F$ 
\item HICs and Generalized no-good cuts remove all $(\hat{x}, y) \in \S$
(feasible or not) for some $\hat{x} \in X$ (once $\Xi(\hat{x})$ is computed
  by solving~\eqref{eqn:computeBestUB}), i.e., all $(\hat{x},
  y) \in \F$ for \textbf{fixed $\hat{x}$}.
\item Benders cuts and ISICs remove all $(x, y) \in \P$ such that $y^*
 \in \P_2(x)$ and $d^2 y > d^2 y^*$, i.e.,, all $(x, y^*)$ for which a
 \textbf{fixed $y^*$} provides a certificate of infeasibility.
\item IDICs remove all $(x, y) \in \P$ such that $\Delta y$ is an improving
feasible direction for $y$, given $x$, i.e., all $(x, y) \in \P$ for
    which \textbf{a fixed $\Delta y$} provides a certificate of infeasibility.
\end{itemize}
Very roughly speaking, we may expect that the set of solutions removed by
IDICs should be larger, given that the set includes solutions without either
$x$ or $y$ fixed to a particular value. The intuition is in fact borne out in
practice, as we show in the second part of the paper. 

\section{Implementation Details \label{sec:implementation}}

\MIBS{} is built on top of the \BLIS{}~\mycitep{BLIS} MILP solver, which is in
turn built on the \CHiPPS{}~\mycitep{RalLadSal04} framework. \BLIS{} provides
the standard control mechanisms and features of a modern implementation of
branch-and-cut for solution of MILPs. Most importantly, this includes
parameters for managing the LP relaxation (primarily, controlling the addition
and removal of cuts). There are, however, some significant differences between
the mechanisms for managing the generation of valid inequalities in the MILP
and MIBLP cases. One of the key takeaways of the results presented here is
that the performance of these inequalities is intricately tied to and affected
by details of the overall solution strategy (particularly branching) of the
underlying MILP solver to a greater extent than in the MILP case. For this
reason, it appears to be much more difficult to tease out the effect on
performance of the individual components of the algorithms, especially cut
generation.
In the remainder of this section, we briefly discuss important control
parameters that all must be tuned to achieve performance. We defer until
Section~\ref{sec:computation} detailed discussion of empirical results that lead
to the current set of default parameters.

\subsection{Cut Generation Parameters \label{sec:cut-strat}}

\paragraph{Standard Control Parameters.} The underlying solver
framework \BLIS{} provides a complete set of control parameters similar to
what one would find in other MILP solvers. We mention here only the most
important ones for our purposes. Specifically, there are parameters that
\begin{itemize}
\item control whether inequalities of a particular class should be
generated at all and specify strategy for generation (various
parameters \t{cut*Strategy}); 
\item limit the number of cut passes (\t{cutPass});
\item limit the number of cuts that can be added to the relaxation in total
(\t{cutFactor}, relative to the original number of constraints);
\item control whether cuts get added based on their dynamism
(\t{scaleConFactor}; see~\mycitet{Dey and Molinaro}{DeyMol18});
\item control whether cuts are added based on their density
(\t{denseConFactor}); and
\item control \emph{tailoff detection} (\t{tailOff}).
\end{itemize}
More will be said about the specific effect of these parameters later. 

\paragraph{MIBLP-specific Control Parameters.}

In addition to the above standard parameters, there are parameters specific to
\MIBS{} that further control cut generation. The determination of
whether or not to attempt to generate inequalities of a particular class in a
given iteration at a given node is also a function of
\begin{enumerate}
\item the integrality status of the current solution (whether the linking
variables and/or second-level variables are integer-valued); and \label{i}
\item whether the second-level problem has been solved with respect to the
first-level part of the solution to the LP relaxation, which is in turn
controlled by another set of parameters. \label{ii}
\end{enumerate}
We discuss Item~\ref{ii} in more detail below in Section~\ref{sec:oracle};
Item~\ref{i} relates to the fact that some inequalities can only be used to
separate a solution to the relaxation from the feasible region when the
solution has certain structure, as described in Section~\ref{sec:classes}. The
structure required for the inequalities described in this paper is as follows.
\begin{itemize}
\item $(\hat{x}, \hat{y}) \in \S$: Integer no-good, Generalized Chv\'atal.
\item $\hat{x}_\J \in \Z^\J$: Benders binary, Benders interdiction,
Generalized no-good, Hypercube IC.
\item $(\hat{x}, \hat{y}) \in \P$: ISIC of types \RNum{1} and \RNum{2}, IDIC
(though separation is not guaranteed when $(\hat{x}, \hat{y}) \not\in \S$, as
discussed earlier).
\end{itemize}
For all classes except for ISICs and IDICs, we have no choice but to separate only
solutions with the specified structure and this means waiting to generate cuts
until a solution to the LP relaxation with such structure arises. It is
possible this will not happen with high frequency, which limits the impact of
generating inequalities of these classes.

The situation is different with ISICs and IDICs. Because these cuts are
relatively expensive to generate (requiring an MILP oracle call in general),
the strategy for generating them should take into account not only the
difficulty of the oracle problem, but the probability of success at separating
the solution. As we have previously mentioned, cut generation is only
guaranteed to succeed when the solution is in $\S$ (fully integral) and this
conservative policy of only attempting to separate fully integer solutions was
the one taken in \MIBS{} \t{1.1}. While this seemed like a logical approach at the
time and we had not questioned it, we realized during the writing of this
paper and the development of the new version of \MIBS{} that restricting
generation of cuts \emph{only} to iterations where the solution is in $\S$
severely limits the frequency of cut generation and therefore the
effectiveness.

\MIBS{} \t{1.2} has additional parameters controlling precisely
what types of points the solver should (attempt to) separate. Solutions in
$\S$ are always separated. Separation of other types of solutions are
controlled by which of the following strategies is chosen. For the present, we
offer six strategies for generating intersection cuts.
\begin{itemize}
\item \t{Always}: Try to separate every iteration.
\item \t{AlwaysRoot}: Try to separate every iteration, but only in the
root node.
\item \t{XYInt}: Try to separate when $(\hat{x}, \hat{y}) \in \S$.
\item \t{LInt}: Try to separate when $\hat{x}_L \in \Z^L$.
\item \t{YInt}: Try to separate when $\hat{y} \in Y$.
\item \t{YLInt}: Try to separate when either $\hat{x}_L \in \Z^L$ or
$\hat{y} \in Y$. 
\end{itemize}
The \t{Always} strategy results in the maximum number of cuts being
generated (and the smallest search trees) but can also result in wasted
effort when cut generation failure rates are high. On the other end of the
spectrum, the \t{XYInt} strategy ensures a failure rate near zero, but
also means many fewer cuts are generated. We analyze the empirical performance
of these strategies in Section~\ref{sec:computation}.

Finally, it should be mentioned that there are situations in which cut
generation is mandatory and parameters settings must be over-ridden. Most
prominently, when branching only on variables with fractional values, we are
required to generate a cut whenever $(\hat{x}, \hat{y}) \in \S \setminus \F$.
In such a case, all parameters that would prevent a cut from being generated
or added are ignored.

\subsection{Controlling Oracle Calls \label{sec:oracle}}

The implementation of the overall bounding loop was described in Algorithm 1
of~\mycitet{Tahernejad et al.}{TahRalDeN20}. In this scheme, parameters control
when the feasibility of a given solution is checked. Unlike in the MILP case,
checking feasibility is an expensive operation and it can be advantageous to
delay it, even if that may results in additional branching.

Similarly, the generation of most classes of inequalities described here
requires a call to an MILP oracle. In the case of Benders binary cuts, Benders
interdiction cuts, Generalized no-good cuts, ISICs of type \RNum{1}, and
Hypercube ICs, the oracle is the very same one that is used to check
feasibility, which ties the feasibility check to cut generation in an
important way.

ISICs of type \RNum{2} and IDICs have their own associated oracle calls, which
currently serve no other purpose and is called only when cuts are to be
generated. The oracle that solves the second-level problem must be called
either when we want to check feasibility of a solution \emph{or} when a cut is
to be generated.
This has the apparently strange consequence that we may sometimes solve the
second-level problem even when $\hat{x} \not\in X$. When $\hat{x} \in X$, the
effect of solving the second-level problem is to produce a point
$y^* \in \RR(\hat{x})$ so that $(\hat{x}, y^*) \in \F$ and can hence be treated
as a new heuristic solution for the purpose of improving the global upper
bound as a side benefit. When $\hat{x} \not\in X$, then $(\hat{x},
y^*) \not\in \F$, but $y^*$ can still be used as an input to generate cuts.

\MIBS{} attempts to generate cuts requiring a
second-level solution only when that solution is already available. This means
that if parameter settings for cut generation are such that points for which
$\hat{x} \not\in X$ are to be separated, we must also ensure that the
second-level problem is solved, though it would not otherwise be.

In~\mycitet{Tahernejad et al.}{TahRalDeN20}, we determined that doing the
feasibility check only when the $(\hat{x}, \hat{y}) \in \S$ or when the values
of all linking variables have been fixed (and we can thus prune the node by
solving~\eqref{eqn:computeBestUB}) was the most advantageous default settings.
However, we noted in this study that this results in certain classes of
inequalities, e.g., Hypercube ICs and ISICs of type \RNum{1}, being rarely or
never generated in many of the instances in our test set. We thus experimented
with solving the second-level problem whenever linking variables were
integer-valued and we report on those experiments in
Section~\ref{sec:computation}.

\subsection{Branching Strategy \label{sec:branching-strat}}

Although this paper is ostensibly about valid inequalities, the branching
strategy employed has an intricate effect both on the generation of valid
inequalities and the overall effectiveness of the algorithm. 
The underlying solver provides a number of the standard strategies and control
parameters for branching that are typically used in solving MILPs.
Version \t{1.2} of \MIBS{} uses a basic pseudo-cost
strategy~\mycitep{benichou.et.al:71} to select the final branching candidate by
default and implements three strategies for determining the set of candidates
for branching.
\begin{enumerate}
\item Select the branching candidate from among all variables with fractional
value in the solution to the current LP relaxation (we refer to this
as \t{fractional} branching). 
\item Prioritize branching on linking variables, only branching on non-linking
variables if the values of all linking variables are fixed by branching
constraints (we refer to this as \t{linking} branching). This means we may
need branch on linking variables whose values are not fixed but are integer in
the solution to the LP relaxation. The mechanisms for this and the reasons why
branching on integer variables may be necessary in solving MIBLPs is discussed
in~\mycitet{Tahernejad et al.}{TahRalDeN20}.
\item Prioritize branching on lower-level variables, as long as some
lower-level variable has a fractional value (the intuition behind the possible
efficacy of this strategy is explained below).
\end{enumerate}
There are several goals that are implicitly being traded off in the choice of
branching strategy and we reserve detailed discussion of how these strategies
behave empirically for Section~\ref{sec:branching-results}.

\section{Empirical Analysis \label{sec:computation}}

In this section, we present detailed empirical analyses. The goals of the
experiments are several-fold. First and foremost, we aim to provide an
objective evaluation of the effectiveness of generating inequalities from the
classes we've discussed. Separation routines for all of the described classes
of valid inequalities have been implemented in the open-source solver \MIBS{}
version \t{1.2.2}. Extensive experiments were conducted in order to gain
insight into the performance of the different classes of valid inequalities
for MIBLPs. In addition, we've tried to gain insight into how the strategy for
cut generation interacts with the other elements of the branch-and-cut
algorithm, particularly the branching strategy. Although the theory presented
parallels that of the MILP case in many respects, the practical aspects of cut
generation for MIBLPs are not as closely aligned as one might expect in a
number of ways. As we've described and further emphasize below, however,
teasing out the overall effect of generating valid inequalities is not easy.

We emphasize that this work should only be considered as the first step
towards a long-term goal of understanding how to manage the cut
generation procedure within a general branch-and-cut algorithm for MIBLPs. For
the MILP case, this maturation process involved many research works and the
accumulation of experience by many solver developers. Though this vast trove
of knowledge about the MILP case gave us a head start on the MIBLP case,
there is good reason to believe that the MIBLP case also requires further deep
investigation and the accumulation of a similar degree of empirical
experience. What we present here is a further step along this path, but likely
raises as many questions as it answers.

The remainder of this section is organized as follows. In
Section~\ref{sec:goals}, we discuss the goals of the experiments and how we
measured performance. This is a bit different and more difficult than in the
MILP case and so requires discussion. Following that, in
Section~\ref{sec:test-set}, we discuss the test set and its properties, after
which we describe the experimental setup in Section~\ref{sec:setup}. Finally,
in Section~\ref{sec:empirical}, we discuss the empirical results.

\subsection{Measures of Performance \label{sec:goals}} 

The effectiveness of generating a particular class of inequalities is more
difficult to measure in the MIBLP setting than in the MILP one, but we first
review a number of common measures of performance for classes of inequalities
that are typical in the MILP setting. Perhaps the most objective measure of
performance, because it is independent of any other details of the algorithm,
is the bound obtained by optimizing over the so-called closure associated with
the given class of valid inequalities. Informally, this closure is the region
described by the initial relaxation strengthened by all inequalities in the
class (see~\mycitet{Lodi et.al}{LodRalWoe14} for a discussion). It is thus a
measure of the maximum amount of improvement in the relaxation bound that can
be achieved in the root node. The closure bound has some drawbacks. First, it
is sometimes not easily computed and can require a separately designed
algorithm (see, e.g.,~\mycitep{BalSaxOptimizing08}). Second, it may or may not be
indicative of performance in practice, since the closure bound may only be
achieved as a result of the interaction of particular sets of inequalities in
the class and may or may not actually be achievable in practice. And of
course, as with all empirical measure, the closure bound is computed for each
instance individually, so results depend on the test set and general
conclusions may be difficult to draw.

The appeal of the closure bound is its objectivity, but in practice, more
subjective measures may be a better indicator of practical effectiveness.
Common empirical measures include (1) the improvement in the bound achieved in
the root node using the standard root processing of the branch-and-cut solver
within which the testing is taking place; (2) the reduction in the size of the
branch-and-bound tree when generating valid inequalities from a given class
over a baseline without; and (3) the reduction in empirical running time with
and without cut generation with respect to a similar baseline. The difficulty
with these measures is that it is much more difficult to separate out the
effects of ostensibly unrelated algorithmic parameter settings, particularly
the branching strategy.

In version \t{1.1} of \MIBS{}, cut generation could only be done once integrality
of either the upper-lever or of both upper- and lower-level solution had been
achieved, which typically required branching. It also resulted in relatively
few cuts being generated and made root bound improvement a largely
meaningless measure. As a result of the possibility of separating arbitrary
fractional solutions, the processing of the root node in \MIBS{} \t{1.2} plays a
role more similar to its role in MILP, where the strategy is generally to
expend substantial effort in improving both the primal and dual bounds before
switching to a phase in which there is more focus on enumeration. Essentially,
root processing can be seen as a pre-processing step in which the initial
formulation is tightened as much as possible.

Nevertheless, some classes of inequality can still only be used (or at least
are most effectively used) to separate points in $\S$ or with specific
structure, while it is typically the case that solutions generated in the root
node are not members of $\S$ (i.e., some integer variables have fractional
values).
For all these reasons, we rely on a ``holistic'' collection of
measures: root bound improvement (when appropriate), the time to optimality
(for solved instances), the size of the branch-and-bound tree (for solved
instances), and the final optimality gap (for unsolved instances). These
measures are summarized using several different kinds of profiles, each of
which reveals a different aspect of the results.
\begin{itemize}
\item \emph{Performance profiles} are empirical cumulative distribution functions
(CDFs) of ratios of a given performance measure against the virtual
best~\mycitep{dolan02}. We use performance profiles primarily for comparing
CPU time across instances that could be solved to optimality.
\item \emph{Baseline profiles} are similar to performance profiles, but present
empirical CDFs of ratios of a given performance measure against a fixed
baseline rather than against the virtual best. We use baseline profiles for
comparing both the reduction in tree size and the reduction in root gap, since
there is a natural baseline (the tree size/root gap when solving with no MIBLP cuts).
\item \emph{Cumulative profiles} are a combination of two profiles. On the left side
of the plot is the empirical CDF of the fraction of instances solved within a
given time interval. On the right side is the empirical CDF of the fraction of
instances that achieved a given final gap within the time limit (the lines on
the two halves of the graph connect because the fraction of instances solved
within the time limit is the same as the fraction of instances with zero gap
at the time limit). We use cumulative profiles to compare performance across
all instances, including those that could not be solved within the time limit. 
\end{itemize}
We filtered the instances as follows.
\begin{itemize}
\item For all profiles, we excluded instances that were solved in less than 1
second by all methods and instances that were solved in less than .01 seconds
by any one method,
\item For performance and baseline profiles, we additionally dropped instances
that could not be solved by any of the methods. 
\item For cumulative profiles, we dropped instances for which no solution was
found by any method.
\end{itemize}

\subsection{Test Set \label{sec:test-set}}

Perhaps as important a factor in performance as the details of the
implementation is the specific classes of problems under consideration. As
with MILPs, specific classes of inequalities can behave very differently on
different data sets. The seven data sets employed in our testing have
differing properties summarized in Table~\ref{tab:dataSetSummary}. (Note that
there is now a standard library of test instances called
BOBILib~\mycitep{ThuKleLjuRalSchBOBILib24} that includes these instances and
more, but this instance library came out after most of these experiments had
been completed).
\begin{table}[tbh]
\begin{center}
\begin{tabular}{c c c c c c c c c c c}
\hline\noalign{\smallskip}
Data Set &
\begin{tabular}[c]{@{}c@{}} \# of \\ instances \end{tabular} &
\begin{tabular}[c]{@{}c@{}} Variable \\ Type \end{tabular} &
\begin{tabular}[c]{@{}c@{}} \# of \\ variables \end{tabular} &
\begin{tabular}[c]{@{}c@{}} \# of \\ constraints \end{tabular} &
\begin{tabular}[c]{@{}c@{}} Degree of \\ alignment \end{tabular} &
\begin{tabular}[c]{@{}c@{}} Original \\ source \end{tabular} \\
\hline
INT-DEN & 300 &
          \begin{tabular}[c]{@{}c@{}} B \\ B \end{tabular} &
          \begin{tabular}[c]{@{}c@{}} 10-40 \\ 10-40 \end{tabular} &
          \begin{tabular}[c]{@{}c@{}} 1 \\ 11-41 \end{tabular} &
          -1 & 
          \begin{tabular}[c]{@{}c@{}} Interdiction \\ \mycitep{DeNegre2011} \end{tabular} \\
\hline
DEN     & 50 &
          \begin{tabular}[c]{@{}c@{}} I \\ I \end{tabular} &
          \begin{tabular}[c]{@{}c@{}} 5-15 \\ 5-15 \end{tabular} &
          \begin{tabular}[c]{@{}c@{}} 0 \\ 20 \end{tabular} &
          Varies & 
          \begin{tabular}[c]{@{}c@{}} \mycitep{DeNegre2011} \end{tabular} \\
\hline
DEN2    & 110 &
          \begin{tabular}[c]{@{}c@{}} I \\ I \end{tabular} &
          \begin{tabular}[c]{@{}c@{}} 5-10 \\ 5-20 \end{tabular} &
          \begin{tabular}[c]{@{}c@{}} 0 \\ 5-15 \end{tabular} &
          Varies & 
          \begin{tabular}[c]{@{}c@{}} \mycitep{DeNegre2011} \end{tabular} \\
\hline
ZHANG   & 30 &
          \begin{tabular}[c]{@{}c@{}} B \\ I \end{tabular} &
          \begin{tabular}[c]{@{}c@{}} 50-80 \\ 70-110 \end{tabular} &
          \begin{tabular}[c]{@{}c@{}} 0 \\ 6-7 \end{tabular} &
          0.6-0.8 & 
          \begin{tabular}[c]{@{}c@{}} \mycitep{ZhaOzaBranchcut17} \end{tabular} \\
\hline
ZHANG2  & 30 &
          \begin{tabular}[c]{@{}c@{}} I \\ I \end{tabular} &
          \begin{tabular}[c]{@{}c@{}} 50-80 \\ 70-110 \end{tabular} &
          \begin{tabular}[c]{@{}c@{}} 0 \\ 6-7 \end{tabular} &
          0.6-0.8 & 
          \begin{tabular}[c]{@{}c@{}} \mycitep{ZhaOzaBranchcut17} \end{tabular} \\
\hline
FIS     & 57 &
          \begin{tabular}[c]{@{}c@{}} B \\ B \end{tabular} &
          Varies &
          Varies &
          -1 &
          \begin{tabular}[c]{@{}c@{}} MIPLIB
            \\ \mycitep{FisLjuMonSinUse18} \end{tabular} \\ 
\hline
XU      & 100  &
          \begin{tabular}[c]{@{}c@{}} I \\ IC \end{tabular} &
          \begin{tabular}[c]{@{}c@{}} 10-460 \\ 4-184 \end{tabular} &
          \begin{tabular}[c]{@{}c@{}} 10-460 \\ 4-184 \end{tabular} &
          $\approx 0$ &
          \begin{tabular}[c]{@{}c@{}} Mixed \\ \mycitep{xuwang14} \end{tabular} \\
\hline
\end{tabular}
\end{center}
\caption{The summary of data sets \label{tab:dataSetSummary}}
\end{table}
The columns of Table~\ref{tab:dataSetSummary} have the following
interpretations. The first column is the name used to denote the data set in
the rest of this paper. The second is the number of instances in the data set.
The third is the type of variables in the first and second levels (B indicates
binary, I indicates general integer and C indicates real-valued/continuous).
The fourth and fifth columns indicate the number of variables and constraints
(upper and lower limits of the range over all instances in the data set). The
sixth column indicates the degree of alignment of the objectives (a normalized
inner product between -1 and 1, with -1 indicating the objectives are exactly
opposed and 1 indicating they are exactly aligned).

There are some important properties that affect problem difficulty, especially
properties that determine whether particular classes of inequalities apply.
Of course, one very important property is whether the instance has continuous
variables, especially at the lower level. Pure integer instances with all
integer coefficients in the constraint matrices have a much wider range of
cuts available. The signs of the coefficients is also an important factor.
Instances in which the lower-level constraints are ``$\leq$'' and have all
non-negative coefficients, for example, tend to be easier, as one might
expect, since this is the case even for classical single-level MILP. Instances
that have coefficients with mixed signs are generally more difficult to solve. 

The degree of alignment of the objectives is one particularly important factor
in determining the difficulty of instances. Instances with alignments near 1
(such as those in ZHANG and ZHANG2) tend to be relatively easy. On other hand,
zero-sum instances (alignment of -1) have specialized cuts and are also easier
than general instances (interdiction problems in particular are relatively
easier to solve than general instances). An alignment of -1 also ensures that
any fractional solution to the LP relaxation must violate the optimality
constraint and that case~\ref{case:goal-cut-1-a}, as described in the
beginning of Section~\ref{sec:classes} cannot arise. As we discuss later,
fractional solutions that satisfy the optimality constraint (those in
case~\ref{case:goal-cut-1-a}) may arise quite frequently when the alignment is
not -1 and this may be an important issue in managing cut generation in
particular cases.

More detailed descriptions of each of the test sets are included in
Appendix~\ref{app:test-set}. 

\subsection{Experimental Setup \label{sec:setup}}

\paragraph{Hardware Platform.}

All computational results we report were generated on compute nodes running
the Linux (Debian 8.7) operating system with dual AMD Opteron 6128 processors
and 32 GB RAM. All experiments were run on a single core with a time limit of
3600 seconds and a memory limit of 8 GB.

\paragraph{Software Versions.} \MIBS{} \t{1.2.2} was employed for
conducting all experiments and \SYMPHONY{} \t{5.7.2} was employed as the
MILP solver for all required subproblems except where noted. Although running
times can be improved to some extent by using either Cbc~\mycitep{Cbc} or a
commercial solver (CPLEX~\mycitep{CPLEX} is the only currently supported
alternative), results presented herein indicated the achieved speedups would
not have substantively changed any of the conclusions and since we are
actively developing methodologies for solving the various subproblems more
efficiently using SYMPHONY's unique capabilities, we prefer to report results
here using SYMPHONY.

\paragraph{Parameters.} The only settings that differed between experiments
were those for turning specific classes of valid inequalities off/on and those
for specifying branching priorities. All parameters of \MIBS{} and \code{Blis}
were left at their defaults. However, we note here that in a number of cases,
default parameter settings in \MIBS{} \t{1.2} were changed from those
in \t{1.1} and some of these now over-ride the default parameters in \BLIS{}
as a result of the extensive testing done here. We report those separately
here.
\begin{itemize}
\item The pseudocost branching strategy was used to choose the best variable
among the branching candidates. We have not put much effort so far into tuning
the branching strategy (see empirical analysis of branching strategy in
Section~\ref{sec:branching-results}).
\item The generation of generic MILP inequalities by the 
Cut Generation Library, which was previously disabled is now enabled (see the
empirical analysis backing this choice in Section~\ref{sec:feas-cuts}).
\item We now do not filter cuts based on dynamism or density
(\t{scaleConFactor} = \t{denseConFactor} = $\infty$), since this sometimes results
in cuts necessary for correctness being rejected. In practice, we did not
encounter any numerical issues as a result of this, at least on these
instances.
\item We generate cuts in every iteration in which they are eligible to be
generated except when tailoff detection forces branching and do not limit the
number of cut passes (\t{cutPass} = $\infty$) or total number of
cuts \t{cutFactor} = $\infty$
\item The tailoff detection mechanism was improved in \MIBS{} \t{1.2} and is now
based on relative improvement in the bound after adding cuts. After testing,
the default cutoff was set to \t{.05}, which seems to provide the best overall
performance across all instances in this set.
\end{itemize}
Scripts for precisely replicating these experiments, raw output files from
\MIBS{}, and spreadsheets with summary statistics are all publicly available in
the Github repository of \MIBS{}~\mycitep{MIBS1.2}.

\subsection{Results \label{sec:empirical}}

We now present the summarized results of the experiments. 

\subsubsection{Branching Strategy \label{sec:branching-results}}

As we mentioned in Section~\ref{sec:branching-strat}, selection of branching
strategy in the case of MIBLP is more intricate than in MILP.
In~\mycitet{Tahernejad et al.}{TahRalDeN20}, we concluded that the
effectiveness of the branching strategy was affected mainly by the ratio of
the number of integer variables at the first and second level. In particular,
with that earlier version of \MIBS{} (version \t{1.1}), when the number of first
level integer variables was smaller than the number of second level integer
variables, performance was better when using \t{linking} branching.

\begin{figure}[tbh]
\begin{subfigure}{0.5\textwidth}
\includegraphics[height=2in]{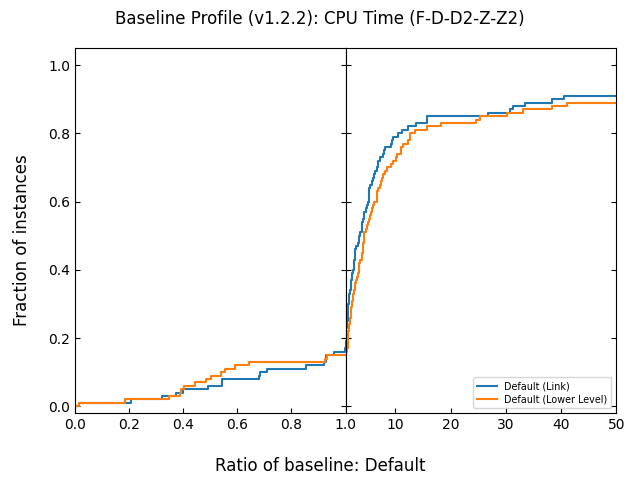}
\caption{Baseline profile of CPU time
  \label{fig:all-branching-baseline-cpu}}
\end{subfigure}
\begin{subfigure}{0.5\textwidth}
\centering
\includegraphics[height=2in]{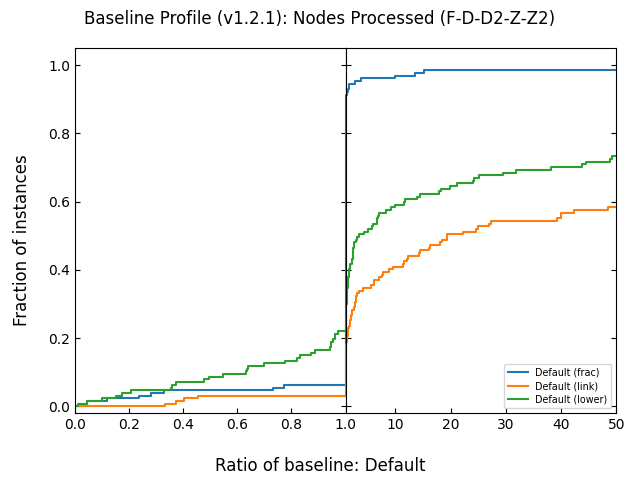}
\caption{Baseline profile of number of processed nodes
  \label{fig:all-branching-baseline-nodes}}
\end{subfigure}
\begin{subfigure}{0.5\textwidth}
\includegraphics[height=2in]{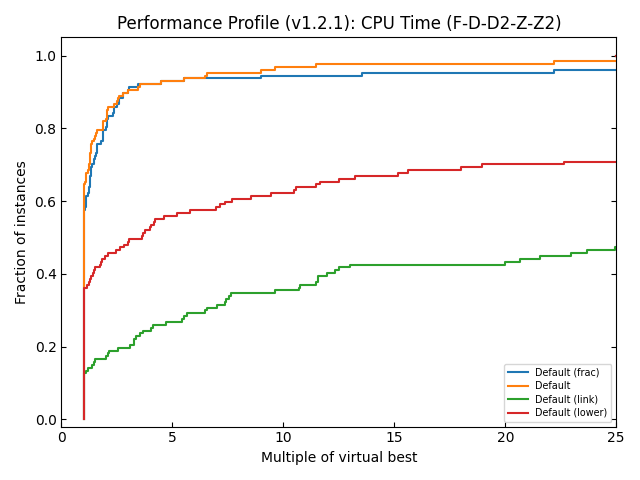}
\caption{Performance profile of CPU time
  \label{fig:all-branching-perf-cpu}}
\end{subfigure}
\begin{subfigure}{0.5\textwidth}
\includegraphics[height=2in]{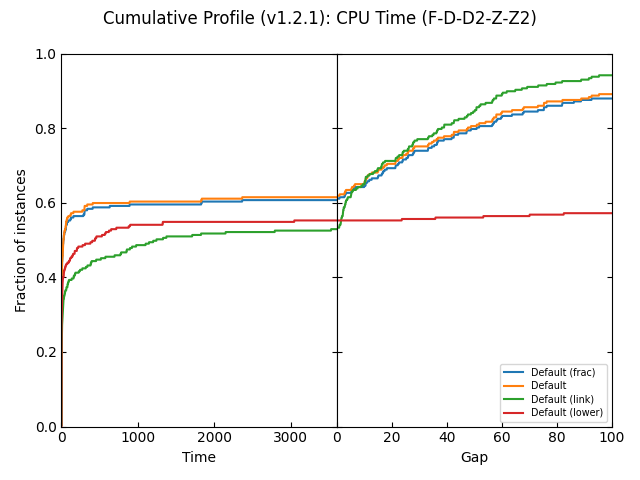}
\caption{Cumulative profile of CPU time
  \label{fig:all-branching-cum-cpu}}
\end{subfigure}
\caption{Impact of branching strategy}
\label{fig:all-branching}
\end{figure}

Here, we have revisited that result in light of the completely revamped cut
generation strategy. We have tried the three general strategies mentioned
already: branching on all fractional variables, prioritizing branching on
linking variables, and prioritizing branching on lower-level variables. Both
cumulative profiles and baseline profiles, with the baseline being
the \t{fractional} branching strategy, across all instances in
our test set (except for the interdiction ones), are shown in
Figure~\ref{fig:all-branching}.
It can be clearly seen that branching on all fractional variables is
advantageous in the vast majority cases.

These results seem to be in conflict with what was reported in our earlier
paper. In particular, branching on all fractional variables is more effective
now than it was previously and as a result, we have changed the default
branching strategy for all non-interdiction instances across the board.

We have investigated the reasons for this change and they are quite intricate.
There are a number of effects that result in complex tradeoffs that must be
navigated.  
\begin{enumerate}
\item As in MILP, the branching variable that is chosen from among the allowed
candidates is the one that is predicted to have the biggest impact in terms of
bound improvement. Restricting the set of branching candidates may reduce the
effectiveness of this strategy by removing what would have been the most
effective candidate for branching from consideration.
\item On the other hand, the logic of branching exclusively on linking
variables is to reduce the problem to a (lexicographic) MILP as quickly as
possible by fixing the values of all linking variables first. Once all linking
variables are fixed, the problem can be handed off to a subsolver. This
effectively moves the branching on second-level variables out of \MIBS{}
altogether and down into the subsolver used to solve the second-level problem,
which is better suited for the job once the problem is reduced to an MILP.
This means an increased emphasis and dependence on oracle calls and perhaps
less emphasis on cut generation. In theory, the fact that we are enumerating
over just the set of linking variables instead of all variables could result
in smaller search trees.
\item A third effect of the choice of branching strategy is that it influences
the structure of the solutions produced by solving the LP relaxation in ways
that may affect our ability to generate strong valid inequalities. Unlike in
MILP, generation of MIBLP cuts is not guaranteed to succeed when the the point
being separated is not integral, as we discussed in
Section~\ref{sec:intersection-cuts}. It can be argued that the
intersection cuts are less likely to fail when the lower-level variables have
integer values because this increases the chance that there exists an
improving direction/solution. One thing is for certain---both the branching
strategy and the separation strategy have a potentially important effect on
the failure rate for cut generation, which in turn affects the overall
efficiency of the algorithm.
\end{enumerate}
Because of these factors, the three branching strategies lead to substantively
different algorithmic behavior and different algorithm emphases. Despite the
dominance of \t{fractional} branching across almost all of our experiments, it is
difficult to extract a simple explanation that provides a universal
explanation for this result.

Table~\ref{tab:branching-stats} shows statistics averaged across all
non-interdiction, pure integer instances with only IDIC inequalities being generated
using the \t{Always} strategy. The columns have the following meanings:
\begin{itemize}
\item The first set of two columns are the averages for total tree size and
total number of oracle calls across all instances.
\item The second set of three columns shows how the CPU time breaks down
between oracle calls and cut generation, with the first column showing total
CPU time and the second and third columns showing time spent on oracle calls
and generation of MIBLP cuts, respectively (the difference between the sum of
the second and third columns and the first column is the time spent in all
other parts of the code---solving LP relaxations, branching, generating MILP
cuts, etc.---and is typically negligible).
\item The final set of three columns shows averages per node for (1) CPU time,
(2) number of calls to the MIBLP cut generator; and (3) total number of cuts
successfully generated.
\end{itemize}

\begin{table}
\begin{center}
  \begin{tabular}{|l|r|r|r|r|r|r|r|r|}
\hline
           & \multicolumn{2}{|c|}{Total Number} & \multicolumn{3}{|c|}{Total
           Time} & \multicolumn{3}{|c|}{Average Per Node} \\
\hline
\begin{tabular}{l} Branching \\ Strategy \end{tabular} & Nodes
& \begin{tabular}{l} Oracle \\ Calls \end{tabular} & CPU
& \begin{tabular}{l} Oracle \\ Calls \end{tabular} & CG   & CPU  & CG Calls & Cuts \\
\hline
Fractional & 2226  & 8      & 52   & 0.4    & 51   & .023 & 1.13     & 0.23  \\
Linking    & 24498 & 1400   & 175  & 23     & 138  & .007 & 0.53     & 0.07  \\
\hline
\end{tabular}
\end{center}
\caption{Statistics for solving with different branching rules \label{tab:branching-stats}}
\end{table}

Here, the dominance of the \t{fractional} branching can be clearly seen. Most
strikingly, the tree size is smaller by an order of magnitude, despite the
much larger search space. This is probably attributable to a combination of the
increased effectiveness of choosing from a larger set of branching candidates
and the increased effectiveness of cut generation. The number of oracle class
is negligible in the case of \t{fractional} branching, while it is substantial in
the case of \t{linking} branching. This can be attributed to dramatic differences
in the number of nodes in which all linking variables have fixed values. This
condition forces an oracle call according to current default parameter
settings in \MIBS{}. This is also reflected in the average number of rounds of
MIBLP cut generation, which is more than one with \t{fractional} branching, but
approximately a half with \t{linking} branching. And finally, the number of cuts
successfully generated and added per node is significantly higher for
\t{fractional} branching. We attribute this to a combination of there simply being
more frequent cut generation and also to a higher rate of success. It is
unclear whether the cuts generated when branching on linking variables are
actually more effective.

Before closing the section, we note that interdiction problems
(and, we conjecture, other zero-sum problem as well) are an exception to the
general rule because of the very close linking of the first- and second-level
variables. For these instances, \t{linking} branching is a clearly dominant
strategy.

\subsubsection{Cut Generation Strategy \label{sec:cut-results}}

To test the effectiveness of our cut generation, we first tested the
effectiveness of generating each type of valid inequality individually under
both of our main branching strategies and, in the case of interdiction cuts,
with the different cut generation strategies described in
Section~\ref{sec:cut-strat}. Based on these results, we also tried
combinations of multiple classes of inequalities in order to determine the
best default values for cut generating parameters.

In each figure that follows, the classes of inequalities are indicated by the following
shorthand. 
\begin{itemize}
\item\t{Benders Interdict}: The Benders interdiction cut.
\item\t{Gen No Good}: The generalized no-good cut.
\item\t{Benders Binary}: The Benders binary cut.
\item\t{Int No Good}: The integer no-good cut.
\item\t{ISIC Type1}: Improving solution ICs of type \RNum{1}. 
\item\t{ISIC Type2}: Improving solution ICs of type \RNum{2}.
\item\t{IDIC}: Improving direction ICs.
\item\t{Hypercube IC}: The hypercube IC.
\item \t{No Cuts}: No MIBLP cuts generated (integrality cuts are still generated). 
\end{itemize}
For intersection cuts, the specific strategy utilized from those described in
Section~\ref{sec:cut-strat} is also indicated as part of the shorthand
description. When \code{link} is in parentheses, the \t{linking} branching
strategy was used and similarly, when \t{frac} is in parentheses,
\t{fractional} branching was used (we primarily show results
with \t{fractional} branching since this strategy was dominant for almost all
classes). We also experimented with prioritizing branching on second-level
variables, but those results were not promising. Note that when no cuts are
generated, we are forced to branch only on linking variables (e.g., to branch
on non-fixed linking variables, even if they have integer values in the
relaxation solution).

In most of the experiments, a single class of inequality was turned on and
tested with each of the two branching strategies and possibly also with
different parameter settings, subject to applicability of each class (based on
instance structure). In the case of intersection cuts, all classes were tested
with all settings described in Section~\ref{sec:cut-strat}, with only the best
settings for each class used in the global comparison between classes. Some
testing with combinations of multiple classes generated simultaneously was
also done in order to determine if there any classes that were
``complementary'' and have an additive effect. We did find some additive
effects, described below, but we did not do extensive testing in this regard.
Because of the (very) large number of different parameters settings that were
tested, as well as the number of different test sets, we present only summary
figures in the main text. Detailed figures for each test set individually are
shown in Appendix~\ref{app:results}.

\begin{figure}[tbh]
\begin{subfigure}{0.5\textwidth}
\includegraphics[height=2.5in]{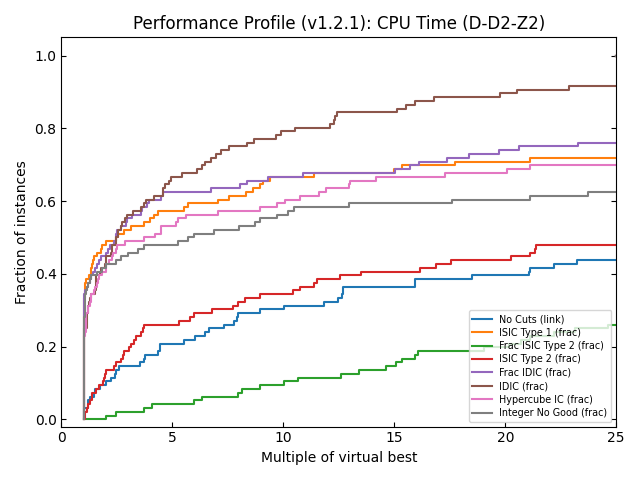}
\caption{Performance profile of CPU time
  \label{fig:pureInteger-perf-cpu}}
\end{subfigure}
\begin{subfigure}{0.5\textwidth}
\centering
\includegraphics[height=2.5in]{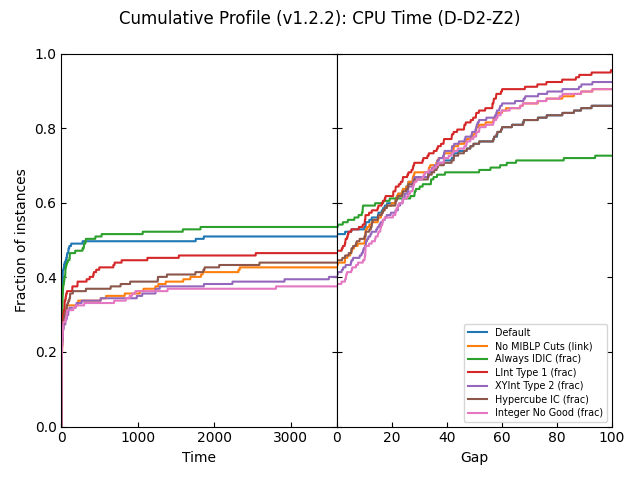}
\caption{Cumulative profile of CPU time
  \label{fig:pureInteger-cum-cpu}}
\end{subfigure}
\begin{subfigure}{0.5\textwidth}
\centering
\includegraphics[height=2.5in]{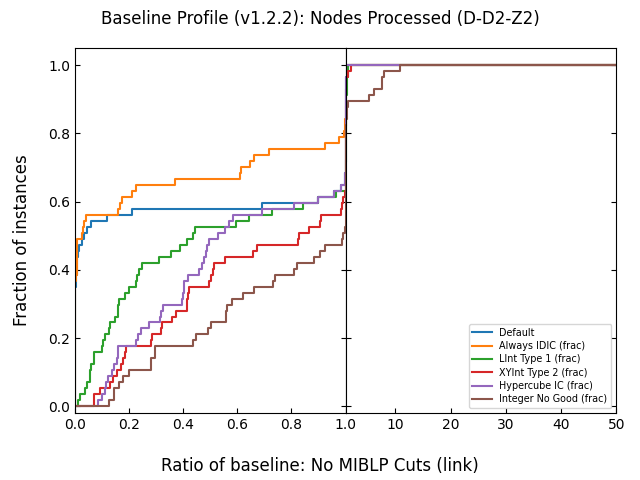}
\caption{Baseline profile of tree size
  \label{fig:pureInteger-nodes}}
\end{subfigure}
\begin{subfigure}{0.5\textwidth}
\centering
\includegraphics[height=2.5in]{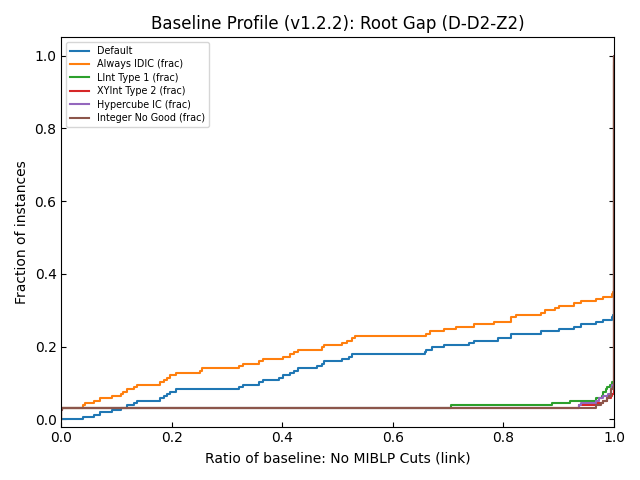}
\caption{Baseline profile of root gap
  \label{fig:pureInteger-root-gap}}
\end{subfigure}
\caption{Comparing the performance of different inequalities on the pure
integer instances}
\label{fig:pureInteger}
\end{figure}

\paragraph{Pure Integer Instances.} Figure~\ref{fig:pureInteger} shows the
results on the non-binary, pure integer instances (DEN, DEN2, and ZHANG2) for
all applicable classes of inequalities, with IDICs generated using
the \t{Always} strategy; type \RNum{1} ISICs generated with the \t{LInt}
strategy; and type \RNum{2} ISICs generated with the \t{XYInt} strategy. We
also tested with what is now the new default strategy, which uses multiple
classes and is described later. The performance profile in
Figure~\ref{fig:pureInteger-perf-cpu} shows that among the instances which
could be solved by one of the methods, generating IDICs with the \t{Always}
strategy is the clear winner in terms of CPU time and
Figure~\ref{fig:pureInteger-nodes} shows that the resulting trees are almost
universally smaller. Somewhat surprisingly, no other class on its own does
much better than no MIBLP cuts at all, with Hypercube ICs and Integer No Goods
faring even worse then the baseline with no MIBLP cuts. On the other hand,
Figure~\ref{fig:pureInteger-nodes} shows that all classes reduce the tree size
substantially over the case of no MIBLP cuts, so it seems that the issue is
that the expense of the cut generation does not pay off.

The fact that IDICs are the most effective inequalities in almost all cases is
expected in general, based on the discussion in Section~\ref{sec:comparing}.
What was less clear before running the experiments was whether aggressive
generation strategies for ICs would pay dividends. The answer to that question
seems to be a resounding yes, despite the high failure rate for generating
these cuts (see further discussion below). Note, however, that
Figure~\ref{fig:iblpDen2-perf} shows that for the single instance set
\t{IBLP-DEN2}, the IDICs vastly under-perform. The reasons for this are not known,
but it can be observed from Figure~\ref{fig:iblpDen2-root-gap} that root gap
closed by IDICs is roughly the same as that closed by Integer No Goods, which
are much cheaper to generate. This seems to indicate that IDICs are not
effective enough to be worth the computational cost in this case. These
instances differ from the others in one obvious way---the constraint matrices
have both positive and negative coefficients---and it's possible that this
makes the cut generation more expensive for this class.

An interesting observation from Figure~\ref{fig:pureInteger-cum-cpu} is that
while IDICs dominate for instances that can be solved within an hour using
IDICs, they are not as effective in general at closing the gap for more
difficult instances. ISICs are better for this. This suggests that a dynamic
scheme in which IDICs are switched off if they're not effective in the early
stages of the algorithm could have an impact.

Based on these results, our default strategy for pure integer instances
in \MIBS{} is to use a combination of IDICs with the \t{Always} strategy along
with ISICs of Type \RNum{1} with the \t{LInt} strategy. Although ISICs alone
do not seems to be effective, the combination of these two classes is superior
to any individual class. There is one exception to this default strategy and
that is for instances with constraint matrices that have both positive and
negative coefficients (like \t{IBLP-DEN2}). For now, we have chosen to use
Hypercube ICs for these instances by default, though this requires more
investigation.

\begin{figure}[tbh]
\begin{subfigure}{0.5\textwidth}
\includegraphics[height=2.5in]{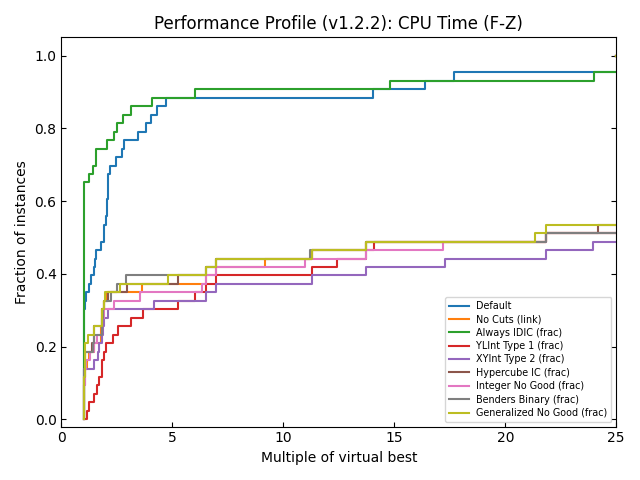}
\caption{Performance profile of CPU time
  \label{fig:binary-perf-cpu}}
\end{subfigure}
\begin{subfigure}{0.5\textwidth}
\centering
\includegraphics[height=2.5in]{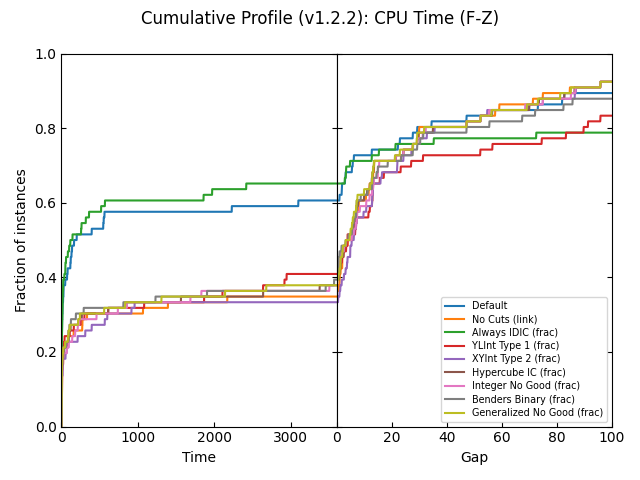}
\caption{Cumulative profile of CPU time
  \label{fig:binary-cum-cpu}}
\end{subfigure}
\begin{subfigure}{0.5\textwidth}
\centering
\includegraphics[height=2.5in]{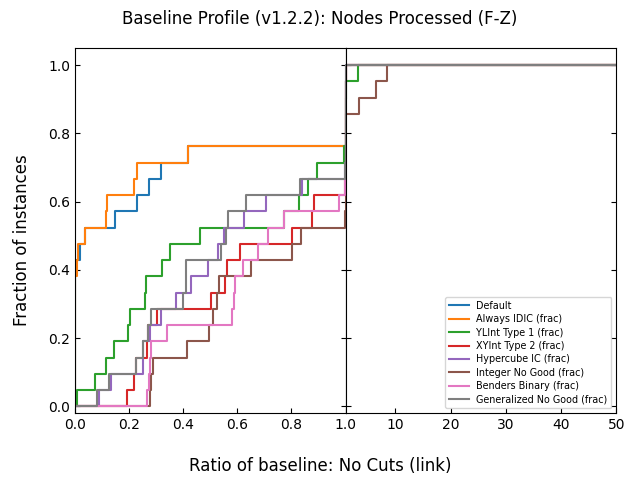}
\caption{Baseline profile of tree size
  \label{fig:binary-nodes}}
\end{subfigure}
\begin{subfigure}{0.5\textwidth}
\centering
\includegraphics[height=2.5in]{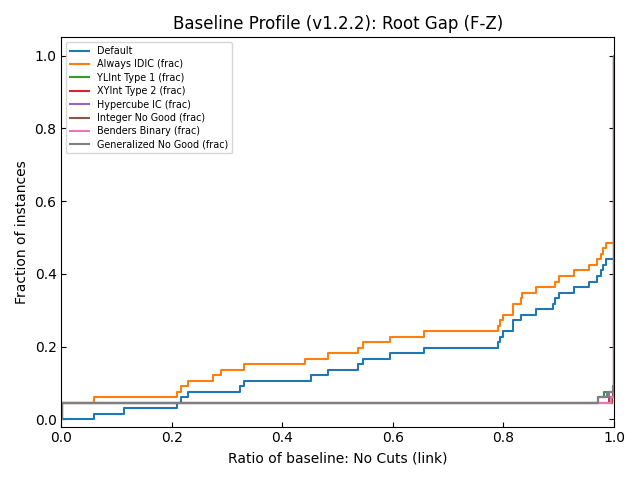}
\caption{Baseline profile of root gap
  \label{fig:binary-root-gap}}
\end{subfigure}
\caption{Comparing the performance of different inequalities on the binary instances}
\label{fig:binary}
\end{figure}

\paragraph{Binary First-Level Instances.} Figure~\ref{fig:binary} shows the
results on the instances with binary first-level variables (FIS and ZHANG) on
all applicable classes of inequalities. Here, we have tested all classes
applicable to only binary instances in addition to the previous inequalities
applicable to pure integer problems, with the only difference being that type
\RNum{1} inequalities are generated using the \t{YLInt} strategy instead of
the \t{Lint}. For these instances, the performance profile in
Figure~\ref{fig:binary-perf-cpu} again shows that among the instances that
could be solved by one of the methods, generating IDICs with the \t{Always}
strategy is the clear winner in terms of CPU time and
Figure~\ref{fig:binary-nodes} shows that the resulting trees are again almost
universally smaller.

As before, almost no other class outperforms the setting in which we generate
no MIBLP cuts in terms of CPU time, but tree size is consistently reduced. The
conclusion once again is that the expense of cut generation is not worthwhile
in many cases. It may seem surprising that the cuts specialized for the binary
case are not as effective as general MIBLP cuts, but this is in line with our
earlier analysis. This is a phenomenon worth further investigation, given the
effectiveness of methods specialized to the binary case in solving MILPs. We
observe once again from Figure~\ref{fig:binary-cum-cpu} that also for these
instances, while IDICs dominate for instances that can be solved within an
hour, they are not as effective in general at closing the gap for other
instances.

Based on these results, our default strategy for pure integer instances
in \MIBS{} is to use a combination of Benders binary cuts, generalized no-good
cuts, ISICs of Type \RNum{1} with the \t{XYInt} strategy, and IDICs with
the \t{Always} strategy. Although the specialized cuts are not as effective as
ICs on their own, they are cheap to generate and do improve performance when
generated in combination with ICs. 

\begin{figure}[tbh]
\begin{subfigure}{0.5\textwidth}
\includegraphics[height=2.5in]{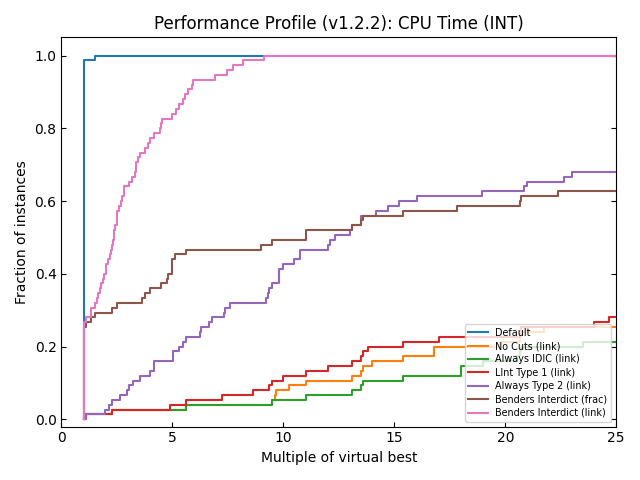}
\caption{Performance profile of CPU time
  \label{fig:interdDen-perf-cpu}}
\end{subfigure}
\begin{subfigure}{0.5\textwidth}
\centering
\includegraphics[height=2.5in]{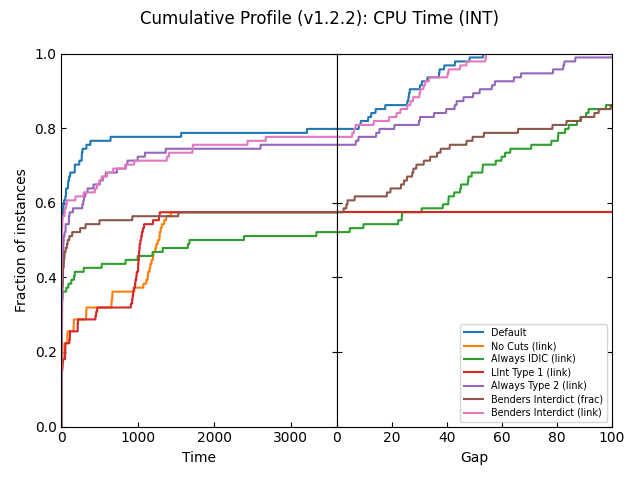}
\caption{Cumulative profile of CPU time
  \label{fig:interdDen-cum-cpu}}
\end{subfigure}
\begin{subfigure}{0.5\textwidth}
\centering
\includegraphics[height=2.5in]{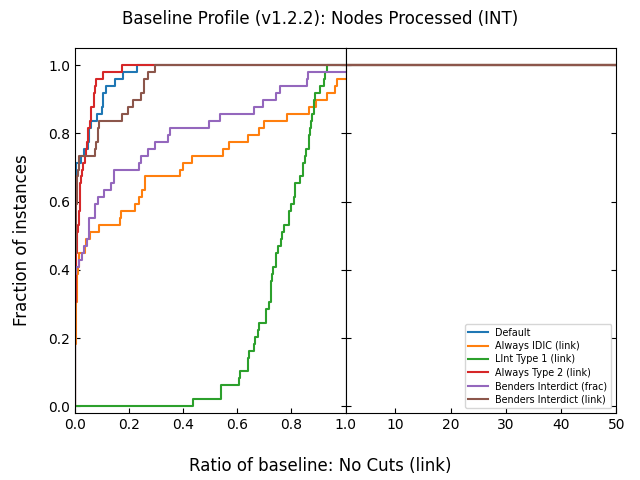}
\caption{Baseline profile of tree size
  \label{fig:interdDen-nodes}}
\end{subfigure}
\begin{subfigure}{0.5\textwidth}
\centering
\includegraphics[height=2.5in]{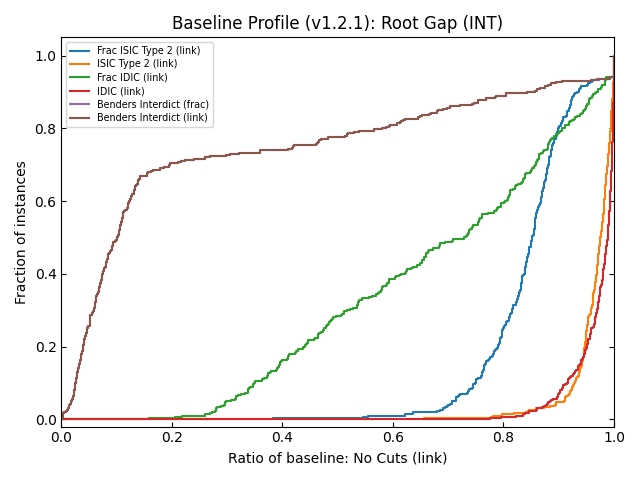}
\caption{Baseline profile of root gap
  \label{fig:interdDen-root-gap}}
\end{subfigure}
\caption{Comparing the performance of different inequalities on the \t{INT-DEN} set}
\label{fig:interdDen}
\end{figure}

\paragraph{Interdiction Instances.} Figure~\ref{fig:interdDen} shows the
results on the interdiction instances from the \t{INT-DEN} set with the
Benders Interdiction cuts, ISICs, and IDICs (other classes of cuts were
clearly dominated and not included to keep the profile from being overly
crowded). For these instances, we are also showing the difference between the
\t{fractional} and \t{linking} branching strategies to highlight that this class of
problems indeed has very different properties when it comes to selecting
branching variables, owing to the obvious connection between the linking and
second-level variables.

For these instances, the performance profile in
Figure~\ref{fig:interdDen-perf-cpu} shows that among the instances that could
be solved by one of the methods, generating Benders Interdiction cuts is the
very clear winner in all respects. The strength of these cuts has been known
for some time and this is not a surprise. \t{Linking} branching is also
clearly dominant. Figure~\ref{fig:interdDen-nodes} shows that the trees
resulting from the generation of cuts of all classes are again universally
smaller and the root gap closed is much larger with Benders Interdiction cuts
than others, as shown in Figure~\ref{fig:interdDen-root-gap}. All classes
close some amount of gap, with the exception of type \RNum{1} ISICs for which
the lower bound remains at its initial value of zero. Unlike in other cases,
the Benders Interdiction cuts are also able to close the gap effectively on
unsolved instances, as shown in Figure~\ref{fig:interdDen-cum-cpu}.

Based on these results, our default strategy for interdiction instances
in \MIBS{} is to use a combination of BendersbBinary cuts, Benders
interdiction cuts, and ISICs of type \RNum{1} with a generation strategy
of \t{Lint}. We do not generate feasibility cuts (MILP cuts) by default, as
discussed below.

\subsubsection{Strength of Inequalities}\label{sec:evaluate-strength-rhs}
We selected randomly 10 instances of the \t{INT-DEN} and \t{IBLP-DEN} sets and
solved each instance with different inequalities as shown in
Tables~\ref{tab:analysis-rhs-interDen} and~\ref{tab:analysis-rhs-iblpDen}. For
each instance, we considered the first generated cut and compared its
right-hand side with the ``best possible'' right-hand side, obtained by
optimizing over $\conv(\F)$ with the coefficient vector of the cut as
the objective function. These values are shown in the \t{Orig RHS} and \t{Best
RHS} columns, respectively. The \t{Obj before cut}, \t{Obj after orig cut}
and \t{Obj after best cut} columns show, respectively, the objective values of
the relaxation problem before adding the first cut, after adding this cut and
after adding this cut with the best right-hand side.

As one can observe in Table~\ref{tab:analysis-rhs-interDen}, the right-hand
sides of the generated Benders cuts are equal to the best possible right-hand
sides for all five instances and it is one of the reasons for the strength of
this cut. Furthermore, the right-hand sides of the generated IDICs are closer
to their best possible values comparing with the integer no-good cuts and it
verifies the results shown in Section~\ref{sec:cut-results}.

Note that in Table~\ref{tab:analysis-rhs-iblpDen}, for the instances 1, 2 and
3, we observed that the bilevel feasible region of the node in which the first
cut is generated is empty and the result is that the best value of the right-hand
side is $\infty$ (because the best right-hand side is obtained by minimizing
the left-hand side over the bilevel feasible region) and the node should be
pruned. This table also shows that the right-hand sides of the IDICs and
hypercube ICs are closer to their best values in comparison with the integer
no-good cuts and it verifies the superiority of these inequalities over the
integer no-good cut. Moreover, since the hypercube IC and generalized no-good
cut may remove a part of the non-improving bilevel feasible solutions, the
right-hand sides of these inequalities may be greater than their best possible
values and it can be observed for the instances 3 and 4.
\begin{table}[h!]
\caption{Analysis of the right-hand sides of different inequalities for the \t{INT-DEN} set}
\label{tab:analysis-rhs-interDen}
\resizebox{\columnwidth}{!}{
\begin{tabular}{cccccccccccccccc}
\hline
\multicolumn{1}{c}{}&
\multicolumn{1}{c}{}&
\multicolumn{4}{c}{Benders binary}&
\multicolumn{1}{c}{}&
\multicolumn{4}{c}{IDIC}&
\multicolumn{1}{c}{}&
\multicolumn{4}{c}{Integer no-good cut}\\
\cline{3-6}
\cline{8-11}
\cline{13-16}
Instance & \begin{tabular}[c]{@{}c@{}}Obj\\before cut\end{tabular} & \begin{tabular}[c]{@{}c@{}}Orig\\RHS\end{tabular} & \begin{tabular}[c]{@{}c@{}}Best\\RHS\end{tabular} & \begin{tabular}[c]{@{}c@{}}Obj after\\orig cut\end{tabular}&\begin{tabular}[c]{@{}c@{}}Obj after\\best cut\end{tabular} & &\begin{tabular}[c]{@{}c@{}}Orig\\RHS\end{tabular} & \begin{tabular}[c]{@{}c@{}}Best\\RHS\end{tabular} & \begin{tabular}[c]{@{}c@{}}Obj after\\orig cut\end{tabular}&\begin{tabular}[c]{@{}c@{}}Obj after\\best cut\end{tabular}& &\begin{tabular}[c]{@{}c@{}}Orig\\RHS\end{tabular} & \begin{tabular}[c]{@{}c@{}}Best\\RHS\end{tabular} &\begin{tabular}[c]{@{}c@{}}Obj after\\orig cut\end{tabular}&\begin{tabular}[c]{@{}c@{}}Obj after\\best cut\end{tabular}  \\
\hline 
1 & 0 & 4520 & 4520 & 424.64 & 424.64 & & 1 & 1.64 & 0 & 514.43 && 1 & 7 & 0 & 0 \\
2 & 0 & 6252 & 6252 & 1270.08 & 1270.08 & & 1 & 1.11 & 0 & 407.73 && 1 & 12 & 0 & 659.91 \\
3 & 0 & 5662 & 5662 & 1285.7 & 1285.7 & & 1 & 1.07 & 0 & 61.21 & & 1 & 14 & 0 & 481.27 \\
4 & 0 & 7113 & 7113 & 1274.73 & 1274.73 & & 1 & 1.15 & 0 & 390.19 & & 1 & 16 & 0 & 701.53 \\
5 & 0 & 9685 & 9685 & 1720.64 & 1720.64 & & 1 & 1.11 & 0 & 406.65 & & 1 & 18 & 0 & 1070.59 \\
\hline
\end{tabular}
}
\end{table}

\begin{table}[h!]
\caption{Analysis of the right-hand sides of different cuts for the \t{IBLP-DEN} set}
\label{tab:analysis-rhs-iblpDen}
\resizebox{\columnwidth}{!}{
\begin{tabular}{cccccccccccccccc}
\hline
\multicolumn{1}{c}{}&
\multicolumn{1}{c}{}&
\multicolumn{4}{c}{IDIC}&
\multicolumn{1}{c}{}&
\multicolumn{4}{c}{Hypercube IC}&
\multicolumn{1}{c}{}&
\multicolumn{4}{c}{Integer no-good cut}\\
\cline{3-6}
\cline{8-11}
\cline{13-16}
Instance & \begin{tabular}[c]{@{}c@{}}Obj\\before cut\end{tabular} & \begin{tabular}[c]{@{}c@{}}Orig\\RHS\end{tabular} & \begin{tabular}[c]{@{}c@{}}Best\\RHS\end{tabular} & \begin{tabular}[c]{@{}c@{}}Obj after\\orig cut\end{tabular}&\begin{tabular}[c]{@{}c@{}}Obj after\\best cut\end{tabular} & &\begin{tabular}[c]{@{}c@{}}Orig\\RHS\end{tabular} & \begin{tabular}[c]{@{}c@{}}Best\\RHS\end{tabular} & \begin{tabular}[c]{@{}c@{}}Obj after\\orig cut\end{tabular}&\begin{tabular}[c]{@{}c@{}}Obj after\\best cut\end{tabular}& &\begin{tabular}[c]{@{}c@{}}Orig\\RHS\end{tabular} & \begin{tabular}[c]{@{}c@{}}Best\\RHS\end{tabular} &\begin{tabular}[c]{@{}c@{}}Obj after\\orig cut\end{tabular}&\begin{tabular}[c]{@{}c@{}}Obj after\\best cut\end{tabular}  \\
\hline 
1 & -669 & 0 & $\infty$ & -623.47  & prune & & 5 & $\infty$ & -645 & prune & & -789 & $\infty$ & -668.5 & prune \\
2 & -688 & 0 & $\infty$ & -544 & prune & & 1 & $\infty$ & -632.46 & prune & & -29 & $\infty$ & -676.36 & prune \\
3 & -766 & 0 & 0 & -286 & -286 & & 3 & 2 & -761.08 & -766 & & -17 & -4 & -761.08& -508.05 \\
4 & -724 & 0 & 0 & -578.63 & -578.63 & & 19 & 18 & -709.75 & -724 & & -407 & -402 & -721 & -706 \\
5 & -659 & 0 & $\infty$ & -656.67 & prune & & 15 & $\infty$ & prune & prune & & -603 & $\infty$ & -656.83 & prune\\ 
\hline
\end{tabular}
}
\end{table} 

\subsubsection{Cut Generation Failures \label{sec:cg-failures}}

As mentioned in Section~\ref{sec:intersection-cuts}, it is possible for cut
generation to fail in the case of both ISICs and IDICs for several reasons.
Here, we provide a brief analysis of how prevalent this is.
Table~\ref{tab:cut-failure} shows overall statistics for each of the test sets
in a variety of scenarios.
\begin{itemize}
\item Lines 1 to 5 in the table are for the cut generation
strategy \t{Always}. The first three of these five lines are the results
when only generating IDICs with both \t{fractional} branching and \t{linking}
branching. The second line shows the effect of disabling the
generation of MILP cuts. Lines 4 and 5 show the results when generating ISICs
of types \RNum{1} and \RNum{2} with the \t{fractional} branching strategy.
\item Lines 6 to 8 show the effect of generating cuts only when the linking
variables are integer (the strategy \t{LInt}) and with the \t{fractional}
branching strategy, with different types of intersection cuts. 
\item Lines 9 to 11 are similar to lines 6 to 8, but with the cut generation
strategy \t{XYInt}.
\item Finally, lines 12 and 13 are both for the default strategy, but one
details failure rates for IDICs and the other for ISIC of type \RNum{1}.  
\end{itemize}
The statistics provide some evidence for several conjectured behaviors,
though more study is needed. First, we conjectured high rates of failure for
ICs for problems in which the objective alignment is near 1 and this seems
to be the case, as per the results on lines 1 to 3 for ZHANG and ZHANG2. We
also conjectured that generation of MILP cuts should help and this also seems
to be the case.

On the other hand, we conjectured that there would be a low failure rate when
the objective alignment is -1. As per lines 1 to 3 for the interdiction
problems, this is not the case. In fact, this seems to be because it can very
easily be the case that the problem of producing an improving direction is
infeasible. For binary problems, anytime all of the free variables have
fractional values in the solution to the relaxation, it is clear there can be
no improving direction (recall that the direction needs to be both feasible
and have integral coefficients). This can very easily happen once some
branching has occurred. The high failure rate indicates that is indeed an
issue. Fortunately, such infeasibility should be easily detected in
pre-processing and should not greatly affect running times.

As expected, failure rates are reduced substantially with the generation
strategy \t{LInt}, but as we have seen, although the lower rates of generation
do reduce time spent generating inequalities in most cases, this does not pay
off in general, as many fewer cuts are successfully generated in the end. As
expected, failure rates are close to zero with the strategy \t{XYInt}, which
was the only strategy available in \MIBS{} \t{1.1} (there can still be some
failures when the intersection with the radial cone and the BFS is empty).

This raises the tricky issue, which remains to be addressed, of the fact that
it is difficult to identify whether a given fractional point is in
case~\ref{case:goal-cut-1-a} or not without trying to separate it. This
separation is expensive and will fail if the point is in
case~\ref{case:goal-cut-1-a}. Hence, we are not only unable to apply strong
cuts, but we are also wasting CPU time with failed separation attempts.

\begin{table}
\begin{center}
  \begin{tabular}{|l|r|r|r|r|r|r|}
    \hline
 Scenario & DEN & DEN2 & FIS & ZHANG & ZHANG2 & INT \\
\hline
Always IDIC (frac)           & .56    & .38    & .53  & .65    & .63    & .46    \\
Always IDIC w/o MILP cuts (frac)  & .57    & .37    & .53  & .78    & .85    & .82    \\
Always IDIC (link)             & .60    & .28    & .75  & .92    & .87    & .50    \\
Always Type 1 (frac)         & .18    & .15    & .20  & .04    & .18    & .01    \\
Always Type 2 (frac)         & .18    & .27    & .45  &  .03   & .14    & .06    \\
LInt IDIC (frac)             & .18    & .16    & .13  & .31    & .38    & .01    \\
LInt Type 1 (frac)           & .07    & .06    & 0    & .03    & .11    & 0      \\
LInt Type 2 (frac)           & .05    & .14    & .14  & .02    & .11    & $<.01$ \\
XYInt IDIC (frac)            & $<.01$ & $<.01$ & .03  & $<.01$ & 0      & $<.01$ \\
XYInt Type 1 (frac)          & 0      & 0      & 0    & 0      & 0      & 0      \\
XYInt Type 2 (frac)          & $<.01$ & .01    &  .03 & $<.01$ & $<.01$ & $<.01$ \\
Default (IDIC)               & .61    & N/A    & .57  & .74    & .67    & N/A \\
Default (ISIC)               & $<.01$ & N/A    & 0    & .05    & .03    & N/A \\
\hline
\end{tabular}
\end{center}
\caption{Cut generation failure rates \label{tab:cut-failure}}
\end{table}

\subsubsection{Effect of Integrality Cuts \label{sec:feas-cuts}}

\begin{figure}[tbh]
\begin{subfigure}{0.5\textwidth}
\centering
\includegraphics[height=1.5in]{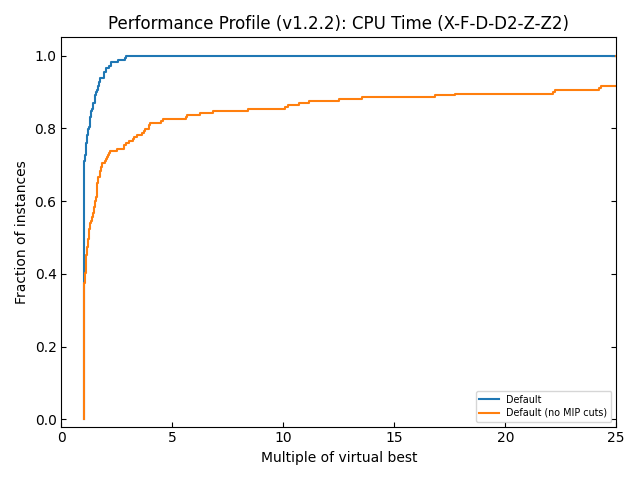}
\caption{Performance profile (all instances) \label{fig:perf_cpu_MIPvNoMIP}}
\end{subfigure}
\begin{subfigure}{0.5\textwidth}
\centering
\includegraphics[height=1.5in]{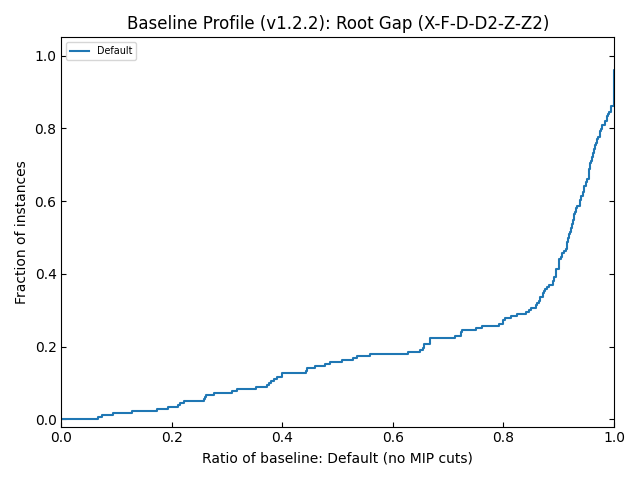}
\caption{Baseline profile  (all instances) \label{fig:base_root_gap_MIPvNoMIP}}
\end{subfigure}
\begin{subfigure}{0.5\textwidth}
\centering
\includegraphics[height=1.5in]{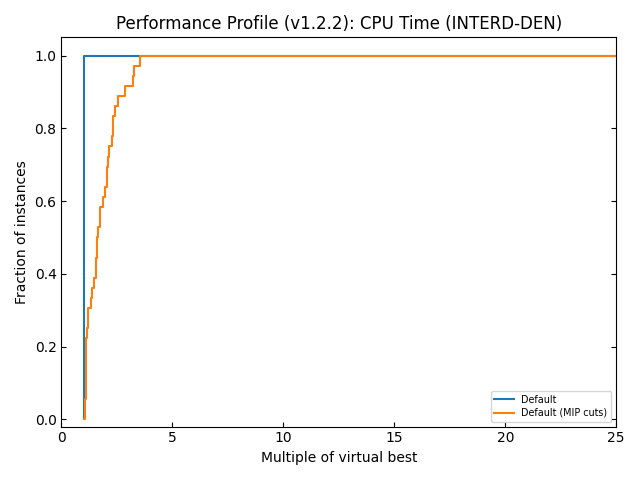}
\caption{Performance profile (interdiction) \label{fig:perf_cpu_MIPvNoMIP_INT}}
\end{subfigure}
\begin{subfigure}{0.5\textwidth}
\centering
\includegraphics[height=1.5in]{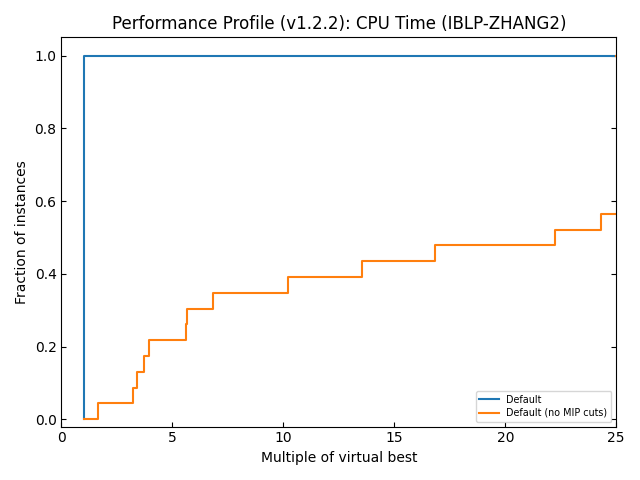}
\caption{Performance profile (ZHANG2) \label{fig:perf_cpu_MIPvNoMIP_ZHANG2}}
\end{subfigure}
\caption{Demonstrating the impact of generating MILP cuts in combination with
MIBLP cuts}
\label{fig:MIPvNoMIP}
\end{figure}

In \MIBS{} \t{1.1}, we decided not to experiment with MILP cuts because we
were focused on MIBLP-specific cuts that we felt would be much stronger. In
developing \MIBS{} \t{1.2} and with a better understanding of the behavior of
the MIBLP cuts, we decided to bring MILP cuts back into the picture. Our
empirical testing has shown that generating MILP cuts is indeed impactful. Our
best understanding of this is that the use of MILP cuts encourages more
integer solutions to be produced and this in turn reduces failure rates
arising because of being in case~\ref{case:goal-cut-1-a} (which can only
happen with fractional solutions), resulting in more frequent generation of
strong ICs. This is born out in the performance profile shown in
Figure~\ref{fig:perf_cpu_MIPvNoMIP}. The baseline profile in
Figure~\ref{fig:base_root_gap_MIPvNoMIP} shows that integrality (MILP) cuts
close a significant amount of gap above what MIBLP cuts already close.

We have also conjectured that MILP cuts should help to a greater extent for
instances with objective alignment near 1 and may not help as much for
instances with objective alignment -1.
Figures~\ref{fig:perf_cpu_MIPvNoMIP_INT}
and~\ref{fig:perf_cpu_MIPvNoMIP_ZHANG2} demonstrate that for the test sets
ZHANG2 and INT-DEN, this is indeed the case. 

\subsubsection{Comparisons to Other Solvers}

\begin{figure}
\begin{subfigure}{0.5\textwidth}
\centering
\includegraphics[height=1.7in]{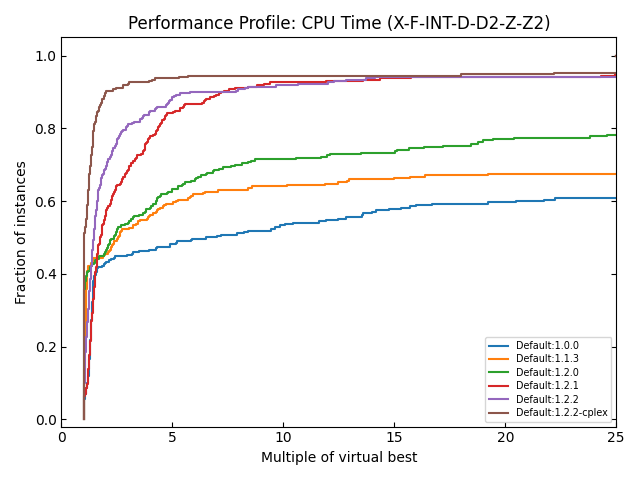}
\caption{Performance profile \label{fig:comparison-MibS}}
\end{subfigure}
\begin{subfigure}{0.5\textwidth}
\centering
\includegraphics[height=1.7in]{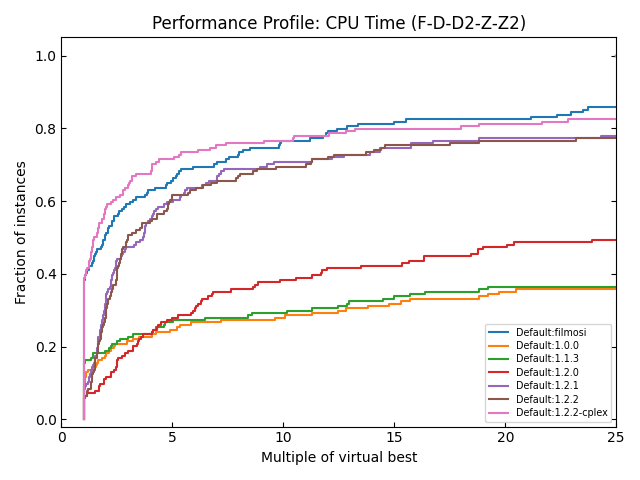}
\caption{Cumulative profile \label{fig:comparison-filmosi}}
\end{subfigure}
\caption{Overall comparison of different version of \MIBS{} to \t{filmosi}}
\label{fig:comparison}
\end{figure}

Finally, we present some overall comparisons that show both the evolution of
\MIBS{} over time and a comparison to \t{filmosi}, the most competitive
alternative solver available~\cite{filmosi}. Figure~\ref{fig:comparison-MibS} show the
comparison of different versions of \MIBS{} over time, including a version
built with CPLEX 12.7 as the subsolver. The impact of the improvements made
between \MIBS{} \t{1.1} and \t{1.2} can be clearly seen. The use of CPLEX
instead of SYMPHONY provides only a modest boost.

Figure~\ref{fig:comparison-filmosi} adds the \t{filmosi} solver to the mix.
The \t{filmosi} solver was also built on top CPLEX 12.7 and uses CPLEX not
only as a subsolver but also uses CPLEX's own branch-and-cut as the base for
the entire algorithm. As such, it is using the powerful engine of CPLEX's
MILP solver to do everything from pre-processing to cut generation to making
branching decisions. Nevertheless, through a combination of all of the
methodologies described in this paper, we are able to achieve performance more
or less on par with \t{filmosi}. It is an interesting question whether \t{filmosi}
could itself be improved by some of the tweaks we have made in \MIBS{} or
whether a version of
\MIBS{} built on top of the CPLEX branch and cut rather than that of \BLIS{} would
be possible and could be effective. 

\section{Conclusions \label{sec:conclusions}} 

Although the generation of valid inequalities has proven to be an invaluable
tool in solving MIBLPs, there remains a large scope for developing the
fundamental theory underlying the branch-and-cut algorithms that have become a
standard solution method. We have taken an initial step towards filling in
some of the gaps, as well as providing a framework within which to view the
known methods, clarifying the relationship to and distinction from similar
methods already studied in the MILP case. We also introduced new general
methodology for generating valid inequalities that applies not only to MIBLPs,
but also potentially to other classes of optimization problem (such as MILPs
themselves). Finally, we provided a comprehensive comparison of known methods
for the generation of valid inequalities. There remains a wide range of
possibilities for innovation in leveraging the decades of research solving
MILPs to derive improved methodology for MIBLPs. It is our hope that this work
motivates future work along these lines, leading to further discoveries in
this important field of research.

There also remain many additional open questions regarding how to efficiently
incorporate the generation of valid inequalities into a branch-and-cut
algorithm. Although one can simply add cut generation to an existing MILP
solver and obtain a solver for MIBLPs in a relatively straightforward way,
many of the control mechanisms that have been well-honed for solving MILPs do
not seem to work as well when naively utilized for solving MIBLPs. The work of
re-designing many of these control mechanisms in a way that works well for
MIBLPs is a rich source of future challenges in this space.

\section*{Acknowledgements}

This research was made possible with support from National Science
Foundation Grants CMMI-1435453, CMMI-0728011, and ACI-0102687, as well
as Office of Naval Research Grant N000141912330.


\appendix

\section{Detailed Descriptions of Test \label{app:test-set}}

Here are detailed descriptions of the test instances.

\begin{itemize}

\item\t{INT-DEN}: This set was generated by~\mycitet{DeNegre}{DeNegre2011} 
and contains 300 knapsack interdiction problems. These problems originate from
the \emph{Multiple Criteria Decision Making library}~\mycitep{mcdm} and has
the same structure as~\eqref{eqn:MIPINT}. The number of first-level variables
($n_1 = n_2$) varies in $n_1\in\left\{10,11,...,19,20,30,40,50\right\}$ and
the number of first- and second-level constraints are 1 and $n_1+1$,
respectively. There are 20 instances corresponding to each level of $n_1$,
except 40 instances for $n_1\in\left\{10,20\right\}$. Due to the difficulty of
the instances with $n_1 = 50$, we excluded them in the experiments.

\item\t{IBLP-DEN}: This set was generated using the publicly available
generator of~\mycitet{DeNegre}{MIBLPGEN} that were used in support of
experiments in~\mycitet{DeNegre}{DeNegre2011}. It contains 50 pure integer
instances with all variables non-negative and upper bounded by 1500. The
number of upper- and lower-level constraints is 0 and 20, respectively. The
constraints are expressed in ``$\geq$'' form with negative coefficients having
absolute value at most 50. Thus, all second-level variables likely have
implicit upper bounds less than the given bound of 1500. The number of upper-
and lower-level variables are as follows. 
\begin{center}
\begin{tabular}{l l l}
\# of instances & $n_1$ & $n_2$ \\
\hline
10 &  5 & 10 \\
10 & 10 & 10 \\
10 & 15 &  5 \\
20 & 15 &  5 \\
\end{tabular}
\end{center}

\item\t{IBLP-DEN2}: This set was generated using the publicly available
generator of~\mycitet{DeNegre}{MIBLPGEN} that were used in support of
experiments in~\mycitet{DeNegre}{DeNegre2011}. It contains 110 pure integer
instances with all variables non-negative and upper bounded by 500. The number
of upper- and lower-level constraints is 0 and the number of lower-level
variables is as indicated in the table below. The constraints are expressed in
``$\leq$'' form with coefficients of mixed sign having absolute value at most
50. Thus, it is once again likely that all second-level variables have
implicit upper bounds less than the given bound of 500. The number of
lower-level constrains and upper- and lower-level variables is as follows.
\begin{center}
\begin{tabular}{l l l l}
\# of instances & $m$ & $n_1$ & $n_2$ \\
\hline
10 & 5 & 5 & 5 \\
10 & 10 & 10 & 10 \\
10 & 10 & 15 & 10 \\
10 & 10 & 20 & 10 \\
10 & 10 & 15 & 15 \\
10 & 10 & 10 & 20 \\
10 & 15 & 10 & 10 \\
10 & 15 & 15 & 10 \\
10 & 15 & 20 & 10 \\
10 & 15 & 15 & 15 \\
10 & 15 & 10 & 20 \\
\end{tabular}
\end{center}

\item\t{IBLP-ZHANG}: This set was generated by~\mycitet{Zhang and
Ozaltin}{ZhaOzaBranchcut17}  
and includes 30 instances with binary first-level variables, integer
second-level variables and no first-level constraints. The number of
first-level variables varies in $n_1\in\left\{50, 60, 70, 80, 90\right\}$ and
the number of second-level variables is set to $n_2 = n_1 + 20$. There are 3
instances with $m_2 = 6$ and 3 instances with $m_2 = 7$ corresponding to
each level of $n_1$. The linear constraints are expressed in ``$\leq$'' form with
positive coefficients on both upper and lower level variables ranging between
0 and 10 while the right-hand side is in the range of 20 to 30. Thus, all
second-level variables have small implicit upper bounds. 

\item\t{IBLP-ZHANG2}: This set is a modification of \t{IBLP-ZHANG} in which the
bounds on first-stage variables are increased from 1 to 10 to make them
non-binary.

\item\t{IBLP-FIS}: This set was generated by~\mycitet{Fischetti et
al.}{FisLjuMonSinUse18}  
and are constructed from well-known instances originating in \t{MILPLIB
3.0}~\mycitep{bixbyceria98} in which all variables are binary and there are no
first-level constraints. The variables are split between upper and lower level
in proportions reflected in the file name of each converted instance. Due to
memory limitations, we did not consider 3 instances of this set in our
experiments. 

\item \t{MIBLP-XU}: \mycitet{Xu and Wang}{xuwang14} generated a set of mixed
integer bilevel instances for the testing of their algorithm. The instances
have $n=n_1=n_2 \in \{10, 60, 110, 160, 210, 260, 310, 360, 410, 460\}$.
The first-level variables are constrained to be integer but some of the
second-level variables are continuous. The number of first-level as well as
second-level constraints is $0.4 n$. All matrices, vectors, right-hand
sides, etc. are uniformly distributed integers.
The constraint matrices have coefficients in $[0, 10]$, while
the objective vectors have entries in $[-50, 50]$. The first-level right-hand
side vector $b_1$ has entries in $[30, 130]$, and the second-level right-hand
side vector $b_2$ has entries in $[10, 110]$.
There are 10~instances for every value of $n$, which gives a total of
100~instances.
\end{itemize}

\section{Results for Individual Test Sets \label{app:results}}

\begin{figure}[tbh]
\begin{subfigure}{0.5\textwidth}
\includegraphics[height=2.5in]{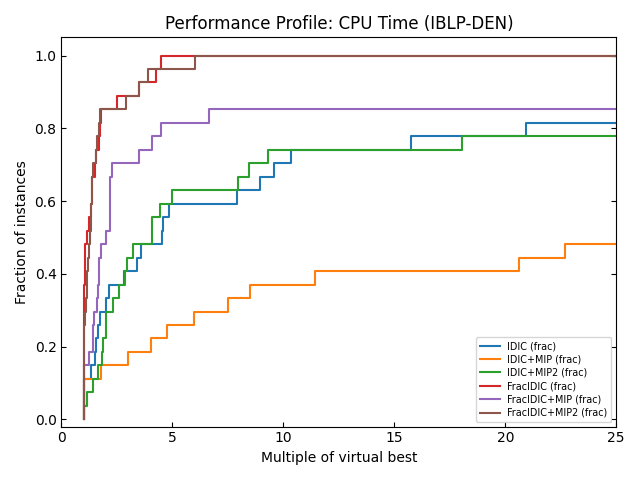}
\caption{Performance profile of CPU time
  \label{fig:iblpDen-perf-cpu}}
\end{subfigure}
\begin{subfigure}{0.5\textwidth}
\centering
\includegraphics[height=2.5in]{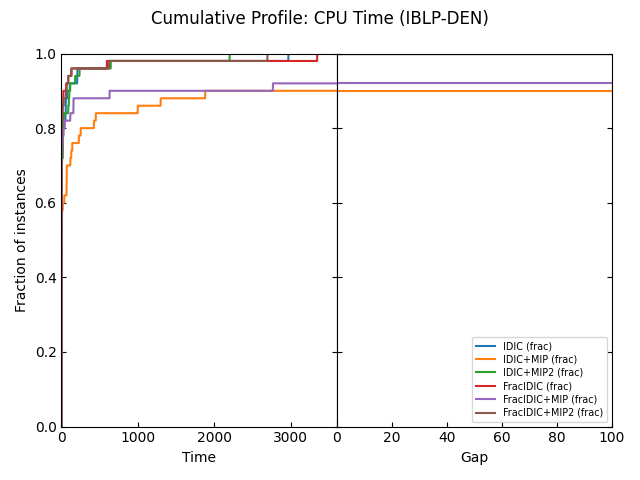}
\caption{Cumulative profile of CPU time
  \label{fig:iblpDen-cum-cpu}}
\end{subfigure}
\begin{subfigure}{0.5\textwidth}
\centering
\includegraphics[height=2.5in]{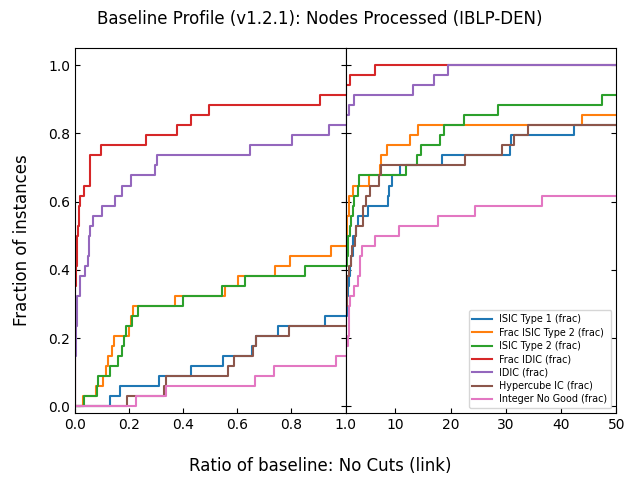}
\caption{Baseline profile of tree size 
  \label{fig:iblpDen-nodes}}
\end{subfigure}
\begin{subfigure}{0.5\textwidth}
\centering
\includegraphics[height=2.5in]{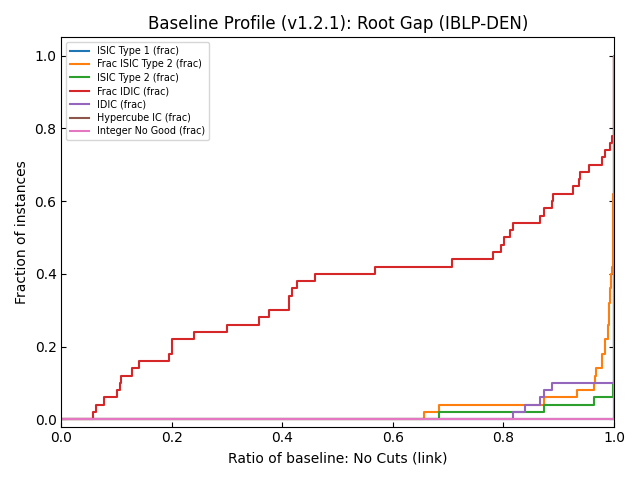}
\caption{Baseline profile of root gap
  \label{fig:iblpDen-root-gap}}
\end{subfigure}
\caption{Comparing the performance of different inequalities on the \t{IBLP-DEN} set}
\label{fig:iblpDen-perf}
\end{figure}

\clearpage

\begin{figure}[tbh]
\begin{subfigure}{0.5\textwidth}
\includegraphics[height=2.5in]{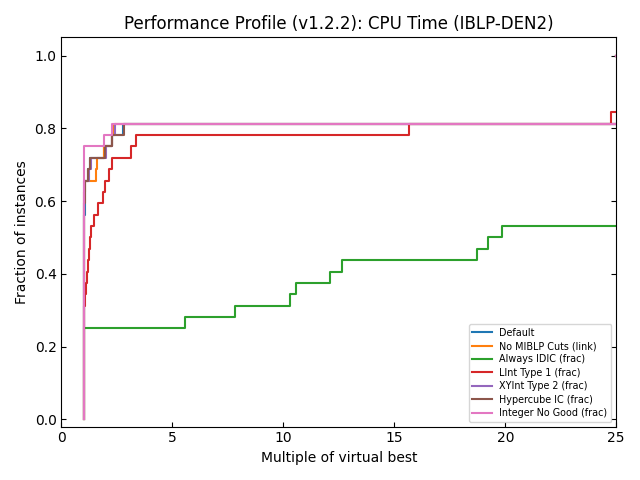}
\caption{Performance profile of CPU time
  \label{fig:iblpDen2-perf-cpu}}
\end{subfigure}
\begin{subfigure}{0.5\textwidth}
\centering
\includegraphics[height=2.5in]{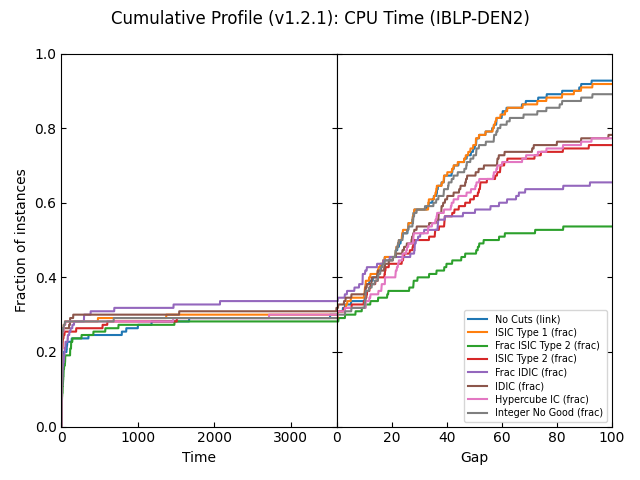}
\caption{Cumulative profile of CPU time
  \label{fig:iblpDen2-cum-cpu}}
\end{subfigure}
\begin{subfigure}{0.5\textwidth}
\centering
\includegraphics[height=2.5in]{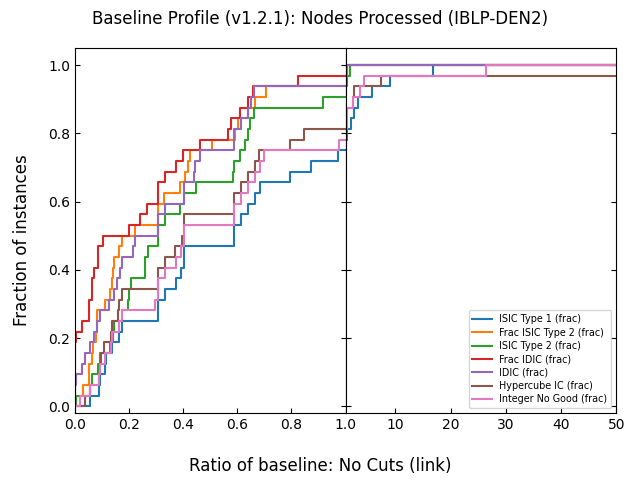}
\caption{Baseline profile of tree size 
  \label{fig:iblpDen2-nodes}}
\end{subfigure}
\begin{subfigure}{0.5\textwidth}
\centering
\includegraphics[height=2.5in]{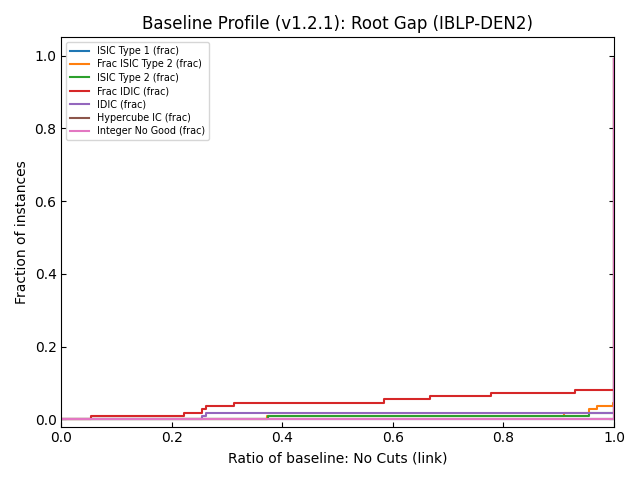}
\caption{Baseline profile of root gap
  \label{fig:iblpDen2-root-gap}}
\end{subfigure}
\caption{Comparing the performance of different inequalities on the \t{IBLP-DEN2} set}
\label{fig:iblpDen2-perf}
\end{figure}

\clearpage

\begin{figure}[tbh]
\begin{subfigure}{0.5\textwidth}
\includegraphics[height=2.5in]{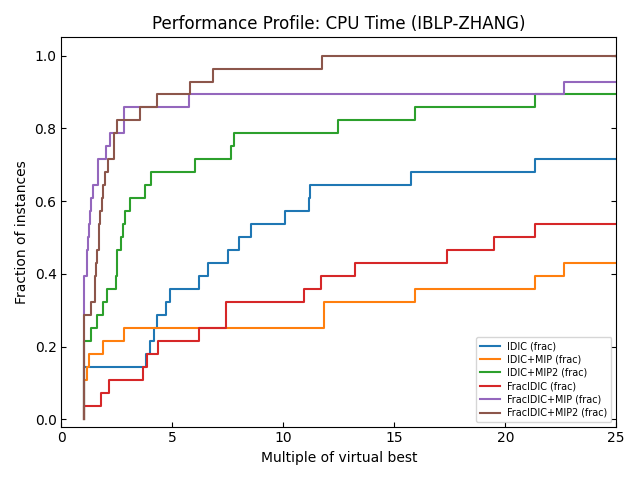}
\caption{Performance profile of CPU time
  \label{fig:iblpZhang-perf-cpu}}
\end{subfigure}
\begin{subfigure}{0.5\textwidth}
\centering
\includegraphics[height=2.5in]{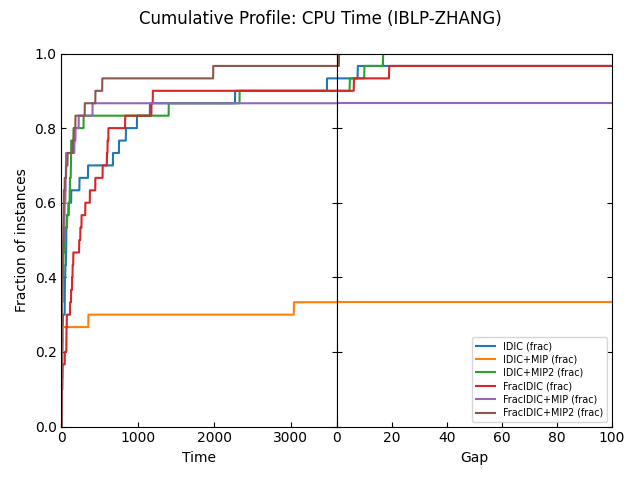}
\caption{Cumulative profile of CPU time
  \label{fig:iblpZhang-cum-cpu}}
\end{subfigure}
\begin{subfigure}{0.5\textwidth}
\centering
\includegraphics[height=2.5in]{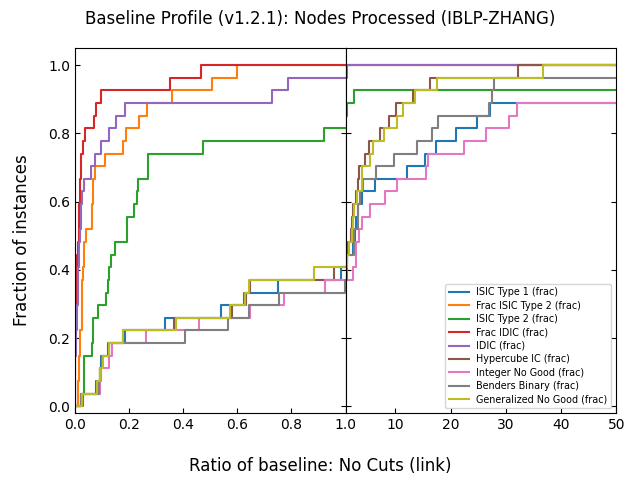}
\caption{Baseline profile of tree size 
  \label{fig:iblpZhang-nodes}}
\end{subfigure}
\begin{subfigure}{0.5\textwidth}
\centering
\includegraphics[height=2.5in]{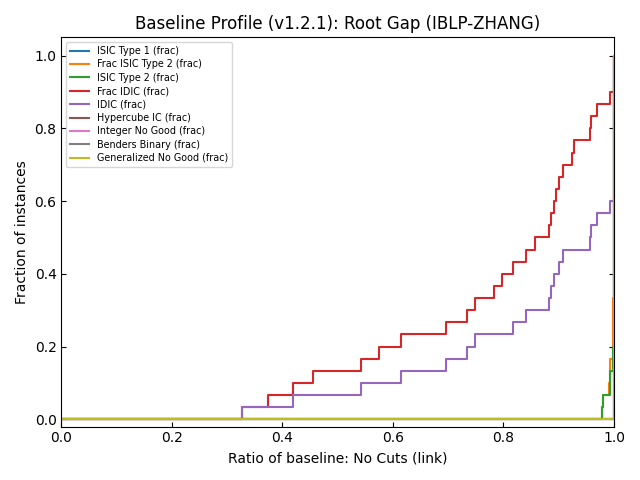}
\caption{Baseline profile of root gap
  \label{fig:iblpZhang-root-gap}}
\end{subfigure}
\caption{Comparing the performance of different inequalities on the \t{IBLP-ZHANG} set}
\label{fig:iblpZhang-perf-prof}
\end{figure}

\clearpage

\begin{figure}[tbh]
\begin{subfigure}{0.5\textwidth}
\includegraphics[height=2.5in]{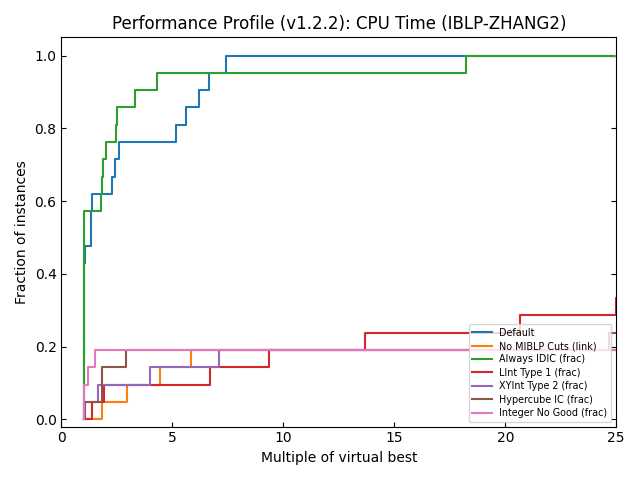}
\caption{Performance profile of CPU time
  \label{fig:iblpZhang2-perf-cpu}}
\end{subfigure}
\begin{subfigure}{0.5\textwidth}
\centering
\includegraphics[height=2.5in]{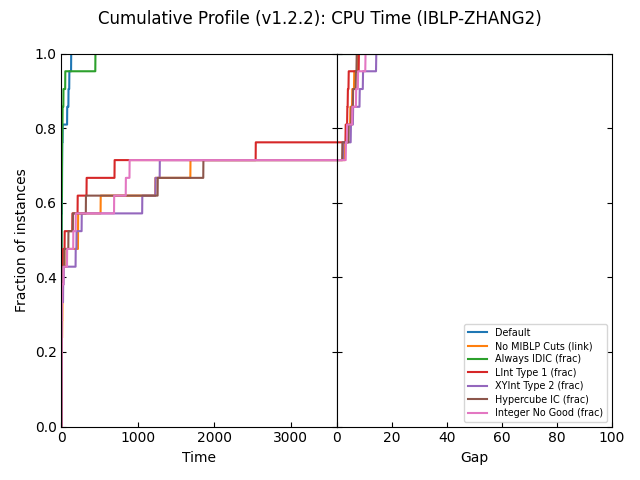}
\caption{Cumulative profile of CPU time
  \label{fig:iblpZhang2-cum-cpu}}
\end{subfigure}
\begin{subfigure}{0.5\textwidth}
\centering
\includegraphics[height=2.5in]{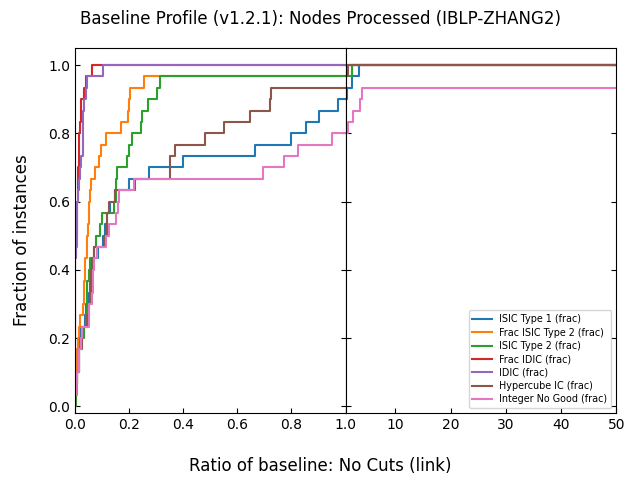}
\caption{Baseline profile of tree size 
  \label{fig:iblpZhang2-nodes}}
\end{subfigure}
\begin{subfigure}{0.5\textwidth}
\centering
\includegraphics[height=2.5in]{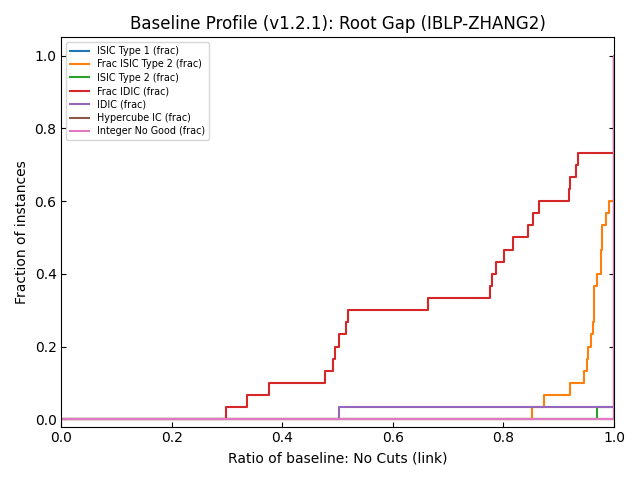}
\caption{Baseline profile of root gap
  \label{fig:iblpZhang2-root-gap}}
\end{subfigure}
\caption{Comparing the performance of different inequalities on the \t{IBLP-ZHANG2} set}
\label{fig:iblpZhang2-perf-prof}
\end{figure}

\clearpage

\begin{figure}[tbh]
\begin{subfigure}{0.5\textwidth}
\includegraphics[height=2.5in]{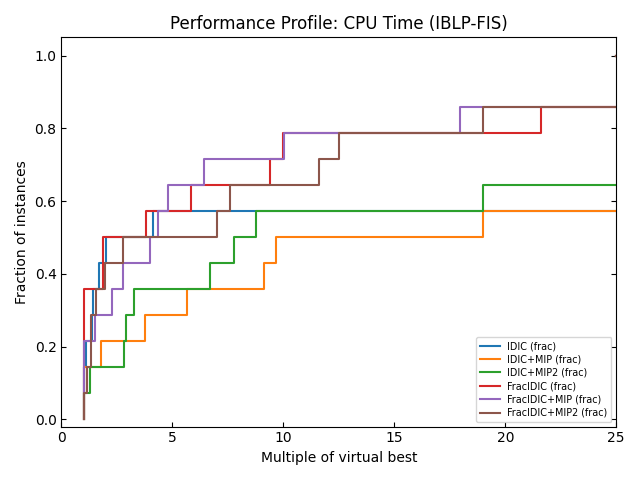}
\caption{Performance profile of CPU time
  \label{fig:iblpFis-perf-cpu}}
\end{subfigure}
\begin{subfigure}{0.5\textwidth}
\centering
\includegraphics[height=2.5in]{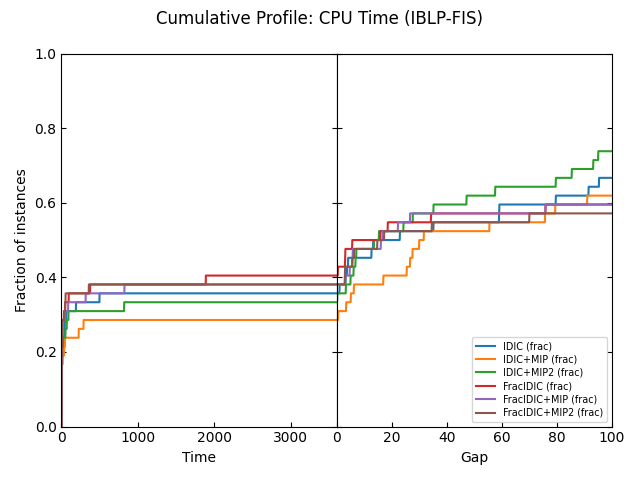}
\caption{Cumulative profile of CPU time
  \label{fig:iblpFis-cum-cpu}}
\end{subfigure}
\begin{subfigure}{0.5\textwidth}
\centering
\includegraphics[height=2.5in]{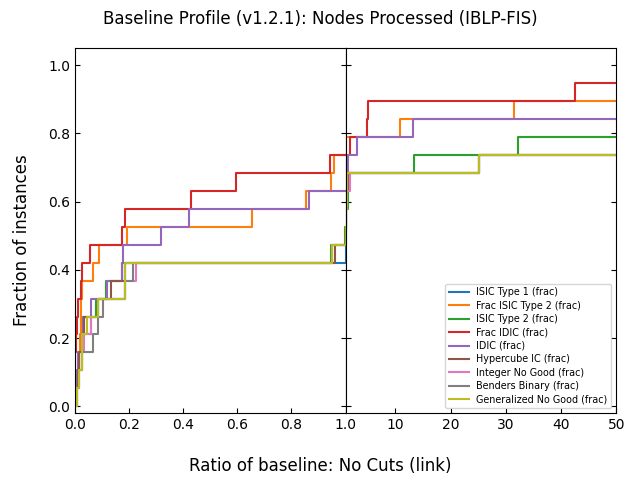}
\caption{Baseline profile of tree size 
  \label{fig:iblpFis-nodes}}
\end{subfigure}
\begin{subfigure}{0.5\textwidth}
\centering
\includegraphics[height=2.5in]{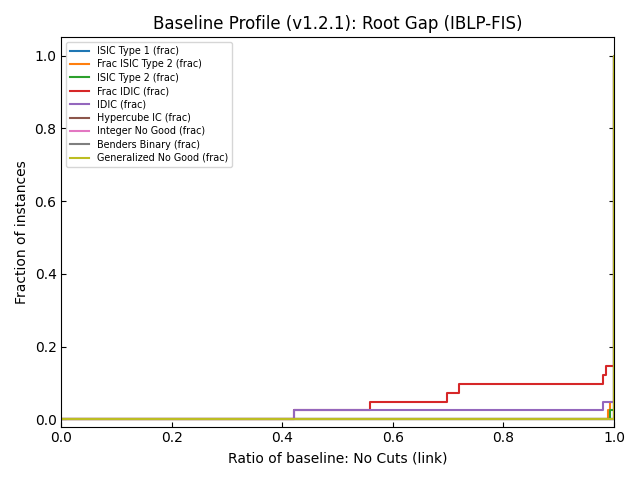}
\caption{Baseline profile of root gap
  \label{fig:iblpFis-root-gap}}
\end{subfigure}
\caption{Comparing the performance of different inequalities on the \t{IBLP-FIS} set}
\label{fig:iblpFis-perf-prof}
\end{figure}

\end{document}